\documentclass{amsart}

\usepackage{amsmath,amsthm,amssymb,amsfonts,enumerate,color}
\usepackage{mathrsfs}
\usepackage{amstext,amsxtra}
\usepackage[margin=1.5in]{geometry}
\usepackage{graphicx}
\usepackage{float,appendix}

\usepackage{tikz}
\usetikzlibrary{matrix,arrows,calc,intersections,fit}
\usetikzlibrary{decorations.markings}
\usepackage{tikz-cd}
\usepackage{textcomp}

\usepackage[colorlinks,urlcolor=black,linkcolor=blue,citecolor=blue,hypertexnames=false]{hyperref}
\allowdisplaybreaks
\usepackage{pgf,tikz}
\usepackage{stmaryrd}
\usepackage{bm}

{\theoremstyle{plain}
\newtheorem{theorem}{Theorem}[section]
\newtheorem{proposition}[theorem]{Proposition}

\newtheorem{lemma}[theorem]{Lemma}
}
{\theoremstyle{definition}
\newtheorem{example}[theorem]{Example}
\newtheorem{definition}[theorem]{Definition}
\newtheorem{construction}[theorem]{Construction}}
{\theoremstyle{remark}
\newtheorem{remark}[theorem]{Remark}}

\newcommand{\mult}{\operatorname{mult}}
\newcommand{\aut}{\operatorname{Aut}}
\newcommand{\conj}{\operatorname{conj}}
\newcommand{\id}{\operatorname{id}}

\newcommand{\stab}{\operatorname{Stab}}
\newcommand{\sym}{\operatorname{Sym}}
\newcommand{\symc}{\operatorname{Symc}}

\newcommand{\nsym}{\operatorname{Nsym}}
\newcommand{\val}{\operatorname{val}}
\newcommand{\vt}{\operatorname{Vert}}
\newcommand{\eg}{\operatorname{Edge}}
\newcommand{\mk}{\operatorname{mark}}
\newcommand{\fix}{\operatorname{fix}}
\newcommand{\rf}{\operatorname{ref}}
\newcommand{\trop}{\operatorname{trop}}
\newcommand{\ffix}{\operatorname{Fix}}
\newcommand{\im}{\operatorname{Im}}
\newcommand{\len}{\operatorname{length}}

\newcommand{\sal}{\mathcal{S}}

\newcommand{\rl}{\mathcal{R}}

\newcommand{\cl}{\mathcal{C}}

\newcommand{\bl}{\mathcal{B}}

\newcommand{\fl}{\mathcal{F}}

\newcommand{\pl}{\mathcal{P}}
\newcommand{\dl}{\mathcal{D}}

\newcommand{\rb}{\mathbb{R}}

\newcommand{~}{\quad}
\newcommand{\cb}{\mathbb{C}}
\newcommand{\nb}{\mathbb{N}}

\newcommand{\pb}{\mathbb{P}}

\newcommand{\undl}{\underline}

\newcommand{\w}{\omega}

\definecolor{cardinalred}{RGB}{140,21,21}
\definecolor{coolgray}{RGB}{77,79,83}
\definecolor{black}{RGB}{0,0,0}
\definecolor{beige}{RGB}{210,194,149}
\definecolor{darkbeige}{RGB}{179,153,93}
\definecolor{darkcardinal}{RGB}{94,48,50}
\definecolor{lightcardinal}{RGB}{141,60,30}
\definecolor{darkpurple}{RGB}{83,40,79}
\definecolor{darkcyan}{RGB}{0,124,146}
\definecolor{skyblue}{RGB}{0,152,219}
\definecolor{seablue}{RGB}{10,100,180}
\definecolor{darkblue}{RGB}{20,80,150}
\definecolor{treegreen}{RGB}{0,155,118}
\definecolor{darkorange}{RGB}{168,101,12}
\definecolor{beigegray}{RGB}{95,87,79}
\definecolor{boxgray}{RGB}{238,235,233}
\definecolor{footergray}{RGB}{199,209,197}

\begin{document}

\title[The uniform asymptotics]{The uniform asymptotics for real double Hurwitz numbers with triple ramification I: the tropical correspondence}

\author{Yanqiao Ding}

\author{Kui Li}

\author{Huan Liu}
\thanks{Corresponding author: Huan Liu, Email: liuhuan@zzu.edu.cn}

\author{Dongfeng Yan}

\address{School of Mathematics and Statistics, Zhengzhou University, Zhengzhou, 450001, China}

\email{yqding@zzu.edu.cn}


\email{likui@zzu.edu.cn}


\email{liuhuan@zzu.edu.cn}


\email{yandongfeng@zzu.edu.cn}



\date{\today}

\begin{abstract}
This is the first of two papers on the uniform asymptotics for real double Hurwitz numbers with triple ramification.
Real double Hurwitz numbers with triple ramification count the number of real ramified coverings of the complex projective line $\cb\pb^1$
by real Riemann surfaces of genus $g$, where the ramification profiles over  $0$ and $\infty$ are $\lambda$ and
$\mu$ respectively, and the ramification profiles over the remaining real branch points consist of either
$(3,1,\ldots,1)$ or $(2,1,\ldots,1)$.

We apply a modified version of the tropical computation framework developed by Markwig and Rau 
for real Hurwitz numbers (Math. Z. 281 (2015), no. 1-2, 501-522) to compute 
the real double Hurwitz numbers with triple ramification. 
The new ingredient in our computation is the application of real simple resolution, 
a technique that enables us to resolve a triple branch point into a pair of simple branch points. 
Using real simple resolution, we establish a correspondence between 
real double Hurwitz numbers with triple ramification and weighted counts of tropical covers. 
This modified tropical correspondence simplifies the asymptotic analysis of real double Hurwitz numbers with triple ramification.
\end{abstract}

\maketitle


\section{Introduction}

Double Hurwitz numbers enumerate degree $d$ ramified coverings of $\cb\pb^1$
by Riemann surfaces of genus $g$ with ramification profiles $\lambda$,
$\mu$ over $0$, $\infty$ respectively, and with simple ramification
over other branch points. Here, $\lambda$ and $\mu$ are two partitions of $d$.
Hurwitz numbers exhibit profound connections with a variety of mathematical fields, including
integrable systems  \cite{okounkov-2000},
intersection numbers over the moduli space of curves \cite{elsv-2001, gjv-2005, op-2006}, 
combinatorics \cite{hurwitz-1891}, and tropical geometry \cite{bbm-2011, cjm-2010}.
The asymptotic behavior of Hurwitz numbers has attracted much attention in the study of Hurwitz theory.
When analyzing the asymptotics of double Hurwitz numbers, two crucial scales come into play: the degree and the genus.
The asymptotics of simple Hurwitz numbers as the genus tends to infinity can be derived from \cite{hurwitz-1891,dyz-2017,yang-2025},
while the large-degree asymptotics of simple Hurwitz numbers were given in \cite{dyz-2017}, 
and those of double Hurwitz numbers have been established in \cite{li-2024}.

Real double Hurwitz numbers are the real counterparts of double Hurwitz numbers.
In comparison to complex double Hurwitz numbers, the computations and properties of 
real double Hurwitz numbers remain relatively obscure.
One reason is the extreme difficulty in computing real double Hurwitz numbers, see \cite{cadoret-2005} for example. 
Another reason lies in their dependence on the distribution of real branch points \cite{gpmr-2015}. Consequently,
it is important to find invariant lower bounds for real double Hurwitz numbers 
and to verify the effectiveness of these lower bounds by comparing their asymptotics with the asymptotics of complex double Hurwitz numbers.
It turns out that the signed count is an efficient method
for defining real invariants \cite{wel2005a,wel2005b,iks2014rw,iks2014,shustin2015}.
The signed count has also been applied to the study of real Hurwitz numbers.
Itenberg and Zvonkine introduced a signed count for normalized real polynomials \cite{iz-2018},
a framework that was later generalized to the problem of counting real simple rational functions $\frac{f(x)}{x-p}$ in \cite{er-2019},
where $f(x)\in\rb[x]$ and $p\in\rb$.

Tropical geometry provides a powerful tool to study real enumerative problems \cite{mikhalkin-2005}.
The tropical computation methods for real double Hurwitz numbers were established in \cite{gpmr-2015,mr-2015}. 
Using the tropical computation framework for real Hurwitz numbers introduced in \cite{mr-2015}, 
Rau \cite{rau2019} obtained invariant lower bounds for real double Hurwitz numbers under a parity condition on the partitions $\lambda,\mu$, 
and proved that real and complex double Hurwitz numbers are logarithmically equivalent as the degree tends to infinity.
The first author of this paper \cite{d-2020} extended Rau's result to double Hurwitz numbers 
where the partitions $\lambda,\mu$ do not satisfy Rau's condition. 
Based on the tropical computation techniques in \cite{mr-2015,rau2019}, 
the large-degree logarithmic asymptotics of real monotone double Hurwitz numbers were studied in \cite{dh-2022}.
All the above comparisons between the asymptotics of real and complex Hurwitz numbers focus on large-degree asymptotics. 
Comparisons based on large-genus asymptotics remain unexplored. 
Moreover, results concerning uniform asymptotics--\textit{i.e.}, 
the asymptotics of double Hurwitz numbers as both the degree and the genus tend to infinity--for 
real and complex double Hurwitz numbers are still lacking.

The primary goal of this paper is to establish a correspondence theorem between real double Hurwitz numbers with triple ramification 
and weighted counts of tropical covers, a result that simplifies the asymptotic analysis of real double Hurwitz numbers with triple ramification.
In a subsequent paper \cite{dlly-26}, we establish the uniform asymptotics for 
real double Hurwitz numbers with triple ramification using the tropical correspondence theorem established in this paper. 
Our uniform asymptotics for real double Hurwitz numbers automatically yields 
a lower bound for the uniform asymptotics of complex double Hurwitz numbers.
A {\it real Riemann surface} $(\Sigma,\tau)$ is a Riemann surface $\Sigma$
endowed with an anti-holomorphic involution $\tau$.
A ramified covering $\pi: (\Sigma,\tau)\to (\cb\pb^1,\conj)$ is {\it real},
if $\pi$ is compatible with the anti-holomorphic involutions, {\it i.e.}
$\pi\circ\tau=\conj\circ\pi$, where $\conj$ is the standard complex conjugation.
A branch point with ramification profile $(3,1,\ldots,1)$ of $\pi$ is called a branch point with {\it triple ramification} or a {\it triple branch point}.
Real double Hurwitz numbers (resp. real double Hurwitz numbers with triple ramification) 
count the number of degree $d$ real ramified coverings
of $\cb\pb^1$ by genus $g$ real Riemann surfaces with ramification profiles $\lambda$,
$\mu$ over $0$, $\infty$ respectively, and with simple ramification 
(resp. with triple or simple ramification) over other branch points.

Related studies of double Hurwitz numbers with triple ramification 
have primarily been conducted within the realm of combinatorics,
that is, to count
the number of minimal transitive factorizations
of a permutation $\sigma$ into $3$-cycles and $2$-cycles.
It is known \cite{dl-2015,irving-2009} that the
number of minimal transitive factorizations of a permutation
of type $(d)$ into $3$-cycles is $d^{\frac{d-3}{2}}$.
Goulden and Jackson conjectured a differential equation for 
the generating function of the numbers of minimal transitive
factorizations of a permutation \cite{gj-2000}.
To the best of our knowledge, a general answer to these related problems
in combinatorics is still unclear.
When $k$-cycles are generalized to completed $k$-cycles, the enumerative problems in question correspond
to double Hurwitz numbers with completed cycles \cite{op-2006}.
Double Hurwitz numbers with completed cycles were studied by
Shadrin, Spitz, and Zvonkine via semi-infinite wedge formalism \cite{ssz-2012},
and they established the strong piecewise polynomiality
and wall-crossing formula for these Hurwitz numbers.

In the paper \cite{mr-2015}, Markwig and Rau gave a tropical computation of real Hurwitz numbers.
In particular, they expressed real double Hurwitz numbers as
a weighted count of coloured tropical covers. Moreover, they distinguished
positive and negative real branch points by applying two different colouring rules of even edges (see also  \cite{gpmr-2015}).
Rau \cite{rau2019} constructed a family of tropical covers where the even edges can be coloured according to any rule.
Such tropical covers are called {\it zigzag covers}. The number of zigzag covers provides a lower bound for
real double Hurwitz numbers.
However, this idea is not enough for us to establish a lower bound for
real double Hurwitz numbers with triple ramification.

In order to obtain a lower bound for real double Hurwitz numbers with triple ramification,
we have to take into account two factors: the signs of the real branch points and the arrangement of triple and simple branch points.
According to Markwig and Rau's tropical computation framework \cite{mr-2015}, real double Hurwitz numbers with triple ramification
should be expressed as a weighted count of real tropical covers from real tropical curves $C$ to a tropical line $D$.
In contrast to the tropical computation for real double Hurwitz numbers,
tropical curves $C$ in this new situation may have $4$-valent vertices and genus one $2$-valent vertices.
The presence of $4$-valent vertices and genus one $2$-valent vertices complicates the characterization of the sign of a real branch point,
necessitating the use of edge colourings of even edges and edge weights of edges to distinguish between positive and negative real branch point.
Consequently, it becomes extremely challenging to analyze the uniform asymptotics for real double Hurwitz numbers with triple ramification.

To address the issue of $4$-valent vertices and genus one $2$-valent vertices in $C$,
we employ a technique called {\it the simple resolution}, which involves resolving these vertices into pairs of $3$-valent vertices
(see Definition \ref{def:rs-resolution} and Construction \ref{def:simple-resolution}). 
After applying the simple resolution and Markwig--Rau's theorem \cite[Theorem 4.1]{mr-2015}, 
real double Hurwitz numbers with triple ramification
are expressed as a weighted count of coloured tropical covers, 
with positive and negative real branch points characterized by two colouring rules for even edges (see Theorem \ref{thm:mr} for more details).

\subsection*{Organization of the paper}
In Section \ref{sec:2}, we review the definition of real double Hurwitz numbers with triple ramification
and some basic facts about tropical real Hurwitz numbers.
The next section is dedicated to introducing the real simple resolution.
Section \ref{sec:4-1} contains the calculation of local Hurwitz numbers. 
The tropical correspondence between real double Hurwitz numbers with
triple ramification and weighted counts of tropical covers is established in
Section \ref{sec:4} by applying Markwig--Rau's theorem and the real simple resolution.
The proofs of several lemmas are provided in Appendix \ref{sec:a},
and we present illustative figures referenced in the preceding sections in Appendix \ref{sec:b}.

\section{Real Hurwitz numbers}
\label{sec:2}

\subsection{Real double Hurwitz numbers with triple ramification}
In this section, we first introduce the definition of real double Hurwitz numbers with triple ramification.
Then we recall the refined real Hurwitz number introduced by Markwig and Rau in \cite{mr-2015}.

Fix two integers $d\geq1$, $g\geq0$, and let $\lambda$, $\mu$ be two partitions of $d$.
The {\it length} $l(\lambda)$ of a partition $\lambda$ is the number of entries in $\lambda$.
We choose two non-negative integers $s,t$ such that $2s+t=l(\lambda)+l(\mu)+2g-2$.
This equality is determined by the Riemann-Hurwitz formula.

Let $\undl b=\{b_1,\ldots,b_{s+t}\}\subset\rb\pb^1\setminus\{0,\infty\}$ be a set
whose elements satisfy $b_1<b_2<\cdots<b_{s+t}$.
A map $r:\undl b\to\{2,3\}$ is called {\it an $(s,t)$-ramification function} on $\undl b$
if $|r^{-1}(3)|=s$ and $|r^{-1}(2)|=t$. A set $\undl b$ together with a ramification function $r$
on it are called an {\it $(s,t)$-branch set}.
An \textit{$(s,t)$-tuple} $\Lambda_{s,t}$ is a collection $(\Lambda_1,\ldots,\Lambda_{s+t})$
of partitions of $d$ such that $s$ entries in $\Lambda_{s,t}$ are the partition $(3,1,\ldots,1)$,
and the remaining $t$ entries are $(2,1,\ldots,1)$.

%
%
%

\begin{definition}
\label{def:compatible}
Let $\undl b$ be an $(s,t)$-branch set whose elements satisfy $b_1<\cdots<b_{s+t}$,
and $\Lambda_{s,t}=(\Lambda_1,\ldots,\Lambda_{s+t})$ be an $(s,t)$-tuple.
The $(s,t)$-tuple $\Lambda_{s,t}$ is \textit{compatible} with the $(s,t)$-branch set $\undl b$,
if the following relation holds:
\begin{itemize}
    \item $\Lambda_i=(3,1,\ldots,1)$ iff $r(b_i)=3$, and $\Lambda_j=(2,1,\ldots,1)$ iff $r(b_j)=2$.
\end{itemize}
Let $\undl b$ be an $(s,t)$-branch set compatible with an $(s,t)$-tuple $\Lambda_{s,t}$.
Let $\Lambda^-_{s,t}=(\Lambda_1,\ldots,\Lambda_{|\undl b\cap\rb^-|})$
and $\Lambda^+_{s,t}=(\Lambda_{1+|\undl b\cap\rb^-|},\ldots,\Lambda_{s+t})$.
The signed splitting $\Lambda_{s,t}=(\Lambda^-_{s,t}, \Lambda^+_{s,t})$ is called a {\it compatible signed splitting}
of the $(s,t)$-branch set $\undl b$.
\end{definition}


\begin{definition}
\label{def:RDH}
Let $\undl b$ be an $(s,t)$-branch set, and $\Lambda_{s,t}=(\Lambda^-_{s,t}, \Lambda^+_{s,t})$ a compatible signed splitting of $\undl b$.
The {\it real double Hurwitz number with triple ramification} is defined as
$$
H^\rb_g(\lambda,\mu;\Lambda^-_{s,t},\Lambda^+_{s,t})=\sum_{[(\pi,\tau)]}
\frac{1}{|\aut^\rb(\pi,\tau)|}.
$$
Here, the sum is taken over all isomorphism classes $[(\pi,\tau)]$ of
ramified coverings $\pi:C\to\cb\pb^1$ satisfying the following conditions:
\begin{enumerate}
    \item $C$ is a connected Riemann surface of genus $g$;
    \item The covering $\pi:C\to\cb\pb^1$ is a degree $d$ holomorphic map;
    \item $\tau:C\to C$ is an anti-holomorphic involution of $C$
    such that $\pi\circ\tau=\conj\circ\pi$, where $\conj$ is the standard complex conjugation;
    \item $\pi$ ramifies with profiles $\lambda$ and $\mu$ over $0$ and $\infty$, respectively;
    \item $\pi$ ramifies with profile $(3,1,\ldots,1)$ over points in $r^{-1}(3)$ and ramifies with simple ramification over points in $r^{-1}(2)$;
    \item $\pi$ is unramified everywhere else.
\end{enumerate}
Note that an isomorphism of two real coverings
$(\pi_1:C_1\to\cb\pb^1,\tau_1)$
and $(\pi_2:C_2\to\cb\pb^1,\tau_2)$
is an isomorphism $\varphi:C_1\to C_2$ of Riemann surfaces
such that $\pi_1=\pi_2\circ\varphi$ and $\varphi\circ\tau_1=\tau_2\circ\varphi$,
and $\aut^{\rb}(\pi,\tau)$ is the automorphism group of  the real ramified covering $(\pi,\tau)$.
\end{definition}

Unlike its complex counterpart, the real double Hurwitz number with triple ramification
$H^\rb_g(\lambda,\mu; \Lambda^-_{s,t},\Lambda^+_{s,t})$ depends on
the positions of the branch points in $\undl b$. More exactly,
when the two partitions $\lambda$, $\mu$ are fixed, $H^\rb_g(\lambda,\mu; \Lambda^-_{s,t},\Lambda^+_{s,t})$
depends on the two sequences of partitions $\Lambda^-_{s,t}$ and $\Lambda^+_{s,t}$ \cite[Section 2]{mr-2015}.
The number $H^\rb_g(\lambda,\mu;\Lambda^-_{s,t},\Lambda^+_{s,t})$
is a generalization of the ordinary real double Hurwitz number.
In particular, when $s=0$, the number $H^\rb_g(\lambda,\mu;\Lambda^-_{0,t},\Lambda^+_{0,t})$ is the
ordinary real double Hurwitz number which depends on
the number of positive branch points $|\undl b\cap\rb^{+}|$ \cite{gpmr-2015,mr-2015}.
In the case that $s=0$ and $|\undl b\cap\rb^{+}|=r\leq t$,
we use $H^\rb_g(\lambda,\mu;r)$ to denote 
$H^\rb_g(\lambda,\mu;\Lambda^-_{0,t},\Lambda^+_{0,t})$.

Hurwitz numbers are equivalently defined as weighted counts of
unramified coverings of topological surfaces with boundaries in \cite{bbm-2011,mr-2015}.
Roughly speaking, after forgetting the complex structure and removing open discs
centered at each branch point as well as their preimages, ramified coverings of $\cb\pb^1$
are transferred to unramified coverings of topological surface with boundaries.
Note that the ramification profile of a point is translated to the degree of the
map between the corresponding boundary circles.
To get a holomorphic ramified coverings of Riemann surface from an unramified coverings
of topological surface with boundaries, one needs to glue discs with a marked point to the boundary circles
of the topological surfaces first, then choose a complex structure on the target surface
which compatible with the orientation.
Once a complex structure on the target surface is chosen, there is a unique complex structure
on the source surface such that the covering map is a holomorphic map.
It turns out that considering open unramified coverings is very convenient when we use
tropical geometry to calculate Hurwitz numbers \cite{bbm-2011,mr-2015}.

We recall two specific notations introduced in \cite{mr-2015}.
Definition \ref{def:mark-H} is an open version of \cite[Definition 2.1]{mr-2015}.

\begin{definition}
\label{def:mark-H}
Let $\sal$ be a connected oriented topological surface of genus $0$ with $n$ boundary circles $B_1,\ldots,B_n$.
Let $\mu=(\mu_1,\ldots,\mu_n)$ be a collection of $n$ partitions of $d$.
The {\it marked Hurwitz number} is the weighted sum
$$
H^{\mk,\cb}_g(\mu)=\sum_{[\pi]}\frac{1}{|\aut^{\mk}(\pi)|},
$$
where we sum over isomorphism classes of unramified marked coverings $\pi:\cl\to\sal$ satisfying the following conditions.
\begin{itemize}
    \item $\cl$ is a connected oriented topological surface of genus $g$ with $\sum_{i=1}^nl(\mu_i)$ boundary circles
$\{B_{1,1},\ldots,B_{1,l(\mu_1)},\ldots,B_{n,1},\ldots,B_{n,l(\mu_n)}\}$;
    \item $\pi$ is a degree $d$ continuous unramified covering map;
    \item Suppose that $\mu_i=(\mu_{i,1},\ldots,\mu_{i,l(\mu_i)})$, where $i\in\{1,\ldots,n\}$.
The $l(\mu_i)$ boundary circles $B_{i,1},\ldots,B_{i,l(\mu_i)}$ in $\cl$ are labelled
    such that $\pi|_{B_{i,j}}:B_{i,j}\to B_i$ is a degree $\mu_{i,j}$ map.
\end{itemize}
Here, an isomorphism of unramified marked coverings $\pi_1:\cl_1\to\sal$ and
$\pi_2:\cl_2\to\sal$ is a homeomorphism $\varphi:\cl_1\to\cl_2$ such that
$\pi_1=\pi_2\circ\varphi$ and $\varphi$ respects the labels of the boundary circles,
and $\aut^{\mk}(\pi)$ is the automorphism group of the unramified marked covering $\pi$.
\end{definition}

Let $\pi:\cl\to\sal$ be a unramified covering of surfaces with boundaries.
Suppose that the surfaces $\cl$ and $\sal$ are equipped with
orientation-reversing involutions $\iota_\cl$ and $\iota_\sal$, respectively,
such that $\pi\circ\iota_\cl=\iota_\sal\circ\pi$.
Assume that the surface $\sal$ is a genus $0$ surface with $n$ boundary circles $\bl=\{B_1,\ldots,B_n\}$
such that $\iota_{\sal}(\bl)=\bl$.
Let $B_{i,1},\ldots,B_{i,l(\mu_i)}$ be the labelled boundary circles in $\cl$ which are mapped to $B_i$.
The involutions $\iota_\cl$ and $\iota_\sal$ induce a map $\iota$ between indices of the boundary circles
such that $\iota_\sal(B_i)=B_{\iota(i)}$ and $\iota_\cl(B_{i,j})=B_{\iota(i,j)}$.
If $B_{i,j}$ is an $\iota_\cl$-invariant boundary circle ({\it i.e.} $\iota(i,j)=(i,j)$),
there are two $\iota_\cl$-fixed points in $B_{i,j}$.
Additionally, if the degree $\mu_{i,j}$ of the map $\pi|_{B_{i,j}}$ is even,
the two fixed points in $B_{i,j}$ are mapped to a same $\iota_{\sal}$-fixed points of $B_i$ (c.f. \cite[Section 2]{mr-2015}).
This gives a map $F:\{(i,j)\in\ffix(\iota)|\mu_{i,j}\equiv 0 \mod 2\}\to\ffix(\iota_{\sal})\cap (B_1\cup\cdots\cup B_n)$.

\begin{definition}[{\cite[Definition 2.5]{mr-2015}}]
\label{def:refined-real-H}
Fix two integers $g\geq0$, $d\geq1$ and a collection $\mu=(\mu_1,\ldots,\mu_{n})$ of $n$ partitions of $d$.
Let $\iota$ be an involution on the set of indices as above with $\mu_{\iota(i,j)}=\mu_{i,j}$.
Let $\sal$ be a surface of genus $0$ with $n$ boundary circles $B_1,\ldots,B_n$.
Suppose that $\iota_{\sal}$ is an orientation-reversing involution on $\sal$ such that $\iota_{\sal}(B_i)=B_{\iota(i)}$.
Fix a map
$$
F:\{(i,j)\in\ffix(\iota)|\mu_{i,j}\equiv 0 \mod 2\}\to\ffix(\iota_{\sal})\cap (B_1\cup\cdots\cup B_n)
$$
with $F(i,j)\in B_i$.
The {\it refined real Hurwitz number} $H^{\rf,\rb}_g((\sal,\iota_{\sal}),\mu,\iota,F)$
is the weighted number of marked (unramified) real coverings $\pi:\cl\to\sal$ of genus $g$
and ramification profile $\mu$ such that
\begin{itemize}
    \item $\iota_\cl(B_{i,j})=B_{\iota(i,j)}$,
    \item $\pi(\ffix(\iota_\cl\cap B_{i,j}))=F(i,j)$, if $\mu_{i,j}$ is even.
\end{itemize}
In particular, every covering map is weighted by one over the number of automorphisms of marked real coverings.
\end{definition}

\subsection{Preliminary on tropical real Hurwitz numbers}
In this section, we briefly review the correspondence theorem between real Hurwitz numbers
and tropical real Hurwitz numbers established by Markwig and Rau in \cite{mr-2015}.
We follow the notations in \cite{mr-2015} closely, and summarize a special form of \cite[Section 3]{mr-2015}.

A {\it tropical curve} $C$ is a connected metric graph with only finitely many points
which are not $2$-valent,
together with a genus function $g:C\to\nb$ which is non-zero at finitely many points
and is zero at all $1$-valent points.
An {\it inner vertex} $v$ of $C$ is a point in $C$ with valence $>2$ or $g(v)>0$.
The $1$-valent points in $C$ are called {\it leaves}.
Denote by $\vt(C)$ the set of inner vertices and leaves of $C$.
An edge $e$ in $C$ with finite length $\ell(e)\in\rb$ is called an {\it inner edge}.
Edges adjacent to leaves are called {\it ends} whose lengths are all infinite.
The set of edges of $C$ is denoted by $\eg(C)$.
An {\it open tropical curve} is a variant of the definition of tropical curve by
requiring all ends to have a finite length.
An isomorphism $\varPhi:C_1\to C_2$ of two tropical curves $C_1$, $C_2$ is
an isometric homeomorphism $\varPhi:C_1^\circ\to C_2^\circ$ respecting the genus function,
where $C_1^\circ$ and $C_2^\circ$ are sub-graphs of
$C_1$ and $C_2$ obtained by removing the $1$-valent vertices, respectively.
Denote by $\val(v)$ the valence of a vertex $v\in C$.
The sum of the first Betti number $b_1(C)$ and $\sum_{v\in C}g(v)$ is called the {\it genus} of $C$.
A automorphism $\iota: C\to C$ with $\iota^2=\id$ is called a {\it prereal structure} on $C$.
We denote by $\vt_{\text{fix}}(C)$ the set of vertices $v\in C$ with $\iota(v)=v$.
The set of edges $e$ with $\iota|_e=\id_e$ is denoted by $\eg_{\text{fix}}(C)$.

\begin{definition}[{\cite[Definition 3.2]{mr-2015}}]
\label{Rtrop-curve}
Let $C$ be a tropical curve with prereal structure $\iota$.
A {\it real structure} on $C$ consists of the following data.
\begin{enumerate}
    \item Every vertex $v\in\vt_{\text{fix}}(C)$ is associated with 
a genus $g(v)$ oriented topological surface $S_v:=S_{g(v),\val(v)}$ with $\val(v)$ boundary circles.
Moreover, the surface $S_{g(v),\val(v)}$ is equipped with
    \begin{itemize}
        \item a labelling of the boundary circles by flags adjacent to the vertex $v$
        (Denote by $B_f$ the boundary circle corresponding to the flag $f$), and
        \item an orientation-reversing involution $\iota_v$ on $S_v$ such that $\iota_v(B_f)=B_{\iota(f)}$ for all flags.
    \end{itemize}
    \item Every edge $e\in\eg_{\text{fix}}(C)$, which is given by two flags $f$, $f'$ adjacent to vertices $v$ and $v'$,
is associated with an identification
between the two fixed points of $\iota_v|_{B_f}$ and the two fixed points of $\iota_{v'}|_{B_{f'}}$.
Then {\it fixed points} of $e$ mean the two identified fixed points of $\iota_v|_{B_f}$ and $\iota_{v'}|_{B_{f'}}$.
We denote by $F_e$ the set of fixed points of $e$.
\end{enumerate}
A {\it real tropical curve} $(C,\iota_C, (S_v)_{v\in\vt(C)})$ is a prereal tropical curve $(C,\iota)$
together with a real structure on it.
\end{definition}

From \cite[Remark 3.3]{mr-2015}, a topological surface with orientation-reversing involution
can be constructed from a real tropical curve $(C,\iota_C, (S_v)_{v\in\vt(C)})$.

\begin{definition}[\cite{mr-2015}]
\label{def:tropical-morphism}
A \textit{morphism} $\varphi:C\to D$ between (open) tropical curves is a continuous
map satisfying the following conditions.
\begin{enumerate}
    \item $\varphi$ is surjective.
    \item $\varphi$ is piecewise integer affine linear, and the slope of $\varphi$ on
    a flag or edge $e$ is called the {\it weight} $\omega(e)\in\nb_{>0}$.
    \item For any point $v\in C$, we choose a
    flag $f'\subset D$ adjacent to $\varphi(v)$.
    Then the integer
    $$
    \deg(\varphi,v):=\sum_{
    \substack{f\text{ flag of } C\\
    v\in f,\varphi(f)\subset f'}}
    \omega(f)
    $$
    does not depend on the choice of $f'$.
    This is called the \textit{balancing condition} or {\it harmonicity}, and $\deg(\varphi,v)$ is called the
    {\it local degree} of $\varphi$ at $v$.
    \item $\varphi$ satisfies the local Riemann-Hurwitz condition: for each vertex $v\in C$, we have
    $$
    2-2g(v)=\deg(\varphi,v)(2-2g(\varphi(v)))-\sum_{
    \substack{e\text{ flag of } C\\
    v\in e}}(\omega(e)-1).
    $$
\end{enumerate}
\end{definition}
Let $\varphi:C\to D$ be a morphism. Edges in $C$ with even and odd weights are called {\it even edges}
and {\it odd edges}, respectively.

\begin{definition}[{\cite[Definition 3.4]{mr-2015}}]
\label{def:local-complex}
Let $D$ be a genus zero tropical curve.
Let $\varphi:C\to D$ be a morphism between two tropical curves $C$ and $D$,
and $v\in C$ be a vertex in $C$. Cut the edges in $C$ adjacent to $v$ and also the edges
in $D$ adjacent to $\varphi(v)$, then one gets a local morphism of open tropical curves
$\varphi_v:C_v\to D_{\varphi(v)}$, where $C_v$ (resp. $D_{\varphi(v)}$)
is the link of $C$ at $v$ (resp. the link of $D$ at $\varphi(v)$), and $\varphi_v=\varphi|_{C_{v}}$.
The data encoded in the local morphism $\varphi_v$ enables us to define a
complex Hurwitz number: for each flag $f$ of $\varphi(v)$, a partition of $\deg(\varphi,v)$
is obtained by collecting the weights of the flags of $v$ mapping to $f$.
The {\it local Hurwitz number} is defined as $H^\cb(\varphi,v):=H^{\mk,\cb}_{g(v)}(\mu)$,
where $\mu$ is the collection of partitions of $\deg(\varphi,v)$ indexed by the flags of $\varphi(v)$.
\end{definition}

\begin{definition}[{\cite[Definition 3.5]{mr-2015}}]
\label{def:real-cover}
Let $D$ be a genus $0$ real tropical curve with prereal structure $\iota_D$.
A {\it real tropical cover} of $D$ is a prereal curve $(C,\iota_C)$
and a morphism $\varphi:C\to D$ satisfying $\iota_D\circ\varphi=\varphi\circ\iota_C$
together with a map
$$
F:\{e\in\eg_{\fix}(C)| \omega(e)\equiv0 \mod 2\}\to\bigcup_{e'\in\eg_{\fix}(D)}F_{e'}
$$
such that $F(e)\in F_{\varphi(e)}$.
\end{definition}
An {\it isomorphism} of real covers $\varphi$, $\varphi'$ is an isomorphism
$\psi:C\to C'$ such that $\psi\circ\iota_C=\iota_{C'}\circ\psi$, $\varphi=\varphi'\circ\psi$
and $F(e)=F'(\psi(e))$ for even edges $e$.

\begin{definition}[{\cite[Definition 3.6]{mr-2015}}]
\label{def:local-real}
Let $(D,\iota_D)$ be a genus $0$ real tropical curve,
and $\varphi:C\to D$ be a real tropical cover with $v\in\vt_{\fix}(C)$.
By cutting the edges in $C$ adjacent to $v$ and edges in $D$
adjacent to $\varphi(v)$, one obtains a local real cover of open tropical curves
$\varphi_v:C_v\to D_{\varphi(v)}$, where $C_v$
is the link of $C$ at $v$, and $\varphi_v=\varphi|_{C_{v}}$.
The data encoded in the local real cover $\varphi_v$ enables us to define
the {\it local Hurwitz number} $H^\rb(\varphi,v)$ as follows:
$H^\rb(\varphi,v):=H^{\rf,\rb}_{g(v)}((\sal_{\varphi(v)}, \iota_{\varphi(v)}),\mu,\iota,F)$,
where $\mu$ is the collection of partitions of $\deg(\varphi,v)$ indexed by the flags of $\varphi(v)$,
and $\iota$ is the map between indices ineduced by $\iota_{\varphi(v)}$ and $\iota_v$.
\end{definition}
Let $(D,\iota_D)$ be a genus $0$ real tropical curve.
For a real tropical cover $\varphi:C\to D$, the {\it real multiplicity} of $\varphi$ is
\begin{equation}
\label{eq:real-mult1}
\mult^\rb(\varphi,F)=\frac{2^{|E(\varphi)|}}{|\aut(\varphi)|}
\prod_{v\in\vt_{\fix}(C)}H^{\rb}(\varphi,v)\prod_{(e,e')\in\eg_{\conj}(C)}\omega_\varphi(e)
\prod_{(v,v')\in\vt_{\conj}(C)}H^\cb(\varphi,v),
\end{equation}
where $E(\varphi)$ is the set of even edges in $\eg_{\fix}(C)$, $\eg_{\conj}(C)$
is the set of unordered $\iota_C$-conjugated pairs of edges of $C$,
and $\vt_{\conj}(C)$ is the set of unordered $\iota_C$-conjugated pairs of vertices of $C$.
\begin{definition}[{\cite[Definition 3.8]{mr-2015}}]
\label{def:real-trop-Hurwitz}
Let $g\geq0$, $d\geq1$ be two integers, and $\mu=(\mu_1,\ldots,\mu_n)$ be a collection of $n$ partitions of $d$.
Let $(D,\iota_D,(S_v)_{v\in\vt(D)})$ be a genus $0$ real tropical curve with $n$ leaves $l_1,\ldots,l_n$.
The {\it real tropical Hurwitz number} is defined as
$$
H_g^{\trop}((D,\iota_D,(S_v)_{v\in\vt(D)}),\mu)
=\sum_{[(\varphi,F)]}\mult^\rb(\varphi,F),
$$
where the sum is taken over the isomorphism classes of real tropical covers $\varphi:C\to D$
of genus $g$, degree $d$ and with ramification $\mu$
(entries in $\mu_i$ are the weights of leaves of $C$ mapping to $l_i$).
\end{definition}

Now we are ready to introduce the correspondence theorem of Markwig and Rau.
\begin{theorem}[{\cite[Theorem 4.1]{mr-2015}}]
\label{thm:markwig-rau}
Fix a genus $0$ real tropical curve $(D,\iota_D,(S_v)_{v\in\vt(D)})$.
Let $(\dl,\iota_\dl,\pl)$ be the associated genus $0$ real topological surface with punctures
according to \cite[Remark $3.3$]{mr-2015}.
Given an integer $d\geq1$, we fix a partition $\mu_i$ of $d$ for each leaf $l_i$ of $D$.
Then for any integer $g\geq0$, we have
$$
H^\rb_g((\dl,\iota_\dl),\pl,\mu)=H^{\trop}_g((D,\iota_D,(S_v)_{v\in\vt(D)}),\mu),
$$
where $\mu$ is the collection of partitions $\mu_i$.
\end{theorem}

\begin{remark}
Markwig--Rau's theorem \cite[Theorem 4.1]{mr-2015} states a correspondence theorem for
arbitrary genus real tropical curve $D$. Since we only care about the tropical computation
of real double Hurwitz numbers with triple ramification, we only need the genus $0$ case.
\end{remark}

\section{Real simple resolution}
\label{sec:3}
In this section, we recall a well-known representation of cycles in the symmetric group $S_d$,
then we use it to introduce the real simple resolution.

\subsection{Circle representation of real cycles}
\label{subsec:circle-rep}

We introduce some notations first.
A $k$-cycle $\sigma\in S_d$ is a permutation $(a_1a_2\ldots a_k)$ of type $(k,1,\ldots,1)$.
The number $k$ is called the length of the $k$-cycle $\sigma$ and is denoted by $\len(\sigma)$.
When $k$ is odd (resp. even), the $k$-cycle $\sigma$ is called an \textit{odd cycle} (resp. \textit{even cycle}).
The $k$-cycle $\sigma=(a_1a_2\ldots a_k)$ can be represented by
an oriented circle labelled with entries $a_1$, $a_2$, $\ldots$, $a_k$.
Note that the entries $a_1$, $a_2$, $\ldots$, $a_k$ are
arranged in the circle such that the orientation of
the circle induces the permutation $\sigma$.
We use the same symbol $\sigma$ to denote the circle representing the cycle $\sigma\in S_d$.
Let $\sigma_1$, $\sigma_2$ be two cycles in $S_d$.
The notation $a\in\sigma_1$ means that $\sigma_1(a)\neq a$, \textit{i.e.}, $a$ is in the support of $\sigma_1$.
$a\in\sigma_1\setminus\sigma_2$ means that $a$ is in the support of $\sigma_1$
but not in the support of $\sigma_2$.
A \textit{blank position} in a circle $\sigma$ is a position in the circle
$\sigma$ which is not located by an entry.
The orientation of the circle $\sigma$ induces an order on
the blank positions and entries of the circle.
A \textit{position before the entry $a_i$} of the circle representing $\sigma=(a_1a_2\ldots a_k)$
is the blank position locating between $a_{i-1}$ and $a_i$ if $1<i\leq k$,
and is the blank position locating between $a_{k}$ and $a_1$ if $i=1$.
If $\sigma$ is an arbitrary permutation in $S_d$, $\sigma$ possesses a disjoint cycle decomposition
$\sigma=\bar\sigma_k\circ\cdots\circ\bar\sigma_1$, unique up to a reordering of the cycles.
The circle representing $\sigma$ consists of circles representing the disjoint cycles $\bar\sigma_1,\ldots,\bar\sigma_k$.
\begin{example}\label{exa:circle1}
The circle representing $\sigma=(12345)$ is depicted in Figure \ref{fig:circle1}.
In the circle in Figure \ref{fig:circle1}, there are five blank positions.
\begin{figure}[ht]
    \centering
    \begin{tikzpicture}
    \draw[decoration={markings, mark=at position 0.18 with {\arrow{<}}},
        postaction={decorate}
        ] (-2,0) circle (1);
        \draw (-2,1) node {$1$} (-2.7,0.7) node{$5$} (-1.3,0.7) node{$2$} (-1.3,-0.7) node{$3$} (-2.7,-0.7) node{$4$};
    \end{tikzpicture}
    \caption{Circle representing $\sigma$.}
    \label{fig:circle1}
\end{figure}
\end{example}

The circle representation of cycles of $S_d$
gives us an intuitive way to interpret
the cut and join operation of a permutation in the symmetric group $S_d$.
\begin{enumerate}
    \item \textbf{Cut:}  let $\sigma$ be a $k$-cycle,
    and let $\tau=(ab)$ be a transposition with $a,b\in\sigma$.
    The composition $\tau\circ\sigma$ consists of two disjoint cycles.
    $\tau\circ\sigma$ is obtained by cutting the cycle $\sigma$ as follows.
    We cut the circle $\sigma$ at the two positions before $a$ and $b$, respectively.
    Then we glue the beginning of the oriented segment containing $a$ (resp. $b$)
    to the end of itself. See Figure \ref{fig:cut-join1}$(1)$.
    \begin{figure}[ht]
    \centering
    \begin{tikzpicture}
    \draw[decoration={markings, mark=at position 0.18 with {\arrow{<}}},
        postaction={decorate},gray
        ] (-1,-4) circle (0.5);
    \draw (-1,-3.5) node {$*$} (-1.5,-3.9) node{$.$} (-0.5,-3.9) node{$.$} (-1.3,-4.4) node{$b$} (-0.7,-4.4) node{$*$};
    \draw[decoration={markings, mark=at position 0.2 with {\arrow{<}}},
        postaction={decorate},gray
        ] (-1,-6) circle (0.5);
    \draw (-1,-6.5) node {$*$} (-1.5,-5.9) node{$\vdots$} (-0.5,-5.9) node{$\vdots$} (-1.3,-5.6) node{$\circ$} (-0.7,-5.6) node{$a$};
    \draw[decoration={markings, mark=at position 0.2 with {\arrow{<}}},  postaction={decorate},gray
        ] (-2,-5) ellipse (0.5 and 1);
    \draw (-2.5,-5) node {\color{red}$-$} (-1.5,-5) node {\color{red}$-$} (-2,-4) node {$*$} (-2,-6) node {$*$} (-2.37,-4.3) node{$.$} (-1.63,-4.3) node{$.$} (-2.53,-4.8) node{$b$} (-1.53,-4.8) node{$*$} (-2.4,-5.5) node{$\vdots$} (-1.6,-5.5) node{$\vdots$} (-2.53,-5.2) node{$\circ$} (-1.53,-5.2) node{$a$} (-1.5,-7) node{\tiny $(1)$ Cut};
    \draw[decoration={markings, mark=at position 0.18 with {\arrow{<}}},
        postaction={decorate},gray
        ] (1,-4) circle (0.5);
    \draw (1,-4.5) node{\color{red}$|$} (1,-3.5) node {$*$} (0.5,-3.9) node{$.$} (1.5,-3.9) node{$.$} (0.7,-4.4) node{$b$} (1.3,-4.4) node{$*$};
    \draw[decoration={markings, mark=at position 0.2 with {\arrow{<}}},
        postaction={decorate},gray
        ] (1,-6) circle (0.5);
    \draw (1,-5.5) node{\color{red}$|$} (0.5,-5.9) node{$\vdots$} (1.5,-5.9) node{$\vdots$} (0.7,-5.6) node{$\circ$} (1.3,-5.6) node{$a$} (0.7,-6.4) node{$*$} (1.3,-6.4) node{$*$};
    \draw[decoration={markings, mark=at position 0.2 with {\arrow{<}}},  postaction={decorate},gray
        ] (2,-5) ellipse (0.5 and 1);
    \draw (2,-4) node {$*$} (1.62,-4.3) node{$.$} (2.38,-4.3) node{$.$} (1.53,-4.8) node{$b$} (2.47,-4.8) node{$*$} (1.6,-5.5) node{$\vdots$} (2.4,-5.5) node{$\vdots$} (1.53,-5.2) node{$\circ$} (2.47,-5.2) node{$a$} (1.72,-5.9) node{$*$} (2.26,-5.9) node{$*$} (1.5,-7) node{\tiny $(2)$ Join};
    \end{tikzpicture}
    \caption{Cut and join of the cycle $\sigma$.}
    \label{fig:cut-join1}
    \end{figure}
    \item \textbf{Join:} let $\sigma=\sigma_1\circ\sigma_2\in S_d$
    be the decomposition of $\sigma$ into disjoint cycles,
    and let $\tau=(ab)$ be a transposition such that
    $a\in\sigma_1$ and $b\in\sigma_2$.
    The composition $\tau\circ\sigma$ is a $k$-cycle, where $k=\len(\sigma_1)+\len(\sigma_2)$.
    $\tau\circ\sigma$ is obtained by joining the cycles $\sigma_1$ and $\sigma_2$ as follows.
    We cut the circle $\sigma_1$ at the position before $a$
    and cut the circle $\sigma_2$ at the position before $b$, respectively.
    Then we glue the beginning  (resp. end) of the oriented segment containing $a$
    to the end (resp. beginning) of the oriented segment containing $b$.
    See Figure \ref{fig:cut-join1}$(2)$.
\end{enumerate}

In the following, we consider the circle representation of real cycles.
Let $\sigma,\sigma_1,\sigma_2$ be three $k$-cycles in the symmetric group $S_d$,
and $\gamma\in S_d$ be an involution, \textit{i.e.}, $\gamma^2=\id$.
Assume that $\gamma\circ\sigma\circ\gamma=\sigma^{-1}$ and
$\gamma\circ\sigma_1\circ\gamma=\sigma_2^{-1}$.
Suppose that $\sigma=(a_1a_2\ldots a_k)$, $\sigma_1=(b_1b_2\ldots b_k)$,
$\sigma_2=(c_1c_2\ldots c_k)$.
From \cite[Lemma $3.12$]{gpmr-2015}, we obtain the following facts:
if $k$ is odd, there is exactly one entry $x\in\sigma\setminus\gamma$;
if $k$ is even, either $\sigma\setminus\gamma$ contains two entries or $\sigma\setminus\gamma$ contains no entry.
Entries in $\sigma\setminus\gamma$ are called fixed points
of the involution $\gamma$ in $\sigma$.
Moreover, it follows from \cite[Lemma $3.12$]{gpmr-2015} that
the real cycle $\sigma$ has a symmetric structure which
is interpreted by circle representation in Figure \ref{fig:sym-circle1}$(1)-(3)$,
and a conjugated pair $(\sigma_1,\sigma_2)$ is depicted in Figure \ref{fig:sym-circle1}$(4)$.
Let $\sigma,\gamma$ be two elements of $S_d$ with $\gamma^2=\id$ and $\gamma\circ\sigma\circ\gamma=\sigma^{-1}$.
Suppose that $\sigma=\bar\sigma_k\circ\cdots\circ\bar\sigma_1$ is the decomposition of $\sigma$
into disjoint cycles. From \cite[Lemma $3.12$]{gpmr-2015},
there is a symmetric axis of the circle representing $\sigma$,
and the circle representing $\sigma$ consists of circles of the type depicted in Figure \ref{fig:sym-circle1}.
A pair of symmetric blank positions of circle $\sigma$
are called a \textit{pair of conjugated positions}.

\begin{figure}[ht]
    \centering
    \begin{tikzpicture}
    \draw[decoration={markings, mark=at position 0.18 with {\arrow{<}}},
        postaction={decorate}
        ] (-2,0) circle (1);
    \draw (-2,1) node {$a_1$} (-2.7,0.7) node{$a_k$} (-1.3,0.7) node{$a_2$}
    (-1,0) node{$\vdots$} (-3,0) node{$\vdots$}
    (-1.3,-0.7) node{$a_{\frac{k+1}{2}}$} (-2.7,-0.7) node{$a_{\frac{k+3}{2}}$};
    \draw [stealth-stealth](-2.5,0.7) -- (-1.5,0.7);
    \draw [stealth-stealth](-2.3,-0.7) -- (-1.7,-0.7);
    \draw [stealth-stealth](-2.5,0) -- (-1.5,0);
    \draw[gray,dashed](-2,1.3) -- (-2,-1.3);
    \draw[decoration={markings, mark=at position 0.18 with {\arrow{<}}},
        postaction={decorate}
        ] (1,0) circle (1);
    \draw (1,1) node {$a_1$} (0.3,0.7) node{$a_k$} (1.7,0.7) node{$a_2$}
    (1.7,-0.7) node{$a_{\frac{k}{2}}$} (0,0) node{$\vdots$} (2,0) node{$\vdots$}
    (0.3,-0.7) node{$a_{\frac{k}{2}+2}$} (1,-1) node{$a_{\frac{k}{2}+1}$};
    \draw [stealth-stealth](0.5,0.7) -- (1.5,0.7);
    \draw [stealth-stealth](0.7,-0.7) -- (1.3,-0.7);
    \draw [stealth-stealth](0.5,0) -- (1.5,0);
    \draw[gray,dashed](1,1.3) -- (1,-1.3);
    \draw[decoration={markings, mark=at position 0.18 with {\arrow{<}}},
        postaction={decorate}
        ] (4,0) circle (1);
    \draw (3,0) node {$\vdots$} (3.3,0.7) node{$a_k$}
    (4.7,0.7) node{$a_1$} (4.7,-0.7) node{$a_{\frac{k}{2}}$}
    (3.3,-0.7) node{$a_{\frac{k}{2}+1}$} (5,0) node{$\vdots$};
    \draw [stealth-stealth](3.5,0.7) -- (4.5,0.7);
    \draw [stealth-stealth](3.5,0) -- (4.5,0);
    \draw [stealth-stealth](3.7,-0.7) -- (4.3,-0.7);
    \draw[gray,dashed](4,1.3) -- (4,-1.3);
    \draw[decoration={markings, mark=at position 0.2 with {\arrow{<}}},  postaction={decorate}] (6.5,0) ellipse (0.5 and 1);
    \draw (7,0) node {$\vdots$} (6.9,0.7) node{$b_1$}
    (6.9,-0.7) node{$b_k$} (7.5,0) node{$\vdots$}
    (7.7,-0.7) node{$c_k$} (7.7,0.7) node{$c_1$};
    \draw[decoration={markings, mark=at position 0.2 with {\arrow{<}}},  postaction={decorate}] (8,0) ellipse (0.5 and 1);
    \draw [stealth-stealth](7.05,0.7) -- (7.45,0.7);
    \draw [stealth-stealth](7.05,0) -- (7.45,0);
    \draw [stealth-stealth](7.05,-0.7) -- (7.45,-0.7);
    \draw[gray,dashed](7.25,1.3) -- (7.25,-1.3);
    \draw (-2,-1.8) node{\tiny $(1)$ $k$ is odd,} (-2,-2) node{\tiny $a_1$ is the fixed point}
    (1,-1.8) node{\tiny $(2)$ $k$ is even,} (1,-2.1) node{\tiny $a_1,a_{\frac{k}{2}+1}$ are fixed points}
    (4,-1.8) node{\tiny $(3)$ $k$ is even,} (4,-2) node{\tiny no fixed point} (7.25,-1.8) node{\tiny $(4)$ a conjugated pair};
    \end{tikzpicture}
    \caption{Symmetric structures of real cycles.}
    \label{fig:sym-circle1}
\end{figure}

Now we consider the circle representation of real cut and real join.
Let $\sigma,\gamma\in S_d$ satisfy $\gamma^2=\id$ and
$\gamma\circ\sigma\circ\gamma=\sigma^{-1}$.
Assume that $\tau=(ab)$ is a transposition such that
$\gamma\circ(\tau\circ\sigma)\circ\gamma=(\tau\circ\sigma)^{-1}$.
Suppose that $\sigma=\sigma_s\circ\cdots\circ\sigma_1$ is the decomposition of
$\sigma$ into disjoint cycles.
If $a,b$ are in one cycle of $\sigma_1,\cdots,\sigma_s$,
the cycle $\tau\circ\sigma$ is a \textit{real cut} of $\sigma$.
Otherwise, $\tau\circ\sigma$ is a \textit{real join} of $\sigma$.
The circle representations of real cuts are depicted in Figure \ref{fig:real-cut1}.
The real join operation is the inverse of real cut,
so the circle representations of real join are obtained by
interchanging the positions of the circles and ellipsoids in Figure \ref{fig:real-cut1}.
\begin{figure}[ht]
    \centering
    \begin{tikzpicture}
    \draw[decoration={markings, mark=at position 0.18 with {\arrow{<}}},
        postaction={decorate},gray
        ] (-1,-0) circle (0.5);
    \draw (-1,0.5) node {$*$} (-1.5,0.1) node{$.$} (-0.5,0.1) node{$.$} (-1.3,-0.4) node{$b$} (-0.7,-0.4) node{$*$};
    \draw[gray,dashed](-1,0.75) -- (-1,-0.85);
    \draw[decoration={markings, mark=at position 0.2 with {\arrow{<}}},
        postaction={decorate},gray
        ] (-1,-2) circle (0.5);
    \draw (-1.28,-2.4) node {$*$} (-0.72,-2.4) node {$*$} (-1.5,-1.9) node{$\vdots$} (-0.5,-1.9) node{$\vdots$} (-1.3,-1.6) node{$\circ$} (-0.7,-1.6) node{$a$};
    \draw[gray,dashed](-1,-1.15) -- (-1,-2.85);
    \draw[decoration={markings, mark=at position 0.2 with {\arrow{<}}},  postaction={decorate},gray
        ] (-2,-1) ellipse (0.5 and 1);
    \draw (-2.5,-1) node {\color{red}$-$} (-1.5,-1) node {\color{red}$-$} (-2,-0) node {$*$} (-1.72,-1.9) node {$*$} (-2.37,-0.3) node{$.$} (-1.63,-0.3) node{$.$} (-2.53,-0.8) node{$b$} (-1.53,-0.8) node{$*$} (-2.4,-1.5) node{$\vdots$} (-1.6,-1.5) node{$\vdots$} (-2.53,-1.2) node{$\circ$} (-1.53,-1.2) node{$a$} (-2.28,-1.9) node{$*$};
    \draw (-1.5,-3) node{\tiny$(1)$ one fixed point};
    \draw[gray,dashed](-2,0.2) -- (-2,-2.2);
    \draw[decoration={markings, mark=at position 0.18 with {\arrow{<}}},
        postaction={decorate},gray
        ] (2,0) circle (0.5);
    \draw (2,0.5) node {$*$} (1.5,0.1) node{$\vdots$} (2.5,0.1) node{$\vdots$} (1.7,-0.4) node{$b$} (2.3,-0.4) node{$*$};
    \draw[gray,dashed](2,0.85) -- (2,-0.85);
    \draw[decoration={markings, mark=at position 0.2 with {\arrow{<}}},
        postaction={decorate},gray
        ] (2,-2) circle (0.5);
    \draw (2,-2.5) node {$*$} (1.5,-1.9) node{$\vdots$} (2.5,-1.9) node{$\vdots$} (1.7,-1.6) node{$\circ$} (2.3,-1.6) node{$a$};
    \draw[gray,dashed](2,-1.15) -- (2,-2.85);
    \draw[decoration={markings, mark=at position 0.2 with {\arrow{<}}},  postaction={decorate},gray
        ] (1,-1) ellipse (0.5 and 1);
    \draw (1.5,-1) node {\color{red}$-$} (0.5,-1) node {\color{red}$-$} (1,0) node {$*$} (1,-2) node {$*$} (0.6,-0.3) node{$\vdots$} (1.4,-0.3) node{$\vdots$} (0.53,-0.8) node{$b$} (1.47,-0.8) node{$*$} (0.6,-1.5) node{$\vdots$} (1.4,-1.5) node{$\vdots$} (0.53,-1.2) node{$\circ$} (1.47,-1.2) node{$a$};
    \draw[gray,dashed](1,0.2) -- (1,-2.2);
    \draw (1.5,-3) node{\tiny$(2)$ two fixed points};
    \draw[decoration={markings, mark=at position 0.18 with {\arrow{<}}},
        postaction={decorate},gray
        ] (5,0) circle (0.5);
    \draw (4.5,0.1) node{$\vdots$} (5.5,0.1) node{$\vdots$} (4.7,-0.4) node{$b$} (5.3,-0.4) node{$*$};
    \draw[gray,dashed](5,0.85) -- (5,-0.85);
    \draw[decoration={markings, mark=at position 0.2 with {\arrow{<}}},
        postaction={decorate},gray
        ] (5,-2) circle (0.5);
    \draw (4.5,-1.9) node{$\vdots$} (5.5,-1.9) node{$\vdots$} (4.7,-1.6) node{$\circ$} (5.3,-1.6) node{$a$};
    \draw[gray,dashed](5,-1.15) -- (5,-2.85);
    \draw[decoration={markings, mark=at position 0.2 with {\arrow{<}}},  postaction={decorate},gray
        ] (4,-1) ellipse (0.5 and 1);
    \draw (4.5,-1) node {\color{red}$-$} (3.5,-1) node {\color{red}$-$} (3.6,-0.3) node{$\vdots$} (4.4,-0.3) node{$\vdots$} (3.53,-0.8) node{$b$} (4.47,-0.8) node{$*$} (3.6,-1.5) node{$\vdots$} (4.4,-1.5) node{$\vdots$} (3.53,-1.2) node{$\circ$} (4.47,-1.2) node{$a$};
    \draw[gray,dashed](4,0.2) -- (4,-2.2);
    \draw[decoration={markings, mark=at position 0.18 with {\arrow{<}}},
        postaction={decorate},gray
        ] (7.5,0) circle (0.5);
    \draw (7.5,0.5) node {$\ldots$} (8,0.1) node{$\vdots$} (7,0.1) node{$\vdots$} (7.2,-0.4) node{$b$} (7.8,-0.4) node{$*$};
    \draw[gray,dashed](7,-1) -- (8,-1);
    \draw[decoration={markings, mark=at position 0.2 with {\arrow{<}}},
        postaction={decorate},gray
        ] (7.5,-2) circle (0.5);
    \draw (7.5,-2.5) node {$\ldots$} (8,-1.9) node{$\vdots$} (7,-1.9) node{$\vdots$} (7.2,-1.6) node{$\circ$} (7.8,-1.6) node{$a$};
    \draw[decoration={markings, mark=at position 0.2 with {\arrow{<}}},  postaction={decorate},gray
        ] (6.5,-1) ellipse (0.5 and 1);
    \draw (6,-1) node {\color{red}$-$} (7,-1) node {\color{red}$-$} (6.5,0) node {$\ldots$} (6.5,-2) node {$\ldots$} (6.1,-0.3) node{$\vdots$} (6.9,-0.3) node{$\vdots$} (6.03,-0.8) node{$b$} (6.97,-0.8) node{$*$} (6.1,-1.5) node{$\vdots$} (6.9,-1.5) node{$\vdots$} (6.03,-1.2) node{$\circ$} (6.97,-1.2) node{$a$};
    \draw (5.75,-3) node{\tiny$(3)$ no fixed point};
    \draw[gray,dashed](5.8,-1) -- (7.1,-1);
    \end{tikzpicture}
    \caption{Real cut of the real cycle $\sigma$.}
    \label{fig:real-cut1}
\end{figure}

\subsection{Real simple resolution}
\label{subsec:simple}

By using monodromy representation, a real cover can be described as a tuple in the symmetric group.
Real double Hurwitz numbers were interpreted as the number of factorizations of identity into
permutations of particular types \cite{cadoret-2005,gpmr-2015,dh-2022}.

Fix two integers $g\geq0$, $d\geq3$, and two partitions $\lambda$, $\mu$
of $d$ with $l(\lambda)+l(\mu)+2g-2=2$.
We have three cases:
\begin{enumerate}[(1)]
\item $(l(\lambda),l(\mu),g)=(1,1,1)$;
\item $(l(\lambda),l(\mu),g)=(1,3,0)$ or $(3,1,0)$;
\item $(l(\lambda),l(\mu),g)=(2,2,0)$.
\end{enumerate}
Let $\rl_g(\lambda,\mu)$ be the set of isomorphism classes of real ramified coverings $(\pi:\Sigma_g\to\cb\pb^1,\iota)$ with ramification profiles $\lambda$,
$(3,1,\ldots,1)$ and $\mu$ over $0$, $p\in\rb_+$ and $\infty$, respectively.

Denote by $\fl^\rb_g(\lambda,\mu)$ the set of tuples $(\sigma_1,\gamma,\tau,\sigma_2)$
of elements of the symmetric group $S_d$ satisfying:
\begin{itemize}
    \item $\sigma_2\circ\tau\circ\sigma_1=\id$;
    \item the cycle types of $\sigma_1, \tau$ and $\sigma_2$ are $\lambda$, $(3,1,\ldots,1)$ and $\mu$, respectively;
    \item the subgroup generated by $\sigma_1,\tau,\sigma_2$ acts transitively on the set $\{1,\ldots,d\}$;
    \item $\gamma$ is an involution (i.e. $\gamma^2=\id$)
    satisfying:
    \begin{equation}\label{eq:real-factor2}
    \gamma\circ\sigma_1\circ\gamma=\sigma_1^{-1} \text{ and } \gamma\circ(\tau\circ\sigma_1)\circ\gamma=(\tau\circ\sigma_1)^{-1}.
    \end{equation}
\end{itemize}

Construction $\ref{const-1}$ is a modification of \cite[Construction 2.4]{gpmr-2015}.
\begin{figure}[ht]
    \centering
    \begin{tikzpicture}
    \draw (0,0)--(7,0);
    \foreach \Point in {(2,0),(4,0)}
    \draw[decoration={markings, mark=at position 0.125 with {\arrow{>}}},
        postaction={decorate}
        ] \Point circle (0.3);
    \foreach \Point in {(1,0), (2,0),(4,0)}
    \draw[fill=black] \Point circle (0.05);
    \draw (1.2,-0.3) node{\tiny{$x_0$}};
    \draw (7.5,0) node{$\rb$};
    \draw[bend left,-]  (1,0) to (2,0.3);
    \draw (1,0) .. controls (1.5,0.7) and (2,0.6) .. (4,0.3);
    \draw (1.5,0.2) node{\tiny{$l_1$}} (2.5,0.6) node{\tiny{$l_2$}} (2,-0.15) node{\tiny{$0$}} (4,-0.15) node{\tiny{$p$}};
    \end{tikzpicture}
    \caption{Generators of $\pi_1(\cb P^1\setminus\{0,p,\infty\},x_0)$.}
    \label{fig:real-loops1}
\end{figure}

\begin{construction}[{\cite[Construction 2.4]{gpmr-2015}}]
\label{const-1}
Choose two fixed points $x_0, p\in\rb\pb^1\setminus\{0,\infty\}\subset\cb\pb^1$ with $x_0<0<p$.
We regard the point $x_0$ as the base point of $\pi_1(\cb\pb^1\setminus\{0, p, \infty\}, x_0)$.
Let $l_1, l_2$ be two loops as depicted in Figure \ref{fig:real-loops1}.
Loops $l_1$ and $l_2$ generate the fundamental group $\pi_1(\cb\pb^1\setminus\{0, p, \infty\}, x_0)$.
The action of the standard complex conjugation on the fundamental group is characterized by the following equations:
$$
\conj(l_1)=l_1^{-1},~ \conj(l_2\cdot l_1)=(l_2\cdot l_1)^{-1}.
$$
Let $[(\pi,\iota)]\in\rl_g(\lambda,\mu)$ be a real ramified covering.
Label the preimage $\pi^{-1}(x_0)$ as $\{q_1,\ldots,q_d\}$.
Denote by $L$ this labelling.
The monodromy representation
$$
\Phi_{(\pi,L)}:\pi_1(\cb\pb^1\setminus\{0, p, \infty\}, x_0)\to S_d
$$
maps loops $l_1$, $l_2$ and $l_2\cdot l_1$ to permutations of type $\lambda$, $(3,1,\ldots,1)$ and $\mu$, respectively.
The involution $\iota$ induces an involution $\gamma\in S_d$ on the indexes of $\pi^{-1}(x_0)$.

Given a tuple $(\sigma_1,\gamma,\tau,\sigma_2)\in\fl^\rb_g(\lambda,\mu)$,
the classical Hurwitz construction
(see \cite{hurwitz-1891} or \cite[Chapter $7$]{cm-2016})
constructs a ramified covering $\pi:\Sigma_g\to\cb P^1$ with ramification
profiles $\lambda$, $(3,1,\ldots,1)$ and $\mu$ over $0$, $p$ and $\infty$, respectively.
Here, $\Sigma_g$ is a genus $g$ Riemann surface.
Moreover, points in $\pi^{-1}(x_0)$ are labelled as $\{q_1,\ldots,q_d\}$,
and the monodromy actions of the loops $l_1,l_2$ are represented by $\sigma_1,\tau$ respectively.
Real involution $\iota$ on the surface $\Sigma_g$ is obtained from $\gamma$ in \cite[Construction 2.4]{gpmr-2015}.
We summarize it as follows.
For an unramified point $y\in \Sigma_g$, we choose a path $h$ in
$\cb P^1\setminus\{0,p,\infty\}$ from $x_0$ to $\pi(y)$.
Let $\tilde h$ be the path in $\Sigma_g$ lifting $h$ with endpoint $y$.
Let $q_k$ be the starting point of $\tilde h$.
Denote by $\conj(h)\subset\cb P^1\setminus\{0,p,\infty\}$ the conjugated path of $h$.
Lift $\conj(h)$ to the path $\widetilde\conj(h)$ in $\Sigma_g$ with starting point $q_{\gamma(k)}$.
We define $\iota(y)$ to be the endpoint of $\widetilde\conj(h)$.
The fourth condition in the definition of $\fl^\rb_g(\lambda,\mu)$ implies that $\iota(y)$ is well-defined.
By standard arguments, the involution $\iota$ can be extended to all of $\Sigma_g$.
Hence, we obtain a real ramified covering $(\pi:\Sigma_g\to\cb\pb^1,\iota)$.
\end{construction}
Consider the action of $S_d$ on $\fl^\rb_g(\lambda,\mu)$ by conjugation (coordinatewise).
Denote by $\fl^\rb_g(\lambda,\mu)/S_d$ the set of orbits of this action.
The following lemma is \cite[Lemma 2.3]{gpmr-2015} essentially.
We just modify \cite[Lemma 2.3]{gpmr-2015} to our case.
\begin{lemma}[{\cite[Lemma 2.3]{gpmr-2015}}]
\label{lem:rep}
Fix two integers $g\geq0$, $d\geq3$, and two partitions $\lambda$, $\mu$
of $d$ with $l(\lambda)+l(\mu)+2g-2=2$.
There is a bijection
$$
\Psi:\fl^\rb_g(\lambda,\mu)/S_d\to\rl_g(\lambda,\mu).
$$
Moreover, $|\stab_{S_d}(A)|=|\aut^\rb(\Psi([A]))|$ for any tuple
$A=(\sigma_1,\gamma,\tau,\sigma_2)\in\fl^\rb_g(\lambda,\mu)$.
\end{lemma}

\begin{proof}
With the help of Construction \ref{const-1}, one can prove this lemma in the same way as
\cite[Lemma 2.3]{gpmr-2015}, so we omit it.
\end{proof}

Let $g,d,\lambda,\mu$ be the same as above.
Denote by $\widetilde\fl^\rb_g(\lambda,\mu)$ the set of tuples $(\sigma_1,\gamma,\tau_1,\tau_2,\sigma_2)$
of elements of the symmetric group $S_d$ satisfying:
\begin{itemize}
    \item $\sigma_2\circ\tau_2\circ\tau_1\circ\sigma_1=\id$;
    \item the cycle types of $\sigma_1$ and $\sigma_2$ are $\lambda$ and $\mu$, respectively,
    and $\tau_1,\tau_2$ are transpositions;
    \item $\tau_2\circ\tau_1$ is a $3$-cycle.
    \item the subgroup generated by $\sigma_1,\tau_1,\tau_2,\sigma_2$ acts transitively on the set $\{1,\ldots,d\}$;
    \item $\gamma$ is an involution (i.e. $\gamma^2=\id$)
    satisfying:
    \begin{equation}\label{eq:real-factor3}
    \begin{aligned}
    \gamma\circ\sigma_1\circ&\gamma=\sigma_1^{-1},\\
    \gamma\circ(\tau_1\circ\sigma_1)\circ\gamma=(\tau_1\circ\sigma_1)^{-1}&\text{ and }
    \gamma\circ(\tau_2\circ\tau_1\circ\sigma_1)\circ\gamma=(\tau_2\circ\tau_1\circ\sigma_1)^{-1}.
    \end{aligned}
    \end{equation}
\end{itemize}

Given any tuple $(\sigma_1,\gamma,\tau_1,\tau_2,\sigma_2)\in\widetilde\fl^\rb_g(\lambda,\mu)$, 
a graph called {\it monodromy graph} can be constructed from
$(\sigma_1,\gamma,\tau_1,\tau_2,\sigma_2)$ according to \cite[Construction 3.13]{gpmr-2015} or \cite[Construction 2.21]{dh-2022}.
Here, we use the colour rule in \cite[Construction 2.21]{dh-2022} to draw the monodromy graph corresponding to a tuple.
For the convenience of readers,
we briefly sketch \cite[Construction $2.21$]{dh-2022} as follows.
Let $(\sigma_1,\gamma,\tau_1,\tau_2,\sigma_2)$ be a tuple in
$\widetilde\fl^\rb_g(\lambda,\mu)$.
Real even cycles in $\sigma_1$, $\tau_1\circ\sigma_1$ and $\tau_2\circ\tau_1\circ\sigma_1$
are represented by edges in colouring,
pairs of conjugated cycles with respect to the involution $\gamma$
are represented by dotted edges,
and the remaining odd cycles are represented by edges in black.
Moreover, all edges are weighted by the length of the corresponding cycle.
Edges are intersected at $3$-valent vertices which correspond
to the transpositions $\tau_1$ and $\tau_2$.
The colouring rule for real even cycles are given as follows:
\begin{itemize}
    \item Even cycles in $\sigma_1$, $\tau_1\circ\sigma_1$ and $\tau_2\circ\tau_1\circ\sigma_1$
    having no fixed point
    with respect to the involution $\gamma$ are represented by edges in blue,
    and those having two fixed points are represented by edges in red;
\end{itemize}
For more details, we refer the readers to \cite[Construction 3.13]{gpmr-2015} or \cite[Construction 2.21]{dh-2022}.

\begin{definition}\label{def:induced-permutation}
Let $\sigma$ be a circle representing a cycle,
and $a_1,a_2,\ldots,a_s$ be $s$ entries locating in the oriented circle $\sigma$.
The permutation of the set $\{a_1,a_2,\ldots,a_s\}$ induced by the orientation
of the circle $\sigma$ is called the \textit{induced permutation} of $a_1,a_2,\ldots,a_s$ in $\sigma$.
\end{definition}

Given a $3$-cycle $\tau=(abc)\in S_d$,
there are three decompositions of $\tau$ into transpositions:
$\tau=(ac)\circ(ab), \tau=(ab)\circ(bc), \tau=(bc)\circ(ac)$.
Lemma $\ref{lem:computation-1}-\ref{lem:computation-3}$ describe the real decomposition of a real $3$-cycle.

\begin{definition}\label{def:standard-decom}
Let $\tau$ be the $3$-cycle in a tuple
$(\sigma_1,\gamma,\tau,\sigma_2)\in\fl^\rb_g(\lambda, \mu)$.
A {\it standard real decomposition} of $\tau$ is a decomposition of $\tau$ into transpositions
$\tau=\bar\tau_2\circ\bar\tau_1$ such that
\begin{itemize}
\item $\gamma\circ(\bar\tau_1\circ\sigma_1)\circ\gamma=(\bar\tau_1\circ\sigma_1)^{-1}$,
\item blank positions in the circle $\sigma_1$ before the two entries in the support
of $\bar\tau_1$ form a pair of conjugated positions,
and the blank position in the circle $\sigma_1$ before the entry in $\tau\setminus\bar\tau_1$ is
an intersection of the circle $\sigma_1$ and its symmetric axis.
\end{itemize}
\end{definition}

\begin{figure}
\centering
    \begin{tikzpicture}
 \draw[line width=0.4mm] (-0.75,-8.5)--(-1.75,-8.5);
    \draw[line width=0.4mm] (0.25,-8.5)--(1.25,-8.5);
    \draw[line width=0.4mm] (-0.25,-8.5) +(0:0.5 and 0.25) arc (0:180:0.5 and 0.25);
    \draw[line width=0.4mm,blue] (-0.25,-8.5) +(0:0.5 and 0.25) arc (0:-180:0.5 and 0.25);
    \draw (-1.25,-8.7) node {\tiny $o$} (-0.25,-8.1) node {\tiny $o_1$} (-0.25,-9) node {\tiny $e_1$} (0.65,-8.7) node {\tiny $o$};
    \draw (-0.25,-9.6) node {(1). \small $1$};
    \draw[line width=0.4mm,blue] (3.05,-8.5)--(2.05,-8.5);
    \draw[line width=0.4mm,blue] (4.05,-8.5)--(5.05,-8.5);
    \draw[line width=0.4mm,blue] (3.55,-8.5) ellipse (0.5 and 0.25);
    \draw (2.45,-8.7) node {\tiny $e$} (3.55,-8.1) node {\tiny $e_1$} (3.55,-9) node {\tiny $e_2$} (4.45,-8.7) node {\tiny $e$};
    \draw (3.55,-9.6) node {(2). \small $4(2)$};
\end{tikzpicture}
    \caption{Two monodromy graphs. When $e_1=e_2$, the number of 3-cycles is in the bracket.}
    \label{fig:1-local-brief}
\end{figure}

\begin{lemma}
\label{lem:computation-1}
\begin{enumerate}[$(1)$]
\item Let $\tau$ be the $3$-cycle in a tuple
$(\sigma_1,\gamma,\tau,\sigma_2)\in\fl^\rb_1((d), (d))$,
then $\tau$ has a unique standard real decomposition $\tau=\bar\tau_2\circ\bar\tau_1$.
\item Given a tuple $(\sigma_1,\gamma,\tau,\sigma_2)\in\fl^\rb_1((d),(d))$,
the monodromy graph corresponding to the tuple
$(\sigma_1,\gamma,\bar\tau_1,\bar\tau_2,\sigma_2)$ is one of the two types depicted in Figure $\ref{fig:1-local-brief}$,
where $\bar\tau_1,\bar\tau_2$ are the unique standard real decomposition of $\tau$ obtained in $(1)$.
\item Let $\sigma\in S_d$ be a permutation of type $(d)$, and $\gamma\in S_d$ be an
involution with $\gamma\circ\sigma\circ\gamma=\sigma^{-1}$.
The number of permutations $\tau$ of type $(3,1,\ldots,1)$ such that
$(\sigma,\gamma,\bar\tau_1,\bar\tau_2,(\tau\circ\sigma_1)^{-1})$ produces the monodromy graph in Figure $\ref{fig:1-local-brief}$ is listed in Figure $\ref{fig:1-local-brief}$, where $\bar\tau_1,\bar\tau_2$ are the unique standard real decomposition of $\tau$ obtained in $(1)$.
\end{enumerate}
\end{lemma}

\begin{proof}
Suppose that $(\sigma_1,\gamma,\tau,\sigma_2)$ is a tuple in $\fl^\rb_1((d),(d))$, and $\tau=(abc)$.
The induced permutation of $a,b,c$ in $\sigma_1$ is $\tau$ or $\tau^{-1}$.
If it is $\tau^{-1}$, the composition $\tau\circ\sigma_1$ is not a $d$-cycle.
Hence, the induced permutation of $a,b,c$ in $\sigma_1$ is $\tau$,
and the induced permutation of $a,b,c$ in $\tau\circ\sigma_1$ is $\tau^{-1}$.
From the circle representation, the composition $\tau\circ\sigma_1$ is obtained from $\sigma_1$ as follows.
    First, cut the circle $\sigma_1$ at the positions before $a,b,c$, respectively,
    and get three oriented segments containing $a,b,c$.
    Then we glue the beginning of the oriented segment containing $c$
    to the end of the segment containing $a$,
    and glue the beginning (resp. end) of the oriented segment containing $b$
    to the end (resp. beginning) of the segment containing $c$ (resp. $a$).
    At last, we obtain the circle representing $\tau\circ\sigma_1$.
    By \cite[Lemma $3.12$]{gpmr-2015}, we have the following three cases.
    \begin{enumerate}[$(a)$]
        \item The cycle $\sigma_1$ is an even cycle with two fixed points.
        The symmetric structure of $\sigma_1$ is broken by the above
        procedure to obtain $\tau\circ\sigma_1$.
        We exclude this case.
        \item The cycle $\sigma_1$ is an even cycle with no fixed point.
        In order to get a real cycle by the above procedure,
        the three positions in $\sigma_1$ where we cut at
        have to be a pair of symmetric positions and
        an intersection position of the symmetric axis and the circle $\sigma_1$.
        Suppose that $a,b$ are the two entries in $\sigma_1$ behind the symmetric
        pair of positions. Then the decomposition $\tau=(ac)\circ(ab)$ of $\tau$ is the
        unique decomposition satisfying the requirement in (1).
        \item The cycle $\sigma_1$ is an odd cycle with one fixed point.
        The proof of this case is the same as that in the above case,
        so we omit it.
    \end{enumerate}
Therefore, the statement $(1)$ holds.

By circle representation of cut and join operation and analysis in case $(b)$ and $(c)$,
we obtain statement (2) as follows.
We first consider case $(b)$, that is, $\sigma_1$ is an even cycle with no fixed point.
We draw an edge $E_1$ in blue, and associate weight $e=\len(\sigma_1)$ to it.
From case $(b)$ the two entries in the support of $\bar\tau_1$ are behind a pair
of symmetric positions in $\sigma_1$, so $\bar\tau_1$ is a cut.
The two cycles in $\bar\tau_1\circ\sigma_1$ are two even cycles without fixed point.
We attach two edges $E_2, E_3$ in blue to an endpoint of $E_1$.
The edges $E_2, E_3$ are weighted by the lengths of the two even cycles in $\bar\tau_1\circ\sigma_1$.
The transposition $\bar\tau_2$ is a join. We join the two endpoints of $E_2, E_3$ and glue
a weighted $e$ blue edge $E_4$ to it. At last we get the monodromy graph in
Figure $\ref{fig:1-local-brief}(2)$. The monodromy graph corresponding to case $(c)$ is
obtained similarly (see Figure $\ref{fig:1-local-brief}(1)$).

When $\sigma,\gamma$ are fixed, we compute the number of $3$-cycles $\tau$ such that
$(\sigma_1,\gamma,\bar\tau_1,\bar\tau_2,(\tau\circ\sigma_1)^{-1})$ produces the monodromy graph in Figure $\ref{fig:1-local-brief}$.
The following computation is carried out by using the circle representation.
If $d$ is odd, we have one way to choose the transposition $\bar\tau_1$
such that $\bar\tau_1\circ\sigma_1$ is real and of type $(o_1,e_1)$ (see Figure \ref{fig:real-cut1}(1) ).
There are $2$ choices for $\bar\tau_2$ such that $\bar\tau_2\circ\bar\tau_1\circ\sigma_1$ is real and of type $(d)$.
But only one in the two choices for $\bar\tau_2$ satisfies that $\bar\tau_2\circ\bar\tau_1$ is a $3$-cycle.
Hence, we have $1$ choice for the $3$-cycle $\tau$.
When $d$ is even, we have two ways (resp. one way) to choose the transposition $\bar\tau_1$
such that $\bar\tau_1\circ\sigma_1$ is real and of type $(e_1,e_2)$ with $e_1\neq e_2$ (resp. $e_1=e_2$).
There are $4$ choices for $\bar\tau_2$ such that $\bar\tau_2\circ\bar\tau_1\circ\sigma_1$ is real and of type $(d)$.
But only $2$ of the $4$ choices for $\bar\tau_2$ satisfies that $\bar\tau_2\circ\bar\tau_1$ is a $3$-cycle.
Hence, we have $4$ choices for the $3$-cycle $\tau$ if $e_1\neq e_2$,
and have $2$ choices for the $3$-cycle $\tau$ if $e_1=e_2$.
\end{proof}

\begin{figure}
    \centering
    \begin{tikzpicture}
    \draw[line width=0.4mm,blue] (4.3,-1)--(5.3,-1)--(7.3,0);
    \draw[line width=0.4mm,blue] (5.3,-1)--(6.3,-1);
    \draw[line width=0.4mm,dotted] (7.3,-1.5)--(6.3,-1)--(7.3,-0.5);
    \draw (7,-0.8) node {\tiny $k$} (7,-1.2) node {\tiny $k$} (5.9,-1.2) node {\tiny $2k$} (5.9,-0.5) node {\tiny $e$} (4.7,-1.2) node {\tiny $e+2k$};
    \draw (5.7,-2) node {(1) \small $2$};
    \draw[line width=0.4mm] (7.8,-1)--(8.8,-1)--(10.8,0);
    \draw[line width=0.4mm,blue] (8.8,-1)--(9.8,-1);
    \draw[line width=0.4mm,dotted] (10.8,-1.5)--(9.8,-1)--(10.8,-0.5);
    \draw (10.5,-0.8) node {\tiny $k$} (10.5,-1.2) node {\tiny $k$} (9.4,-1.2) node {\tiny $2k$} (9.4,-0.5) node {\tiny $o$} (8.2,-1.2) node {\tiny $o+2k$};
    \draw (9.2,-2) node {(2) \small $1$};
    \end{tikzpicture}
    \caption{two monodromy graphs of type $((d),(x,y,z))$ with $y=z$.}
    \label{fig:2-local-brief}
\end{figure}

\begin{lemma}
\label{lem:computation-2}
Fix an integer $d$, and a partition $(x,y,z)$ of $d$.
Let $\sigma\in S_d$ be a permutation of type $(d)$, and $\gamma\in S_d$ be an
involution with $\gamma\circ\sigma\circ\gamma=\sigma^{-1}$.
Then, we have the following results.
\begin{enumerate}[$(1)$]
\item Let $\tau$ be the $3$-cycle in a tuple
$(\sigma_1,\gamma,\tau,\sigma_2)\in\fl^\rb_0((d), (x,y,z))$,
then $\tau$ has a unique standard real decomposition $\tau=\bar\tau_2\circ\bar\tau_1$.
\item Given a tuple $(\sigma_1,\gamma,\tau,\sigma_2)\in\fl^\rb_0((d),(x,y,z))$,
the monodromy graph corresponding to
$(\sigma_1,\gamma,\bar\tau_1,\bar\tau_2,\sigma_2)$ is one of the two types depicted in Figure $\ref{fig:2-local-brief}$.
\item Suppose that $d$ and $x,y,z$ are weights of left and right ends 
of a monodromy graph depicted in Figure $\ref{fig:2-local-brief}$, respectively.
Then, the number of permutations $\tau$ of type $(3,1,\ldots,1)$ such that
$(\sigma_1,\gamma,\tau,(\tau\circ\sigma_1)^{-1})\in\fl^\rb_0((d),(x,y,z))$ is listed in Figure $\ref{fig:2-local-brief}$.
\end{enumerate}
\end{lemma}

\begin{proof}
Suppose that $(\sigma_1,\gamma,\tau,\sigma_2)$ is a tuple in $\fl^\rb_0((d),(x,y,z))$, and $\tau=(abc)$. 
The induced permutation of $a,b,c$ in $\sigma_1$ must be $\tau^{-1}$.
By using the circle representation
of cut and join operation, the cycle
$\tau\circ\sigma_1$ is obtained as follows.
We first cut $\sigma_1$ at the three positions before $a$, $b$ and $c$ respectively,
and we get three oriented segments containing $a,b,c$ respectively.
Then we glue the beginning of each oriented segment to the end of itself.
Suppose that $\tau\circ\sigma_1=\bar\sigma_3\circ\bar\sigma_2\circ\bar\sigma_1$ is the disjoint cycle decomposition of $\tau\circ\sigma_1$.
By \cite[Lemma $3.12$]{gpmr-2015}, we have the following three cases.
\begin{enumerate}[$(a)$]
    \item The cycle $\sigma_1$ is an even cycle without fixed point.
    Then the three cycles $\bar\sigma_1$, $\bar\sigma_2$
    and $\bar\sigma_3$ have no fixed point.
    We have two possibilities.
    \begin{itemize}
        \item $\bar\sigma_1$, $\bar\sigma_2$ and $\bar\sigma_3$
        are three even real cycles without fixed point.
        Since the symmetric structures of the cycles are not preserved in this case,
        we exclude this possibility.
        \item There is an even real cycle without fixed point and
        a conjugated pair of cycles in the set $\{\bar\sigma_1,\bar\sigma_2,\bar\sigma_3\}$.
        The three positions, where we perform cuts at, in the circle $\sigma_1$
        are a pair of symmetric blank positions and a blank position
        which is the intersection of the symmetric axis and the circle.
        Suppose that $a$ and $b$ are entries after the pair of symmetric positions respectively.
        Then $\tau$ has a unique decomposition
        $\tau=(ac)\circ(ab)$ satisfying the conditions in (1).
        \end{itemize}
    \item The cycle $\sigma_1$ is an even cycle with $2$ fixed points.
    The above procedure to get circles
    $\bar\sigma_1$, $\bar\sigma_2$ and $\bar\sigma_3$
    from the circle $\sigma_1$ breaks the symmetric structures of the circles,
    so $\sigma_1$ and $\bar\sigma_1\circ\bar\sigma_2\circ\bar\sigma_3$
    can not be both real.
    We exclude this case.
    \item The cycle $\sigma_1$ is an odd cycle with one fixed point.
    To keep the symmetric structures of the circles,
    the three circles $\bar\sigma_1$, $\bar\sigma_2$ and $\bar\sigma_3$
    have to be an odd real cycle and a pair of conjugated cycles.
    The three positions, where we perform cuts at, in the circle $\sigma_1$
    are a pair of symmetric blank positions and a blank position
    which is the intersection of the symmetric axis and the circle.
    We obtain the required unique decomposition of $\tau$ as in the case $(a)$.
\end{enumerate}
Therefore, the statement $(1)$ holds. Moreover, the circle representation
and analysis in case $(a)$ and $(c)$
imply that statement $(2)$ holds. One can draw monodromy graphs
in Figure $\ref{fig:2-local-brief}$ in a similar way to the Figure $\ref{fig:1-local-brief}$,
so we omit the details.

Suppose that $d$ and $x,y,z$ are weights of left and right ends of 
a monodromy graph depicted in Figure $\ref{fig:2-local-brief}$, respectively.
Assume that $\sigma_1$ and $\gamma$ are fixed.
Now we compute the number of $3$-cycles $\tau$ such that 
$(\sigma_1,\gamma,\tau,(\tau\circ\sigma_1)^{-1})\in\fl^\rb_0((d),(x,y,z))$
via circle representation of real cut and join.
When $d$ is even, the symmetric axis intersects at the circle $\sigma_1$ 
at two blank positions, so we have two ways to choose the entry in $\tau$ 
which locates after the intersection of the symmetric axis and the circle $\sigma_1$. 
Since the partition $(x,y,z)$ is fixed, we only have one way to choose 
the pair of conjugated position where we perform the cut.
Hence, we have $2$ choices for the $3$-cycle $\tau$.
If $d$ is odd, the symmetric axis intersects at the circle $\sigma_1$ 
at a unique blank position, so we have $1$ choice for the $3$-cycle $\tau$.
\end{proof}

\begin{figure}[ht]
    \centering
    \begin{tikzpicture}
    \draw[line width=0.4mm,red] (-3,-3.5)--(-2,-3.5);
    \draw[line width=0.4mm] (-2,-3.5)--(0,-2.7);
    \draw[line width=0.4mm] (-2,-3.5)--(-1,-4);
    \draw[line width=0.4mm] (-3,-4.3)--(-1,-4);
    \draw[line width=0.4mm,red] (-1,-4)--(0,-4);
    \draw (-2.8,-4.1) node{\tiny $o$} (-2.8,-3.3) node{\tiny $e$}
    (-1.5,-3.1) node{\tiny $o_1$} (-1.5,-3.6) node{\tiny $o_2$} (-0.5,-3.85) node{\tiny $o_2+o$}
    (-1.6,-4.8) node {(1) 2};
    \draw[line width=0.4mm,blue] (0.5,-3.5)--(1.5,-3.5);
    \draw[line width=0.4mm,blue] (3.5,-2.7)--(1.5,-3.5)--(2.5,-4);
    \draw[line width=0.4mm] (0.5,-4.3)--(2.5,-4)--(3.5,-4);
    \draw (0.7,-4.1) node{\tiny $o$} (0.7,-3.3) node{\tiny $e$}
    (2,-3.1) node{\tiny $e_1$} (2,-3.6) node{\tiny $e_2$} (3,-3.85) node{\tiny $e_2+o$};
    \draw (1.9,-4.8) node {(2) 2};
    \draw[line width=0.4mm,blue] (7.3,-2.7)--(5.3,-3.5);
    \draw[line width=0.4mm] (4.3,-3.5)--(5.3,-3.5)--(6.3,-4)--(7.3,-4);
    \draw[line width=0.4mm,blue] (6.3,-4)--(4.3,-4.3);
    \draw (4.5,-4.1) node{\tiny $e$} (4.5,-3.3) node{\tiny $o$}
    (5.8,-3.1) node{\tiny $e_1$} (5.8,-3.6) node{\tiny $o_1$} (6.8,-3.85) node{\tiny $o_1+e$};
    \draw (5.7,-4.8) node {(3) 2};
    \draw[line width=0.4mm] (7.8,-3.5)--(8.8,-3.5)--(10.8,-2.7);
    \draw[line width=0.4mm,blue] (8.8,-3.5)--(9.8,-4);
    \draw[line width=0.4mm,blue] (7.8,-4.3)--(9.8,-4)--(10.8,-4);
    \draw (8,-4.1) node{\tiny $e$} (8,-3.3) node{\tiny $o$}
    (9.3,-3.1) node{\tiny $o_1$} (9.3,-3.6) node{\tiny $e_1$} (10.3,-3.85) node{\tiny $e_1+e$};
    \draw (9.2,-4.8) node {(4) 2};
    \end{tikzpicture}
    \caption{Four types of partitions $\lambda,\mu$ when $d$ is odd.}
    \label{fig:3-local-brief}
\end{figure}

\begin{figure}[ht]
    \centering
    \begin{tikzpicture}
    \draw[line width=0.4mm,blue] (-3,-3.5)--(-2,-3.5)--(0,-2.7);
    \draw[line width=0.4mm,blue] (-2,-3.5)--(-1,-4);
    \draw[line width=0.4mm,blue] (-3,-4.3)--(-1,-4)--(0,-4);
    \draw (-2.8,-4.1) node{\tiny $e_2$} (-2.8,-3.3) node{\tiny $e_1$}
    (-1.5,-3.1) node{\tiny $e_3$} (-1.5,-3.6) node{\tiny $e_4$} (-0.5,-3.85) node{\tiny $e_2+e_4$}
    (-1.6,-4.8) node {(1) 4};
    \draw[line width=0.4mm,red] (0.5,-3.5)--(1.5,-3.5);
    \draw[line width=0.4mm] (3.5,-2.7)--(1.5,-3.5)--(2.5,-4)--(3.5,-4);
    \draw[line width=0.4mm,blue] (0.5,-4.3)--(2.5,-4);
    \draw (0.7,-4.1) node{\tiny $e_2$} (0.7,-3.3) node{\tiny $e_1$}
    (2,-3.1) node{\tiny $o_1$} (2,-3.6) node{\tiny $o_2$} (3,-3.85) node{\tiny $e_2+o_2$};
    \draw (1.9,-4.8) node {(2) 4};
    \draw[line width=0.4mm] (4.3,-3.5)--(5.3,-3.5);
    \draw[line width=0.4mm,blue] (7.3,-2.7)--(5.3,-3.5);
    \draw[line width=0.4mm] (5.3,-3.5)--(6.3,-4);
    \draw[line width=0.4mm] (6.3,-4)--(4.3,-4.3);
    \draw[line width=0.4mm,red] (6.3,-4)--(7.3,-4);
    \draw (4.5,-4.1) node{\tiny $o_2$} (4.5,-3.3) node{\tiny $o_1$}
    (5.8,-3.1) node{\tiny $e_1$} (5.8,-3.6) node{\tiny $o_3$} (6.8,-3.85) node{\tiny $o_2+o_3$};
    \draw (5.7,-4.8) node {(3) 1};
    \draw[line width=0.4mm] (7.8,-3.5)--(8.8,-3.5)--(10.8,-2.7);
    \draw[line width=0.4mm,blue] (8.8,-3.5)--(9.8,-4);
    \draw[line width=0.4mm] (7.8,-4.3)--(9.8,-4)--(10.8,-4);
    \draw (8,-4.1) node{\tiny $o_2$} (8,-3.3) node{\tiny $o_1$}
    (9.3,-3.1) node{\tiny $o_3$} (9.3,-3.6) node{\tiny $e_1$} (10.3,-3.85) node{\tiny $e_1+o_2$};
    \draw (9.2,-4.8) node {(4) 1};
    \end{tikzpicture}
    \caption{Four types of partitions $\lambda,\mu$ when $d$ is even. When $e_1=e_2$ or $o_1=o_2$, the number of $\tau$ is listed in the bracket.}
    \label{fig:4-local-brief}
\end{figure}

\begin{definition}\label{def:relation-graph}
If there is a monodromy graph $G$ in Figure \ref{fig:3-local-brief} or Figure \ref{fig:4-local-brief}
such that its left ends and right ends can be weighted by $a_1,b_1$ and $a_2,b_2$ respectively,
we say $a_1,b_1,a_2,b_2$ satisfy {\it the relation encoded in $G$}.
\end{definition}

\begin{lemma}
\label{lem:computation-3}
Fix an integer $d$, and two partitions $(a_1,b_1)$, $(a_2,b_2)$ of $d$.
Let $\sigma_1\in S_d$ be a permutation of type $(a_1,b_1)$, and $\gamma\in S_d$ be an
involution with $\gamma\circ\sigma_1\circ\gamma=\sigma_1^{-1}$.
Suppose that $\sigma_1=\bar\sigma_2\circ\bar\sigma_1$ is the disjoint cycle decomposition of $\sigma_1$.
We have the following results.
\begin{enumerate}[$(1)$]
\item Let $\tau$ be the $3$-cycle in a tuple
$(\sigma_1,\gamma,\tau,\sigma_2)\in\fl^\rb_0((a_1,b_1), (a_2,b_2))$,
then $\tau$ has a unique standard real decomposition.
\item Given a tuple $(\sigma_1,\gamma,\tau,\sigma_2)\in\fl^\rb_0((a_1,b_1), (a_2,b_2))$,
the monodromy graph corresponding to
$(\sigma_1,\gamma,\bar\tau_1,\bar\tau_2,\sigma_2)$ is one of the four types depicted in Figure $\ref{fig:3-local-brief}$
(resp. Figure $\ref{fig:4-local-brief}$) when $d$ is odd (resp. even).
\item Suppose that $a_1,b_1$ (resp. $a_2,b_2$) are weights of left (resp. right) ends of one monodromy graph
depicted in Figure $\ref{fig:3-local-brief}$ or Figure $\ref{fig:4-local-brief}$.
Label $\bar\sigma_1$ and $\bar\sigma_2$ to the left ends of one monodromy graph
depicted in Figure $\ref{fig:3-local-brief}$ or Figure $\ref{fig:4-local-brief}$.
Then, the number of $3$-cycles $\tau$ of type $(3,1,\ldots,1)$ such that
$(\sigma_1,\gamma,\tau,(\tau\circ\sigma_1)^{-1})\in\fl^\rb_0((a_1,b_1), (a_2,b_2))$ is listed in
Figure $\ref{fig:3-local-brief}$ and Figure $\ref{fig:4-local-brief}$.
\end{enumerate}
\end{lemma}

\begin{proof}
Suppose that $(\sigma_1,\gamma,\tau,\sigma_2)$ is a tuple in $\fl^\rb_0((a_1,b_1), (a_2,b_2))$.
Assume $\tau=(abc)$.
From proofs of Lemma \ref{lem:computation-1} and Lemma \ref{lem:computation-2},
the three entries in the support of $\tau$ are not totally contained in the support of one cycle among $\bar\sigma_1$, $\bar\sigma_2$.
Assume that $\bar\sigma_1$ contains two elements $a,b$ in the support of $\tau$,
and $\bar\sigma_2$ contains the remaining one number $c$.
By using the circle representation
of cut and join operation, the cycle
$\tau\circ\bar\sigma_2\circ\bar\sigma_1$ is obtained as follows.
We first cut $\bar\sigma_1$ at the two positions before $a$ and $b$ respectively,
and cut $\bar\sigma_2$ at the position before $c$.
Now we achieve three oriented segments containing $a,b,c$ respectively.
Then we glue the beginning (resp. end) of the oriented segment containing $c$
to the end (resp. beginning) of the oriented segment containing $a$.
At last, we glue the beginning of the oriented segment containing $b$
to the end of itself.

Recall that both $\sigma_1$ and $\tau\circ\sigma_1$ are real.
It follows from \cite[Lemma $3.12$]{gpmr-2015} that
if $\gamma\circ\bar\sigma_1\circ\gamma\neq\bar\sigma_1^{-1}$
or $\gamma\circ\bar\sigma_2\circ\gamma\neq\bar\sigma_2^{-1}$,
the above procedure to get $\tau\circ\bar\sigma_1\circ\bar\sigma_2$ can not
produce real cycles.
Hence, $\gamma\circ\bar\sigma_i\circ\gamma=\bar\sigma_i^{-1}$ for $i=1,2$.
By \cite[Lemma $3.12$]{gpmr-2015}, we have the following three cases:
\begin{enumerate}[(a)]
    \item $\bar\sigma_2$ is a real even cycle with two fixed points.
    The above procedure to get $\tau\circ\bar\sigma_1\circ\bar\sigma_2$
    breaks the symmetric structure of the cycle in
    $\tau\circ\bar\sigma_1\circ\bar\sigma_2$ which contains the
    entries of $\bar\sigma_2$,
    so $\tau\circ\sigma$ is not real which contradicts to the assumption.
    We exclude the case that $\bar\sigma_2$ is a real even cycle with two fixed points.
    \item $\bar\sigma_2$ is a real even cycle with no fixed point.
    In order to preserve the symmetric structure of the cycle,
    which contains the entries of $\bar\sigma_2$,
    in $\tau\circ\bar\sigma_1\circ\bar\sigma_2$,
    the position in $\bar\sigma_2$ where we cut at
    has to be the intersection of the symmetric axis and the circle $\bar\sigma_2$.
    Once the position to take a cut in $\bar\sigma_2$ is fixed,
    the two positions in $\bar\sigma_1$
    must be symmetric according to the symmetric axis of the circle $\bar\sigma_1$.
    These two positions in $\bar\sigma_1$ determine a unique
    decomposition of $\tau=(ac)\circ(ab)$ into transpositions satisfying the conditions in (1).
    \item $\bar\sigma_2$ is an odd cycle with one fixed point.
    The proof of this case is the same as that in the case that
    $\bar\sigma_2$ is a real even cycle with no fixed point,
    so we omit it.
\end{enumerate}
Therefore, statement $(1)$ holds. Moreover,
the description of the circle representation of cut and join operation in the first paragraph
and analysis in case $(b)$ and $(c)$
imply that statement $(2)$ holds.
Monodromy graphs
in Figure $\ref{fig:3-local-brief}$ is drawn similarly to Figure $\ref{fig:1-local-brief}$.

Now we show that the statement $(3)$ is true.
Note that $\sigma_1\in S_d$ and the involution $\gamma\in S_d$ are fixed now.
Assume that $\tau=(abc)$ is a permutation such that
$(\sigma_1,\gamma,\tau,(\tau\circ\sigma_1)^{-1})\in\fl^\rb_0((a_1,b_1), (a_2,b_2))$.
Suppose that two entries $a,b$ are contained in $\bar\sigma_1$, and $c$ is contained in $\bar\sigma_2$.
If $a_1,b_1,a_2,b_2$ satisfy the relation encoded in the graph $(1)$ or $(2)$ in Figure \ref{fig:3-local-brief},
$\bar\sigma_2$ is an odd cycle and $\bar\sigma_1$ is an even cycle.
There is only one choice for the number $c$
which is the entry in $\bar\sigma_2$ after the intersection of the symmetric axis and the circle.
There are two choices for the symmetric positions in $\bar\sigma_1$ to take a cut,
so we have two choices for $\tau$.
If $a_1,b_1,a_2,b_2$ satisfy the relation encoded in the graph $(3)$ or $(4)$ in Figure \ref{fig:3-local-brief},
$\bar\sigma_2$ is an even cycle with no fixed point and $\bar\sigma_1$ is an odd cycle.
There are two choices for the number $c$
which are entries in $\bar\sigma_2$ after the two intersections of the symmetric axis and the circle.
There is only one choice for the symmetric position in $\bar\sigma_1$ to take a cut,
so we have two choices for $\tau$ too.
If $a_1,b_1,a_2,b_2$ satisfy the relation encoded in the graph $(1)$ or $(2)$ in Figure \ref{fig:4-local-brief},
$\bar\sigma_2$ is an even cycle with no fixed point and $\bar\sigma_1$ is an even cycle.
There are two choices for the number $c$
which are entries in $\bar\sigma_2$ after the two intersections of the symmetric axis and the circle.
There are two choices for the symmetric positions in $\bar\sigma_1$ to take a cut,
so we have four choices for $\tau$.
If $a_1,b_1,a_2,b_2$ satisfy the relation encoded in the graph $(3)$ or $(4)$ in Figure \ref{fig:4-local-brief},
$\bar\sigma_2$ and $\bar\sigma_1$ are two odd cycles.
There is only one choice for the number $c$
which is the entry in $\bar\sigma_2$ after the intersection of the symmetric axis and the circle.
There is only one choice for the symmetric position in $\bar\sigma_1$ to take a cut,
so we have one choice for $\tau$.
\end{proof}

Let $\widetilde\rl_g(\lambda,\mu)$ be the set of isomorphism classes of real ramified coverings $(\pi,\iota)$
corresponding to tuples in $\widetilde\fl^\rb_g(\lambda,\mu)$ which are
obtained by  \cite[Construction 2.4]{gpmr-2015}. 
From \cite[Construction 2.4]{gpmr-2015},  
real ramified coverings $(\pi,\iota)\in\widetilde\rl_g(\lambda,\mu)$
are ramified at $0$ and $\infty$ with ramification profiles $\lambda$ and $\mu$, respectively,
and ramified at two simple positive branch points.

\begin{proposition}\label{prop:resolution1}
Fix two integers $g\geq0$, $d\geq3$, and two partitions $\lambda$, $\mu$
of $d$ with $l(\lambda)+l(\mu)+2g-2=2$.
Then, there is a bijection
$$
\Phi:\rl_g(\lambda,\mu)\to\widetilde\rl_g(\lambda,\mu).
$$
\end{proposition}

\begin{proof}
From the proof of \cite[Lemma 2.3]{gpmr-2015},
there is a bijection
$$
\widetilde\Psi:\widetilde\fl^\rb_g(\lambda,\mu)/S_d\to\widetilde\rl_g(\lambda,\mu).
$$
Lemma \ref{lem:computation-1}--\ref{lem:computation-3} imply that there is a bijection $\Psi_1$ between the set
$\fl^\rb_g(\lambda,\mu)$ and $\widetilde\fl^\rb_g(\lambda,\mu)$.
The bijection $\Psi_1$ induces a bijection $\Psi_2:\fl^\rb_g(\lambda,\mu)/S_d\to\widetilde\fl^\rb_g(\lambda,\mu)/S_d$.
The bijection $\Phi$ is $\widetilde\Psi\circ\Psi_2\circ\Psi^{-1}$,
where $\Psi$ is the bijection in Lemma \ref{lem:rep}.
\end{proof}

\begin{definition}\label{def:rs-resolution}
Let $(\pi,\iota)\in\rl_g(\lambda,\mu)$ be a real ramified covering.
The covering $\Phi(\pi,\iota)\in\widetilde\rl_g(\lambda,\mu)$ is called the {\it real simple resolution} of $(\pi,\iota)$.
We use $(\tilde\pi,\tilde\iota)$ to denote $\Phi(\pi,\iota)$.
\end{definition}

The following lemma is very useful for computing local Hurwitz numbers.
\begin{lemma}
\label{lem:fixed-points}
Let $\pi:\Sigma_g\to\cb\pb^1$ be a real ramified covering in $\rl_g(\lambda,\mu)$ 
which corresponds to a tuple $A=(\sigma_1,\gamma,\tau,\sigma_2)\in\fl^\rb_g(\lambda,\mu)$ under the bijective map $\Psi$.
Let $\tilde 0\in\pi^{-1}(0)$ (resp. $\tilde\infty\in\pi^{-1}(\infty)$) be a real ramified point with
even ramification index $\lambda_i$ (resp. $\mu_j$).
Suppose that real cycles $\sigma_1'$ and $\sigma_2'$ are assigned to $\tilde 0$ and $\tilde\infty$,
respectively, under the monodromy representation.
Let $B, B'\subset\Sigma_g$ be two small real circles around $\tilde 0$ and $\tilde\infty$ respectively.
Then, the following statements hold.
\begin{enumerate}[$(1)$]
    \item If $\sigma_1'$ (resp. $\sigma_2'$) has two fixed points, the two real points in $B$ (resp. $B'$)
    are mapped to the negative real point in $\pi(B)$ (resp. $\pi(B')$). Otherwise,
    the two real points in $B$ (resp. $B'$) are mapped to the positive real point in $\pi(B)$ (resp. $\pi(B')$).
    \item The preimage $\pi^{-1}(p)$ consists of $k$ real points and $\frac{d-k-2}{2}$
    pairs of conjugated points, where $k$ is the number of elements of $\{1,\ldots,d\}$
    which are not in the support of $\sigma_1\circ\gamma$, and $p$ is the fixed point in
    Construction $\ref{const-1}$.
\end{enumerate}
\end{lemma}

\begin{proof}
Suppose that the real circle $B$ is small enough such that the real points
in $\pi(B)$ are both larger than the base point $x_0$ (see Construction $\ref{const-1}$).
Let $x_1$ be the negative real point in $\pi(B)$.
Denote by $s$ the real line segment from $x_0$ to $x_1$.
Note that in Construction \ref{const-1} the set $\pi^{-1}(x_0)$ is labelled as $\{q_1,\ldots,q_n\}$.
Lift the path $s$ to a path $\tilde s_i$ in $\Sigma_g$ with starting point $q_i$, where $i=1,\ldots,d$.
We denote by $y_i$ the end point of $\tilde s_i$.
Hence, the path $s$ induces a labelling on $\pi^{-1}(x_1)$.
The real involution $\iota$ on $\Sigma_g$ induces an involution
$\gamma'$ on the labelled set $\pi^{-1}(x_1)=\{y_1,\ldots,y_d\}$.
By Construction \ref{const-1}, the involution $\gamma'$ is equal to $\gamma$.
If $i\in\{1,\ldots,d\}$ is not in the support of $\gamma$,
the unramified point $y_i$ is a real point, \textit{i.e.} $\iota(y_i)=y_i$.
Consider the monodromy representation
$$
\Phi_{x_1}:\pi_1(\cb\pb^1\setminus\{0,p,\infty\},x_1)\to S_d.
$$
Note that under the monodromy representation $\Phi_{x_1}$,
the point $\tilde0\in\pi^{-1}(0)$ is assigned to the cycle $\sigma_1'$ too.
Assume that $\sigma_1'$ has two fixed points with respect to $\gamma'$, 
then the two fixed points of $\sigma_1'$ are the labels of the two real points in $B\cap\pi^{-1}(x_1)$.
Hence, the two real points in $B$ are mapped to $x_1$ (see Figure \ref{fig:map-F}).
If $\sigma_1'$ has no fixed point, the two real points in $B$ must be mapped to
the positive real point of $\pi(B)$.
The relation between real points in $B'$ and $\pi(B')$ can be proved similarly,
so we omit it.
We obtain the statement $(1)$ immediately.
\begin{figure}[H]
    \centering
    \begin{tikzpicture}
    \draw (0,0)--(5,0);
    \draw[decoration={markings, mark=at position 0.125 with {\arrow{>}}},
        postaction={decorate}
        ] (2,0) circle (0.3);
    \foreach \Point in {(1,0),(1.7,0), (2,0),(4,0)}
    \draw[fill=black] \Point circle (0.05);
    \draw (1,-0.3) node{\tiny{$x_0$}} (1.7,-0.3) node{\tiny{$x_1$}} (2.4,-0.4) node{\tiny{$\pi(B)$}};
    \draw (5.5,0) node{$\rb$};
    \draw (2,-0.15) node{\tiny{$0$}} (1.4,-0.15) node{\tiny{$s$}} (4,-0.15) node{\tiny{$p$}};
    \draw (0.2,1.7) .. controls (0.7,2.7) and (2,2.3) .. (2.5,2);
    \draw[bend right,-] (2.5,2) to (3,2);
    \draw (3,2) .. controls (3.5,2.3) and (4.5,2.7) .. (5,1.7);
    \draw[bend right,-] (0.2,1.7) to (0.2,1.3);
    \draw (0.2,1.3) .. controls (0.6,0.5) and (2,0.7) .. (2.5,1);
    \draw[bend left,-] (2.5,1) to (3,1);
    \draw (3,1) .. controls (3.5,0.7) and (4.5,0.7) .. (5,1.3);
    \draw (5,1.7) .. controls (5.1,1.5) and (5.1,1.45) .. (5,1.3);
    \draw[bend left,-] (3.5,1.4) to (4.5,1.4);
    \draw[bend right,-] (3.5,1.5) to (4.5,1.5);
    \draw (1.5,1.95) circle (0.3);
    \draw[fill=black] (1.5,1.95) circle (0.05);
    \foreach \Point in {(1.5,2.25),(1.2,1.95), (1.8,1.95),(1.5,1.65),(0.8,1.5),(2.5,1.5),(3.5,2),(4.5,1.9)}
    \draw[fill=red] \Point circle (0.03);
    \foreach \Point in {(1.3,2.25),(1,1.9), (2,1.9),(1.2,1.65),(0.8,1.3),(2.5,1.3),(3.5,1.8),(4.5,1.7)}
    \draw[fill=blue] \Point circle (0.03);
    \draw (0,1.7) node{\tiny$\Sigma_g$} (1.65,1.95) node{\tiny$\tilde0$} (2,2.1) node{\tiny$B$};
    \end{tikzpicture}
    \caption{Covering map $\pi:\Sigma_g\to\cb\pb^1$: $\pi^{-1}(x_0)$ (resp. $\pi^{-1}(x_1)$) are colored in blue (resp. red).}
    \label{fig:map-F}
\end{figure}

Denote by $x_2$ the positive real point in $\pi(B)$.
In the affine plane $\cb\pb^1\setminus\{\infty\}$, the real line $\rb^1$ divides $\pi(B)$ into two parts.
Let $t$ be the upper half circle of $\pi(B)$ from $x_1$ to $x_2$.
The path $t\cdot s$ induces a labelling on the preimage $\pi^{-1}(x_2)=\{z_1,\ldots,z_d\}$.
In fact, let $\tilde{t}_j$ be the lift of $t\cdot s$ in $\Sigma_g$ with starting point $q_j$,
then its endpoint is labelled as $z_j$.
Let $\gamma''$ be the involution on the labelled set $\pi^{-1}(x_2)$ 
determined by the relation $\iota(z_i)=z_{\gamma''(i)}$.
Let $\widetilde{t}_i$ be the lift of $t\cdot s$ with endpoint $z_i$.
By the labelling rule, the path $\widetilde{t}_i$ is started at $q_i$.
Denote by $\conj(t\cdot s)$ the conjugation of the path $t\cdot s$.
It is easy to see $\conj(t\cdot s)=\conj(t)\cdot s$.
Let $\widetilde\conj(t)_i$ be the lift of $\conj(t\cdot s)$ with starting point $q_{\gamma(i)}$.
By the Construction \ref{const-1}, $\iota(z_i)$ is the endpoint of $\widetilde\conj(t)_i$.
Suppose that $\iota(z_i)$ is labelled as $z_j$.
Since $\pi(\widetilde{t}_j^{-1}\cdot\widetilde\conj(t)_i)$ is a positively oriented loop around $0$ with base point $x_0$,
we have $q_j=q_{\sigma_1\circ\gamma(i)}$ (see Figure \ref{fig:involution}).
Therefore, the involution $\gamma''$ is equal to $\sigma_1\circ\gamma$.
\begin{figure}[H]
    \centering
    \begin{tikzpicture}
    \draw (0,0)--(5,0);
    \draw[decoration={markings, mark=at position 0.125 with {\arrow{>}}},
        postaction={decorate}
        ] (2,0) circle (0.3);
    \foreach \Point in {(1,0),(1.7,0), (2,0), (2.3,0),(4,0)}
    \draw[fill=black] \Point circle (0.05);
    \draw (1,-0.3) node{\tiny{$x_0$}} (1.5,0.1) node{\tiny{$x_1$}} (2.5,0.1) node{\tiny{$x_2$}} (2.55,-0.2) node{\tiny{$\pi(B)$}}  (2,0.4) node{\tiny{$t$}} (1.9,-0.4) node{\tiny{$\conj(t)$}};
    \draw (5.5,0) node{$\rb$};
    \draw (2,-0.15) node{\tiny{$0$}} (1.4,-0.15) node{\tiny{$s$}} (4,-0.15) node{\tiny{$p$}};
    \draw (0.2,1.7) .. controls (0.7,2.7) and (2,2.3) .. (2.5,2);
    \draw[bend right,-] (2.5,2) to (3,2);
    \draw (3,2) .. controls (3.5,2.3) and (4.5,2.7) .. (5,1.7);
    \draw[bend right,-] (0.2,1.7) to (0.2,1.3);
    \draw (0.2,1.3) .. controls (0.6,0.5) and (2,0.7) .. (2.5,1);
    \draw[bend left,-] (2.5,1) to (3,1);
    \draw (3,1) .. controls (3.5,0.7) and (4.5,0.7) .. (5,1.3);
    \draw (5,1.7) .. controls (5.1,1.5) and (5.1,1.45) .. (5,1.3);
    \draw[bend left,-] (3.5,1.4) to (4.5,1.4);
    \draw[bend right,-] (3.5,1.5) to (4.5,1.5);
    \draw[bend left,-,decoration={markings, mark=at position 0.4 with {\arrow{>}}},
        postaction={decorate}] (2.5,1.5) to (3.5,2);
    \draw[bend right,-,decoration={markings, mark=at position 0.5 with {\arrow{>}}},
        postaction={decorate}] (3.2,1.5) to (3.5,2);
    \draw (1.5,1.95) circle (0.3);
    \draw[fill=black] (1.5,1.95) circle (0.05);
    \foreach \Point in {(1.5,2.25),(1.2,1.95), (1.8,1.95),(1.5,1.65),(0.8,1.5),(2.2,1.3),(3.5,2),(4.5,1.9)}
    \draw[fill=green] \Point circle (0.03);
    \foreach \Point in {(1.3,2.25),(1,1.9), (2,1.9),(1.2,1.65),(0.8,1.3),(2.5,1.5),(3.2,1.5),(4.6,1.7)}
    \draw[fill=blue] \Point circle (0.03);
    \draw (0,1.7) node{\tiny$\Sigma_g$} (1.65,1.95) node{\tiny$\tilde0$}
    (2,2.1) node{\tiny$B$} (2.5,1.3) node{\tiny$q_j$} (2.8,1.55) node{\tiny$\tilde t_j$}
    (3.2,1.3) node{\tiny$q_{\gamma(i)}$} (3.9,1.7) node{\tiny$\widetilde\conj(t)_i$}
    (3.6,2.1) node{\tiny$z_j$};
    \end{tikzpicture}
    \caption{Covering map $\pi:\Sigma_g\to\cb\pb^1$: $\pi^{-1}(x_0)$ (resp. $\pi^{-1}(x_2)$) are colored in blue (resp. green).}
    \label{fig:involution}
\end{figure}
If $i\in\{1,\ldots,d\}$ is not in the support of $\gamma''$,
the unramified point $z_i$ is a real point, \textit{i.e.} $\iota(z_i)=z_i$.
The preimage $\pi^{-1}(x_2)$ consists of $k$ real points and $\frac{d-k}{2}$ pairs of conjugated points,
where $k$ is the number of elements of $\{1,\ldots,d\}$ which are not in the support of $\sigma_1\circ\gamma$.
Let $l$ be the real line segment from $x_2$ to $p$.
Note that the preimage of any point $x\in l^\circ$ consists of $k$ real points and $\frac{d-k}{2}$ pairs of conjugated points.
Suppose $x\in l^\circ$ closes to $p$ enough.
Under the monodromy representation
$$
\Phi_{x}:\pi_1(\cb\pb^1\setminus\{0,p,\infty\},x)\to S_d,
$$
a real circle $l_p$ around $p$ with based point $x$ corresponds to a $3$-cycle $\tau'$
satisfying $\gamma''\circ\tau'\circ\gamma''=(\tau')^{-1}$.
From \cite[Lemma $3.12$]{gpmr-2015}, there is a real point and a pair of conjugated points of $\pi^{-1}(x)$
which are contained in a small real circle around the ramified point in $\pi^{-1}(p)$.
Hence, when $x\to p$ along the real line segment $l$, one real point and a pair of conjugated points in $\pi^{-1}(x)$
converge to the real ramified point in $\pi^{-1}(p)$,
and the second statement of Lemma \ref{lem:fixed-points} holds.
\end{proof}

\section{Computation of local Hurwitz numbers}
\label{sec:4-1}

In this section we use monodromy representation and circle
representation of real cycles to compute local Hurwitz numbers.

Let $T$ be an open three-valent line, and $p\in\rb\pb^1\setminus\{0,\infty\}$ be a positive point.
Equip $T$ with the standard real structure $(\iota_T, \sal)$, where $\iota_T=\id$,
and $\sal$ is $\cb\pb^1$ with three real boundary circles around punctures at $0,p,\infty\in\rb\pb^1$.
We consider the case that the ramification profiles associated to the three leaves of $T$ are
$(d)$, $(d)$ and $(3,1,\ldots,1)$.
The Riemann-Hurwitz condition implies that any genus $1$ real tropical cover $\varphi:C\to T$ has
one inner vertex $v$ to which all $d$ ends are adjacent.
Let $l$ be the leaf of $T$ with ramification $(3,1,\ldots,1)$.
Denote by $\alpha(\varphi)$ the number of real automorphisms of $\varphi$ which are non-trivial only on $\varphi^{-1}(l)$.

\begin{figure}
    \centering
    \begin{tikzpicture}
    \draw (0,0.5)--(2,0)--(4,0.5);
    \draw (2,0)--(0.5,-1);
    \draw (1,-1)--(2,0)--(1.4,-1);
    \draw (2.6,-1)--(2,0)--(3,-1);
    \draw[fill=black] (2,0) circle (0.03);
    \draw[line width=0.9pt,dotted] (1.8,-0.5)--(2.2,-0.5);
    \draw[<->,bend right] (0.9,-1.2) to (1.4,-1.2);
    \draw[<->,bend right] (2.6,-1.2) to (3.1,-1.2);
    \draw (0.5,0.25) node{\tiny$d$} (2,0.2) node{\tiny$g=1$} (3.5,0.25) node{\tiny$d$} (0.5,-1.1) node{\tiny$3$} (1,-1.1) node{\tiny$1$} (1.4,-1.1) node{\tiny$1$} (2.6,-1.1) node{\tiny$1$} (3,-1.1) node{\tiny$1$};
    \foreach \Point in {(6,0),(8,0)}
    \draw \Point ellipse (0.2 and 0.5);
    \draw (7,-1.1)  ellipse(0.4 and 0.2);
    \draw[bend right] (6,0.5) to (8,0.5);
    \draw (6,-0.5)..controls(6.5,-0.5) and (6.6,-0.5)..(6.6,-1.1);
    \draw (7.4,-1.1)..controls(7.4,-0.5) and (7.5,-0.5)..(8,-0.5);
    \draw (5.6,0) node{\tiny$(d)$} (8.4,0) node{\tiny$(d)$} (8.1,-1.2) node{\tiny$(3,1,\ldots,1)$} (4.6,-1.6) node{$(1)$ $d$ is odd};
    \draw (0,-2)--(2,-2.5)--(4,-2);
    \draw (2,-2.5)--(0.5,-3.5);
    \draw (2,-2.5)--(0.7,-3.5);
    \draw (1,-3.5)--(2,-2.5)--(1.4,-3.5);
    \draw (2.6,-3.5)--(2,-2.5)--(3,-3.5);
    \draw[fill=black] (2,-2.5) circle (0.03);
    \draw[line width=0.9pt,dotted] (1.8,-3)--(2.2,-3);
    \draw[<->,bend right] (0.9,-3.7) to (1.4,-3.7);
    \draw[<->,bend right] (2.6,-3.7) to (3.1,-3.7);
    \draw (0.5,-2.25) node{\tiny$d$} (2,-2.3) node{\tiny$g=1$} (3.5,-2.25) node{\tiny$d$} (0.5,-3.6) node{\tiny$3$} (0.7,-3.6) node{\tiny$1$} (1,-3.6) node{\tiny$1$} (1.4,-3.6) node{\tiny$1$} (2.6,-3.6) node{\tiny$1$} (3,-3.6) node{\tiny$1$};
    \foreach \Point in {(6,-2.5),(8,-2.5)}
    \draw \Point ellipse (0.2 and 0.5);
    \draw (7,-3.6)  ellipse(0.4 and 0.2);
    \draw[bend right] (6,-2) to (8,-2);
    \draw (6,-3)..controls(6.5,-3) and (6.6,-3)..(6.6,-3.6);
    \draw (7.4,-3.6)..controls(7.4,-3) and (7.5,-3)..(8,-3);
    \foreach \Point in {(6,-3),(8,-3)}
    \draw[fill=red] \Point circle (0.05);
    \draw (5.6,-2.5) node{\tiny$(d)$} (8.4,-2.5) node{\tiny$(d)$} (8.1,-3.7) node{\tiny$(3,1,\ldots,1)$} (4.6,-4) node{$(2)$ $d$ is even};
    \end{tikzpicture}
    \caption{The map $F$ used to define $H^\rb(\varphi,v)$ is encoded by marking the points in $\im(F)$ in red.}
    \label{fig:1-local}
\end{figure}

\begin{figure}
    \centering
    \begin{tikzpicture}
    \draw (0,0)--(1.5,0);
    \draw (2,0) ellipse (0.5 and 0.2);
    \draw (4,0)--(2.5,0);
    \draw (0,-1)--(1.5,0)--(0.3,-1);
    \draw (0.6,-1)--(1.5,0)--(0.9,-1);
    \draw (1.5,-1)--(1.5,0)--(1.8,-1);
    \draw[line width=0.9pt,dotted] (1.2,-0.5)--(1.5,-0.5);
    \draw[<->,bend right] (0.6,-1.2) to (0.9,-1.2);
    \draw[<->,bend right] (1.5,-1.2) to (1.8,-1.2);
    \draw (2,-1)--(2.5,0)--(2.3,-1);
    \draw (2.6,-1)--(2.5,0)--(2.9,-1);
    \draw (3.5,-1)--(2.5,0)--(3.8,-1);
    \draw[line width=0.9pt,dotted] (2.7,-0.5)--(3.0,-0.5);
    \draw[<->,bend right] (2.6,-1.2) to (2.9,-1.2);
    \draw[<->,bend right] (3.5,-1.2) to (3.8,-1.2);
    \draw (0.5,0.15) node{\tiny$d$} (2,0.3) node{\tiny$o$} (2,-0.3) node{\tiny$e$} (3.5,0.15) node{\tiny$d$} (0,-1.1) node{\tiny$2$} (0.3,-1.1) node{\tiny$1$} (0.6,-1.1) node{\tiny$1$} (0.9,-1.1) node{\tiny$1$} (1.5,-1.1) node{\tiny$1$} (1.8,-1.1) node{\tiny$1$} (2,-1.1) node{\tiny$2$} (2.3,-1.1) node{\tiny$1$} (2.6,-1.1) node{\tiny$1$} (2.9,-1.1) node{\tiny$1$} (3.5,-1.1) node{\tiny$1$} (3.8,-1.1) node{\tiny$1$};
    \draw[decoration={markings, mark=at position 0.5 with {\arrow{<}}},        postaction={decorate},bend right,blue] (6,0.5) to(6,-0.5);
    \draw[bend left,dotted,blue] (6,0.5) to(6,-0.5);
    \draw[decoration={markings, mark=at position 0.5 with {\arrow{<}}},        postaction={decorate},bend right,blue] (8,0.5) to(8,-0.5);
    \draw[bend left,dotted,blue] (8,0.5) to(8,-0.5);
    \draw[decoration={markings, mark=at position 0.25 with {\arrow{>}}},        postaction={decorate},blue] (7,-1.1)  ellipse(0.4 and 0.2);
    \draw[bend right,blue] (6,0.5) to (8,0.5);
    \draw (6,-0.5)..controls(6.5,-0.5) and (6.6,-0.5)..(6.6,-1.1);
    \draw (7.4,-1.1)..controls(7.4,-0.5) and (7.5,-0.5)..(8,-0.5);
    \draw (5.6,0) node{\tiny$(d)$} (8.4,0) node{\tiny$(o,e)$} (10.4,0) node{\tiny$(d)$} (7.1,-1.4) node{\tiny$(2,1,\ldots,1)$} (9.1,-1.4) node{\tiny$(2,1,\ldots,1)$} (4.6,-1.6) node{$(1)$ $d$ is odd};
    \draw (10,0) ellipse (0.2 and 0.5);
    \draw (9,-1.1)  ellipse(0.4 and 0.2);
    \draw[bend right] (8,0.5) to (10,0.5);
    \draw (8,-0.5)..controls(8.5,-0.5) and (8.6,-0.5)..(8.6,-1.1);
    \draw (9.4,-1.1)..controls(9.4,-0.5) and (9.5,-0.5)..(10,-0.5);
    \draw[fill=black] (8,0.5) circle (0.05);
    \draw[fill=red] (8,-0.5) circle (0.05);
    \draw[fill=red] (7.4,-1.1) circle (0.05);
    \draw[fill=red] (8.6,-1.1) circle (0.05);
    \draw[blue] (8,0.5)--(6.8,-0.95);
    \draw (8,0.65) node{\tiny$x_0$} (5.8,0.4) node{\tiny$l_1$} (6.8,-0.7) node{\tiny$l_2$} (7.7,-0.3)node{\tiny$B$};
    \draw (0,-2.5)--(1.5,-2.5);
    \draw (2,-2.5) ellipse (0.5 and 0.2);
    \draw (4,-2.5)--(2.5,-2.5);
    \draw (0,-3.5)--(1.5,-2.5)--(0.2,-3.5);
    \draw (1.5,-2.5)--(0.4,-3.5);
    \draw (0.6,-3.5)--(1.5,-2.5)--(0.9,-3.5);
    \draw (1.5,-3.5)--(1.5,-2.5)--(1.8,-3.5);
    \draw[line width=0.9pt,dotted] (1.2,-3)--(1.5,-3);
    \draw[<->,bend right] (0.6,-3.7) to (0.9,-3.7);
    \draw[<->,bend right] (1.5,-3.7) to (1.8,-3.7);
    \draw (2,-3.5)--(2.5,-2.5)--(2.2,-3.5);
    \draw (2.5,-2.5)--(2.4,-3.5);
    \draw (2.6,-3.5)--(2.5,-2.5)--(2.9,-3.5);
    \draw (3.5,-3.5)--(2.5,-2.5)--(3.8,-3.5);
    \draw[line width=0.9pt,dotted] (2.7,-3)--(3.0,-3);
    \draw[<->,bend right] (2.6,-3.7) to (2.9,-3.7);
    \draw[<->,bend right] (3.5,-3.7) to (3.8,-3.7);
    \draw (0.5,-2.35) node{\tiny$d$} (2,-2.2) node{\tiny$e_1$} (2,-2.8) node{\tiny$e_2$} (3.5,-2.35) node{\tiny$d$} (0,-3.6) node{\tiny$2$} (0.2,-3.6) node{\tiny$1$} (0.4,-3.6) node{\tiny$1$} (0.6,-3.6) node{\tiny$1$} (0.9,-3.6) node{\tiny$1$} (1.5,-3.6) node{\tiny$1$} (1.8,-3.6) node{\tiny$1$} (2,-3.6) node{\tiny$2$} (2.2,-3.6) node{\tiny$1$} (2.4,-3.6) node{\tiny$1$} (2.6,-3.6) node{\tiny$1$} (2.9,-3.6) node{\tiny$1$} (3.5,-3.6) node{\tiny$1$} (3.8,-3.6) node{\tiny$1$};
    \draw (10,-2.5) ellipse (0.2 and 0.5);
    \draw[bend right] (6,-2) to(6,-3);
    \draw[bend left,dotted] (6,-2) to(6,-3);
    \draw[bend right] (8,-2) to(8,-3);
    \draw[bend left,dotted] (8,-2) to(8,-3);
    \draw (7,-3.6)  ellipse(0.4 and 0.2);
    \draw (9,-3.6)  ellipse(0.4 and 0.2);
    \draw[bend right] (6,-2) to (8,-2);
    \draw[bend right] (8,-2) to (10,-2);
    \draw (6,-3)..controls(6.5,-3) and (6.6,-3)..(6.6,-3.6);
    \draw (7.4,-3.6)..controls(7.4,-3) and (7.5,-3)..(8,-3);
    \draw (8,-3)..controls(8.5,-3) and (8.6,-3)..(8.6,-3.6);
    \draw (9.4,-3.6)..controls(9.4,-3) and (9.5,-3)..(10,-3);
    \foreach \Point in {(6,-3),(8,-3),(7.4,-3.6),(8.6,-3.6),(10,-3)}
    \draw[fill=red] \Point circle (0.05);
    \draw (5.6,-2.5) node{\tiny$(d)$} (10.4,-2.5) node{\tiny$(d)$} (8.3,-2.5) node{\tiny$B$} (8,-3.2) node{\tiny\color{red}$2\times$} (7.1,-3.9) node{\tiny$(2,1,\ldots,1)$} (9.1,-3.9) node{\tiny$(2,1,\ldots,1)$} (4.6,-4) node{$(2)$ $d$ is even};
    \end{tikzpicture}
    \caption{The map $F$ used to define real tropical cover is encoded by marking the points in $\im(F)$ in red. If fixed points of two fixed even leaves are mapped to a same fixed point of $B$, we use the symbol `` $2\times$" to indicate that.}
    \label{fig:1-local-res}
\end{figure}

\begin{lemma}
\label{lem:1-local}
Let $\varphi:C\to T$ be a genus $1$ real tropical cover of an open $3$-valent line $T$ with ramification
profiles $(d)$, $(d)$ and $(3,1,\ldots,1)$ over the leaves.
All non-zero local Hurwitz numbers $H^\rb(\varphi,v)$ are listed in the following.
\begin{enumerate}[$(1)$]
    \item If $d$ is odd,
    $H^\rb(\varphi,v)=\frac{d-1}{2}\cdot\alpha(\varphi)=\frac{d-1}{2}\cdot2^{\frac{d-3}{2}}\cdot\left(\frac{d-3}{2}\right)!$,
    where the involution $\iota$ is depicted in Figure $\ref{fig:1-local}(1)$.
    Moreover, any real ramified covering $(\tilde\pi,\tilde\iota)\in\widetilde\rl_1((d),(d))$ defines a real tropical cover depicted in Figure $\ref{fig:1-local-res}(1)$.
    \item If $d$ is even, $H^\rb(\varphi,v)=\frac{d-2}{2}\cdot\alpha(\varphi)=\frac{d-2}{2}\cdot2^{\frac{d-4}{2}}\cdot\left(\frac{d-4}{2}\right)!$,
    where the map $F$, $\iota$ are depicted in Figure $\ref{fig:1-local}(2)$, and $\im(F)$ is coloured in red.
    Moreover, any real ramified covering $(\tilde\pi,\tilde\iota)\in\widetilde\rl_1((d),(d))$ defines a real tropical cover depicted in Figure $\ref{fig:1-local-res}(2)$.
\end{enumerate}
\end{lemma}

\begin{proof}
From Lemma \ref{lem:rep}, we only need to calculate
the number of tuples in $\fl^\rb_1((d),(d))$.
There are $(d-1)!$ permutations of type $(d)$ in the symmetric group $S_d$,
so we have $(d-1)!$ possibilities to choose $\sigma_1$.
From \cite[Lemma $3.12$]{gpmr-2015} or Figure \ref{fig:sym-circle1},
we have $d$ possibilities to choose $\gamma$ if $d$ is odd, and have $\frac{d}{2}$ choices for $\gamma$ if $d$ is even.
There are $\frac{d-1}{2}$ (resp. $\left\lceil\frac{d-2}{4}\right\rceil$) monodromy graphs of type $(1)$ (resp. type $(2)$) in Figure \ref{fig:1-local-brief}. From Lemma \ref{lem:computation-1}(3),
there is $1$ permutation $\tau$
of type $(3,1,\ldots,1)$ for the monodromy graph in Figure \ref{fig:1-local-brief}(1)  when $d$ is odd.
When $d$ is even and $e_1\neq e_2$ (resp. $e_1=e_2$), there are $4$ (resp. $2$) permutations $\tau$
of type $(3,1,\ldots,1)$ for the monodromy graph in Figure \ref{fig:1-local-brief}(2).
Hence, $|\fl^\rb_1((d),(d))|=\frac{d-1}{2}\cdot d!$ if $d$ is odd.
When $d$ is even, $|\fl^\rb_1((d),(d))|=\frac{d-2}{2}\cdot d!$.
Recall that $S_d$ acts on $\fl^\rb_1((d),(d))$ by conjugation (coordinatewise),
so the stabilizer of any tuple $(\sigma_1,\gamma,\tau,\sigma_2)\in\fl^\rb_1((d),(d))$ is trivial.
An orbit of the $S_d$ action contains $d!$ tuples of $\fl^\rb_1((d),(d))$,
and an orbit of the $S_d$ action corresponds to a real ramified covering
in the set $\rl_1((d),(d))$ with trivial automorphism group (see Lemma \ref{lem:rep}).
Therefore, there are $\frac{d-1}{2}$ (resp. $\frac{d-2}{2}$) real ramified coverings in $\rl_1((d),(d))$ if $d$ is odd (resp. even).
From Lemma \ref{lem:fixed-points}(2), the set $\pi^{-1}(p)$ consists of one real point (resp. two real points)
and $\frac{d-3}{2}$ (resp. $\frac{d-4}{2}$) pairs of conjugated points for any real cover $\pi\in\rl_1((d),(d))$ if $d$ is odd (resp. even).
When $d$ is odd (resp. even), the number of labels of the leaves in $C$ mapping to the leaf
$l\subset T$ with ramification profile $(3,1,\ldots,1)$ by $\varphi$ is
$\alpha(\varphi)=2^{\frac{d-3}{2}}\cdot\left(\frac{d-3}{2}\right)!$ (resp. $\alpha(\varphi)=2^{\frac{d-4}{2}}\cdot\left(\frac{d-4}{2}\right)!$).

Let $A=(\sigma_1,\gamma,\tau,\sigma_2)$ be a tuple in $\fl^\rb_1((d),(d))$, and $(\pi,\iota)=\Psi(A)$.
Suppose that $\tilde A=(\sigma_1,\gamma,\bar\tau_1,\bar\tau_2,\sigma_2)\in\widetilde\fl^\rb_1((d),(d))$
is the tuple corresponding to $A$ via Lemma \ref{lem:computation-1}(1).
From Lemma \ref{lem:computation-1}(2), the composition $\bar\tau_1\circ\sigma_1$ consists of two disjoint real cycles: an even cycle $\bar\sigma_1'$ with no fixed point and an odd cycle $\bar\sigma_1''$ (resp. two even cycles $\bar\sigma_1'$ and $\bar\sigma_1''$ with no fixed point) if $d$ is odd (resp. if $d$ is even).
Let $(\tilde\pi:\widetilde\cl\to\widetilde\dl,\tilde\iota)$ be the simple resolution of $(\pi:\cl\to\dl,\iota)$.
Let $B\subset\widetilde\dl$ be the real circle depicted in blue in Figure $\ref{fig:1-local-res}(1)$.
Let $x_0\in B$ be the negative real point in $B$, and $l_1,l_2$ be the oriented loops depicted in blue in Figure $\ref{fig:1-local-res}(1)$.
From \cite[Construction 2.4]{gpmr-2015},
the monodromy representations of $l_1,l_2$ are $\sigma_1,\bar\tau_1$ respectively,
so the monodromy representation of $B$ is $\bar\tau_1\circ\sigma_1=\bar\sigma_1'\circ\bar\sigma_1''$.
Therefore, $\tilde\pi^{-1}(B)$ consists of two circles $B_1,B_2$, and degrees of $\tilde\pi$ restricted on $B_1,B_2$ are $l(\bar\sigma_1'), l(\bar\sigma_1'')$.
The circle $B$ divides the punctured sphere $\widetilde\dl$ into two parts $\widetilde\dl_1$, $\widetilde\dl_2$.
At least one of $\pi^{-1}(\widetilde\dl_1)$ and $\pi^{-1}(\widetilde\dl_2)$ is connected.
Suppose that $\pi^{-1}(\widetilde\dl_1)$ is connected.
Apply Riemann-Hurwitz formula to $\pi:\pi^{-1}(\widetilde\dl_1)\to\widetilde\dl_1$.
We obtain that the genus of $\pi^{-1}(\widetilde\dl_1)$ is 0.
Then $\pi^{-1}(\widetilde\dl_2)$ has to be connected and the genus of $\pi^{-1}(\widetilde\dl_2)$ is also $0$.
The circles $B_1$ and $B_2$ divide the torus $\widetilde\cl$ into two genus 0 components $\widetilde\cl_1,\widetilde\cl_2$.
From \cite[Lemma 5.5]{mr-2015}, the real tropical cover corresponding to $(\tilde\pi,\tilde\iota)$ is depicted in Figure \ref{fig:1-local-res}(1) if $d$ is odd. When $d$ is even the real tropical cover corresponding to $(\tilde\pi,\tilde\iota)$
is obtained similarly (see Figure \ref{fig:1-local-res}(2)).
\end{proof}

The proofs of Lemma \ref{lem:2-local}--Lemma \ref{lem:4-local-brief-version} are
similar to that of Lemma \ref{lem:1-local}, so we postpone them to Appendix \ref{sec:a}.
Let $T$ be an open three-valent line equipped with the standard real structure $(\iota_T, \sal)$ as above.
The ramification profiles assigned to the three leaves of $T$ are
$(d)$, $(a,b,c)$ and $(3,1,\ldots,1)$.
By the Riemann-Hurwitz formula, any genus $0$ real tropical cover $\varphi:C\to T$ has
one inner vertex $v\in C$.
Let $l$ be the leaf of $T$ with ramification $(3,1,\ldots,1)$.
Denote by $\alpha(\varphi)$ the number of real automorphisms of $\varphi$ which are non-trivial only on $\varphi^{-1}(l)$.

\begin{lemma}
\label{lem:2-local}
Let $\varphi:C\to T$ be a genus $0$ real tropical cover of an open $3$-valent line $T$ with ramification
profiles $(d)$, $(a, b, c)$ and $(3,1,\ldots,1)$.
All non-zero local Hurwitz numbers $H^\rb(\varphi,v)$ are listed in the following.
\begin{enumerate}[$(1)$]
    \item If $d$ is odd and two entries in $(a,b,c)$ are equal,
    $
    H^\rb(\varphi,v)=2\cdot\alpha(\varphi)=2^{\frac{d-1}{2}}\cdot\left(\frac{d-3}{2}\right)!,
    $
    where the involution $\iota$ is depicted in Figure $\ref{fig:2-local}(1)$.
    Moreover, $|\widetilde\rl_0((d),(a,b,c))|=1$ and the only real ramified covering $(\tilde\pi,\tilde\iota)\in\widetilde\rl_0((d),(a,b,c))$
    defines a real tropical cover depicted in Figure $\ref{fig:2-local-res}(1)$.
    \item If $d$ is even and two entries in $(a,b,c)$ are equal,
    $
H^\rb(\varphi,v)=2\cdot\alpha(\varphi)=2^{\frac{d-2}{2}}\cdot\left(\frac{d-4}{2}\right)!,
    $
    where the map $F$ and the involution $\iota$ are depicted in Figure $\ref{fig:2-local}(2)$.
    Moreover, the number $|\widetilde\rl_0((d),(a,b,c))|$ is $1$ and the only real ramified covering
    $(\tilde\pi,\tilde\iota)\in\widetilde\rl_0((d),(a,b,c))$
    defines a real tropical cover depicted in Figure $\ref{fig:2-local-res}(2)$.
\end{enumerate}
\end{lemma}

\begin{figure}
    \centering
    \begin{tikzpicture}
    \draw (0,0.5)--(2,0)--(4,0.5);
    \draw (2,0)--(0.5,-1);
    \draw (4,0.1)--(2,0)--(4,-0.3);
    \draw (1,-1)--(2,0)--(1.4,-1);
    \draw (2.6,-1)--(2,0)--(3,-1);
    \draw[line width=0.9pt,dotted] (1.8,-0.5)--(2.2,-0.5);
    \draw[<->,bend right] (0.9,-1.2) to (1.4,-1.2);
    \draw[<->,bend right] (2.6,-1.2) to (3.1,-1.2);
    \draw[<->,bend left] (4.1,0.1) to (4.1,-0.3);
    \draw (0.5,0.25) node{\tiny$d$} (3.3,0.4) node{\tiny$a$} (3.6,0) node{\tiny$b$} (3.5,-0.4) node{\tiny$c$} (0.5,-1.1) node{\tiny$3$} (1,-1.1) node{\tiny$1$} (1.4,-1.1) node{\tiny$1$} (2.6,-1.1) node{\tiny$1$} (3,-1.1) node{\tiny$1$};
    \foreach \Point in {(6,0),(8,0)}
    \draw \Point ellipse (0.2 and 0.5);
    \draw (7,-1.1)  ellipse(0.4 and 0.2);
    \draw[bend right] (6,0.5) to (8,0.5);
    \draw (6,-0.5)..controls(6.5,-0.5) and (6.6,-0.5)..(6.6,-1.1);
    \draw (7.4,-1.1)..controls(7.4,-0.5) and (7.5,-0.5)..(8,-0.5);
    \draw (5.6,0) node{\tiny$(d)$} (8.6,0) node{\tiny$(a,b,c)$} (8.1,-1.2) node{\tiny$(3,1,\ldots,1)$} (4.6,-1.6) node{$(1)$ $d$ is odd and $b=c$};
    \draw (0,-2)--(2,-2.5)--(4,-2);
    \draw (4,-2.4)--(2,-2.5)--(4,-2.8);
    \draw (2,-2.5)--(0.5,-3.5);
    \draw (2,-2.5)--(0.7,-3.5);
    \draw (1,-3.5)--(2,-2.5)--(1.4,-3.5);
    \draw (2.6,-3.5)--(2,-2.5)--(3,-3.5);
    \draw[line width=0.9pt,dotted] (1.8,-3)--(2.2,-3);
    \draw[<->,bend right] (0.9,-3.7) to (1.4,-3.7);
    \draw[<->,bend right] (2.6,-3.7) to (3.1,-3.7);
    \draw[<->,bend left] (4.1,-2.4) to (4.1,-2.8);
    \draw (0.5,-2.25) node{\tiny$d$} (3.3,-2.1) node{\tiny$a$} (3.6,-2.5) node{\tiny$b$} (3.5,-2.9) node{\tiny$c$} (0.5,-3.6) node{\tiny$3$} (0.7,-3.6) node{\tiny$1$} (1,-3.6) node{\tiny$1$} (1.4,-3.6) node{\tiny$1$} (2.6,-3.6) node{\tiny$1$} (3,-3.6) node{\tiny$1$};
    \foreach \Point in {(6,-2.5),(8,-2.5)}
    \draw \Point ellipse (0.2 and 0.5);
    \draw (7,-3.6)  ellipse(0.4 and 0.2);
    \draw[bend right] (6,-2) to (8,-2);
    \draw (6,-3)..controls(6.5,-3) and (6.6,-3)..(6.6,-3.6);
    \draw (7.4,-3.6)..controls(7.4,-3) and (7.5,-3)..(8,-3);
    \foreach \Point in {(6,-3),(8,-3)}
    \draw[fill=red] \Point circle (0.05);
    \draw (5.6,-2.5) node{\tiny$(d)$} (8.6,-2.5) node{\tiny$(a,b,c)$} (8.1,-3.7) node{\tiny$(3,1,\ldots,1)$} (4.6,-4.1) node{$(2)$ $d$ is even and $b=c$};
    \end{tikzpicture}
    \caption{The case $b=c$: points in $\im(F)$ are coloured in red.}
    \label{fig:2-local}
\end{figure}

\begin{figure}
    \centering
    \begin{tikzpicture}
    \draw (0,0)--(1.5,0)--(4,0.5);
    \draw (1.5,0)--(2.5,0);
    \draw (2.5,0.2)--(2.4,0.5);
    \draw (2.6,0.5)--(2.5,0.2)--(2.9,0.5);
    \draw (3.3,0.5)--(2.5,0.2)--(3.6,0.5);
    \draw (2.4,0.6) node{\tiny$1$} (2.6,0.6) node{\tiny$1$}
    (2.9,0.6) node{\tiny$1$} (3.3,0.6) node{\tiny$1$} (3.6,0.6) node{\tiny$1$};
    \draw[<->,bend left] (2.6,0.7) to (2.9,0.7);
    \draw[<->,bend left] (3.3,0.7) to (3.6,0.7);
    \draw[line width=0.9pt,dotted] (2.95,0.5)--(3.25,0.5);
    \draw (4,0)--(2.5,0)--(4,-0.5);
    \draw (0,-1)--(1.5,0)--(0.3,-1);
    \draw (0.6,-1)--(1.5,0)--(0.9,-1);
    \draw (1.5,-1)--(1.5,0)--(1.8,-1);
    \draw[line width=0.9pt,dotted] (1.2,-0.5)--(1.5,-0.5);
    \draw[<->,bend right] (0.6,-1.2) to (0.9,-1.2);
    \draw[<->,bend right] (1.5,-1.2) to (1.8,-1.2);
    \draw (2,-1)--(2.5,0);
    \draw (2.6,-1)--(2.5,0)--(2.9,-1);
    \draw (3.5,-1)--(2.5,0)--(3.8,-1);
    \draw[line width=0.9pt,dotted] (2.7,-0.5)--(3.0,-0.5);
    \draw[<->,bend right] (2.6,-1.2) to (2.9,-1.2);
    \draw[<->,bend right] (3.5,-1.2) to (3.8,-1.2);
    \draw[<->,bend left] (4.1,0) to (4.1,-0.5);
    \draw (0.5,0.15) node{\tiny$d$} (2,0.3) node{\tiny$a$} (2,-0.2) node{\tiny$2b$}
    (3.5,0.15) node{\tiny$b$} (3.8,0.4) node{\tiny$a$} (3.5,-0.2) node{\tiny$c$} (0,-1.1) node{\tiny$2$}
    (0.3,-1.1) node{\tiny$1$} (0.6,-1.1) node{\tiny$1$} (0.9,-1.1) node{\tiny$1$}
    (1.5,-1.1) node{\tiny$1$} (1.8,-1.1) node{\tiny$1$} (2,-1.1) node{\tiny$2$}
    (2.6,-1.1) node{\tiny$1$} (2.9,-1.1) node{\tiny$1$}
    (3.5,-1.1) node{\tiny$1$} (3.8,-1.1) node{\tiny$1$};
    \draw[bend right] (6,0.5) to(6,-0.5);
    \draw[bend left,dotted] (6,0.5) to(6,-0.5);
    \draw[bend right] (8,0.5) to(8,-0.5);
    \draw[bend left,dotted] (8,0.5) to(8,-0.5);
    \draw (7,-1.1)  ellipse(0.4 and 0.2);
    \draw[bend right] (6,0.5) to (8,0.5);
    \draw (6,-0.5)..controls(6.5,-0.5) and (6.6,-0.5)..(6.6,-1.1);
    \draw (7.4,-1.1)..controls(7.4,-0.5) and (7.5,-0.5)..(8,-0.5);
    \draw (5.6,0) node{\tiny$(d)$} (8.3,0) node{\tiny$B$}
    (10.4,0) node{\tiny$(a,b,c)$} (7.1,-1.4) node{\tiny$(2,1,\ldots,1)$} (9.1,-1.4) node{\tiny$(2,1,\ldots,1)$} (4.6,-1.5) node{$(1)$ $a$ is odd, $b=c$};
    \draw (10,0) ellipse (0.2 and 0.5);
    \draw (9,-1.1)  ellipse(0.4 and 0.2);
    \draw[bend right] (8,0.5) to (10,0.5);
    \draw (8,-0.5)..controls(8.5,-0.5) and (8.6,-0.5)..(8.6,-1.1);
    \draw (9.4,-1.1)..controls(9.4,-0.5) and (9.5,-0.5)..(10,-0.5);
    \draw[fill=red] (8,-0.5) circle (0.05);
    \draw[fill=red] (7.4,-1.1) circle (0.05);
    \draw[fill=red] (8.6,-1.1) circle (0.05);
    \draw (0,-2.5)--(1.5,-2.5)--(4,-2);
    \draw (1.5,-2.5)--(2.5,-2.5);
    \draw (2.1,-2)--(2.5,-2.3)--(2.4,-2);
    \draw (2.6,-2)--(2.5,-2.3)--(2.9,-2);
    \draw (3.3,-2)--(2.5,-2.3)--(3.6,-2);
    \draw (2.1,-1.9) node{\tiny$1$} (2.4,-1.9) node{\tiny$1$} (2.6,-1.9) node{\tiny$1$}
    (2.9,-1.9) node{\tiny$1$} (3.3,-1.9) node{\tiny$1$} (3.6,-1.9) node{\tiny$1$};
    \draw[<->,bend left] (2.6,-1.8) to (2.9,-1.8);
    \draw[<->,bend left] (3.3,-1.8) to (3.6,-1.8);
    \draw[line width=0.9pt,dotted] (2.95,-2)--(3.25,-2);
    \draw (4,-2.5)--(2.5,-2.5)--(4,-3);
    \draw (0,-3.5)--(1.5,-2.5)--(0.3,-3.5);
    \draw (0.6,-3.5)--(1.5,-2.5)--(0.9,-3.5);
    \draw (0.15,-3.5)--(1.5,-2.5);
    \draw (1.5,-3.5)--(1.5,-2.5)--(1.8,-3.5);
    \draw[line width=0.9pt,dotted] (1.2,-3)--(1.5,-3);
    \draw[<->,bend right] (0.6,-3.7) to (0.9,-3.7);
    \draw[<->,bend right] (1.5,-3.7) to (1.8,-3.7);
    \draw (2.5,-2.5)--(2.3,-3.5);
    \draw (2.6,-3.5)--(2.5,-2.5)--(2.9,-3.5);
    \draw (3.5,-3.5)--(2.5,-2.5)--(3.8,-3.5);
    \draw[line width=0.9pt,dotted] (2.7,-3)--(3.0,-3);
    \draw[<->,bend right] (2.6,-3.7) to (2.9,-3.7);
    \draw[<->,bend right] (3.5,-3.7) to (3.8,-3.7);
    \draw[<->,bend left] (4.1,-2.5) to (4.1,-3);
    \draw (0.5,-2.35) node{\tiny$d$} (2,-2.2) node{\tiny$a$} (2,-2.6) node{\tiny$2b$} (3.5,-2.35) node{\tiny$b$} (3.8,-2.1) node{\tiny$a$} (3.5,-2.65) node{\tiny$c$} (0,-3.6) node{\tiny$2$} (0.3,-3.6) node{\tiny$1$} (0.6,-3.6) node{\tiny$1$} (0.15,-3.6) node{\tiny$1$} (0.9,-3.6) node{\tiny$1$} (1.5,-3.6) node{\tiny$1$} (1.8,-3.6) node{\tiny$1$}  (2.3,-3.6) node{\tiny$2$} (2.6,-3.6) node{\tiny$1$} (2.9,-3.6) node{\tiny$1$} (3.5,-3.6) node{\tiny$1$} (3.8,-3.6) node{\tiny$1$};
    \draw (10,-2.5) ellipse (0.2 and 0.5);
    \draw[bend right] (6,-2) to(6,-3);
    \draw[bend left,dotted] (6,-2) to(6,-3);
    \draw[bend right] (8,-2) to(8,-3);
    \draw[bend left,dotted] (8,-2) to(8,-3);
    \draw (7,-3.6)  ellipse(0.4 and 0.2);
    \draw (9,-3.6)  ellipse(0.4 and 0.2);
    \draw[bend right] (6,-2) to (8,-2);
    \draw[bend right] (8,-2) to (10,-2);
    \draw (6,-3)..controls(6.5,-3) and (6.6,-3)..(6.6,-3.6);
    \draw (7.4,-3.6)..controls(7.4,-3) and (7.5,-3)..(8,-3);
    \draw (8,-3)..controls(8.5,-3) and (8.6,-3)..(8.6,-3.6);
    \draw (9.4,-3.6)..controls(9.4,-3) and (9.5,-3)..(10,-3);
    \foreach \Point in {(6,-3),(8,-3),(7.4,-3.6),(8.6,-3.6),(10,-3)}
    \draw[fill=red] \Point circle (0.05);
    \draw (5.6,-2.5) node{\tiny$(d)$} (10.4,-2.5) node{\tiny$(a,b,c)$}(8.3,-2.5) node{\tiny$B$}
    (8,-3.2) node{\tiny\color{red}$2\times$} (7.1,-3.9) node{\tiny$(2,1,\ldots,1)$}
    (9.1,-3.9) node{\tiny$(2,1,\ldots,1)$} (4.6,-4) node{$(2)$ $a$ is even, $b=c$};
    \end{tikzpicture}
    \caption{The map $F$ used to define real tropical cover is encoded by marking the points in $\im(F)$ in red. If fixed points of two fixed even leaves are mapped to a same fixed point of $B$, we use the symbol `` $2\times$" to indicate that.}
    \label{fig:2-local-res}
\end{figure}

Let $T$ be an open three-valent line equipped with the standard real structure $(\iota_T, \sal)$ as above.
Let $(a_1,b_1)$ and $(a_2,b_2)$ be two partitions of $d$.
We consider genus $0$ real tropical cover $\varphi:C\to T$ with ramification profiles $(a_1,b_1)$, $(a_2,b_2)$ and $(3,1,\ldots,1)$.
From the Riemann-Hurwitz formula we know that any genus $0$ real tropical cover $\varphi:C\to T$ has
one inner vertex $v\in C$.
Let $l$ be the leaf of $T$ with ramification $(3,1,\ldots,1)$.
Denote by $\alpha(\varphi)$ the number of real automorphisms of $\varphi$ which are non-trivial only on $\varphi^{-1}(l)$.
The local Hurwitz number $H^\rb(\varphi,v)$ depends on the relations between the ramification indexes $a_1,b_1,a_2,b_2$.
In Lemma \ref{lem:3-local-brief-version} and Lemma \ref{lem:4-local-brief-version} we show
that only if relations between $a_1,b_1,a_2,b_2$ can be described by monodromy graphs in
Figure \ref{fig:3-local-brief} or Figure \ref{fig:4-local-brief} (see Definition \ref{def:relation-graph}), the number $H^\rb(\varphi,v)$ is non-trivial.

\begin{figure}
    \centering
    \begin{tikzpicture}
    \draw (0,0.5)--(2,0)--(4,0.5);
    \draw (2,0)--(0.5,-1);
    \draw (0,-0.2)--(2,0)--(4,-0.3);
    \draw (1,-1)--(2,0)--(1.4,-1);
    \draw (2.6,-1)--(2,0)--(3,-1);
    \draw[line width=0.9pt,dotted] (1.8,-0.5)--(2.2,-0.5);
    \draw[<->,bend right] (0.9,-1.2) to (1.4,-1.2);
    \draw[<->,bend right] (2.6,-1.2) to (3.1,-1.2);
    \draw (0.5,0.25) node{\tiny$a_1$} (3.3,0.45) node{\tiny$a_2$} (0.2,0) node{\tiny$b_1$} (3.5,-0.4) node{\tiny$b_2$} (0.5,-1.1) node{\tiny$3$} (1,-1.1) node{\tiny$1$} (1.4,-1.1) node{\tiny$1$} (2.6,-1.1) node{\tiny$1$} (3,-1.1) node{\tiny$1$};
    \foreach \Point in {(6,0),(8,0)}
    \draw \Point ellipse (0.2 and 0.5);
    \draw (7,-1.1)  ellipse(0.4 and 0.2);
    \foreach \Point in {(6,0.5),(8,0.5)}
    \draw[fill=red] \Point circle (0.05);
    \draw[bend right] (6,0.5) to (8,0.5);
    \draw (6,-0.5)..controls(6.5,-0.5) and (6.6,-0.5)..(6.6,-1.1);
    \draw (7.4,-1.1)..controls(7.4,-0.5) and (7.5,-0.5)..(8,-0.5);
    \draw (5.2,0) node{\tiny$(a_1,b_1)$} (8.7,0) node{\tiny$(a_2,b_2)$} (8.1,-1.2) node{\tiny$(3,1,\ldots,1)$} (4.6,-1.6) node{$(1)$ $a_1>b_2$ and $b_1<a_2$.};
    \draw (0,-2)--(2,-2.5)--(4,-2);
    \draw (0,-2.7)--(2,-2.5)--(4,-2.8);
    \draw (2,-2.5)--(0.5,-3.5);
    \draw (2,-2.5)--(0.7,-3.5);
    \draw (2,-2.5)--(0.2,-3.5);
    \draw (1,-3.5)--(2,-2.5)--(1.4,-3.5);
    \draw (2.6,-3.5)--(2,-2.5)--(3,-3.5);
    \draw[line width=0.9pt,dotted] (1.8,-3)--(2.2,-3);
    \draw[<->,bend right] (0.9,-3.7) to (1.4,-3.7);
    \draw[<->,bend right] (2.6,-3.7) to (3.1,-3.7);
    \draw (0.5,-2.25) node{\tiny$a_1$} (3.3,-2.05) node{\tiny$a_2$} (0.4,-2.5) node{\tiny$b_1$} (3.5,-2.9) node{\tiny$b_2$} (0.5,-3.6) node{\tiny$1$} (0.2,-3.6) node{\tiny$3$} (0.7,-3.6) node{\tiny$1$} (1,-3.6) node{\tiny$1$} (1.4,-3.6) node{\tiny$1$} (2.6,-3.6) node{\tiny$1$} (3,-3.6) node{\tiny$1$};
    \foreach \Point in {(6,-2.5),(8,-2.5)}
    \draw \Point ellipse (0.2 and 0.5);
    \draw (7,-3.6)  ellipse(0.4 and 0.2);
    \draw[bend right] (6,-2) to (8,-2);
    \draw (6,-3)..controls(6.5,-3) and (6.6,-3)..(6.6,-3.6);
    \draw (7.4,-3.6)..controls(7.4,-3) and (7.5,-3)..(8,-3);
    \foreach \Point in {(6,-3),(8,-3)}
    \draw[fill=red] \Point circle (0.05);
    \draw (5.2,-2.5) node{\tiny$(a_1,b_1)$} (8.7,-2.5) node{\tiny$(a_2,b_2)$} (8.1,-3.7) node{\tiny$(3,1,\ldots,1)$} (4.6,-4.1) node{$(2)$ $a_1>a_2$, or $b_1>a_2$, or $b_1>b_2$.};
    \end{tikzpicture}
    \caption{$d$ is odd and $a_1,a_2$ are even. The points in $\im(F)$ are coloured in red.}
    \label{fig:3-local1}
\end{figure}

\begin{lemma}
\label{lem:3-local-brief-version}
Fix an odd integer $d\geq3$, and two partitions $(a_1,b_1)$ and $(a_2,b_2)$ of $d$.
Let $\varphi:C\to T$ be a genus $0$ real tropical cover of an open $3$-valent line $T$ with ramification
profiles $(a_1,b_1)$, $(a_2,b_2)$ and $(3,1,\ldots,1)$.
Then, all non-zero local Hurwitz numbers $H^\rb(\varphi,v)$ are listed as follows.
\begin{enumerate}[$(1)$]
  \item If $a_1,b_1,a_2,b_2$ satisfy the relation encoded in Figure $\ref{fig:3-local-brief}(1)$,
there is only one real covering in $\rl_0((a_1,b_1), (a_2,b_2))$ having non-trivial contribution to $H^\rb(\varphi,v)$,
and
$
H^\rb(\varphi,v)=\alpha(\varphi)=2^{\frac{d-3}{2}}\cdot\left(\frac{d-3}{2}\right)!.
$
The involution $\iota$ and the image of $F$ are depicted in Figure $\ref{fig:3-local1}(1)$.
Moreover, the real simple resolution of the real ramified covering, which contributes to $H^\rb(\varphi,v)$,
defines a real tropical cover depicted in Figure $\ref{fig:3-local1-res}$.
\item If $a_1,b_1,a_2,b_2$ satisfy the relation encoded in Figure $\ref{fig:3-local-brief}(2)$, $(3)$ or $(4)$,
there exist $i(v)$ real coverings in $\rl_0((a_1,b_1), (a_2,b_2))$ that have non-trivial contribution to $H^\rb(\varphi,v)$.
Here, $i(v)$ denotes the number of monodromy graphs in Figure $\ref{fig:3-local-brief}(2)-(4)$ that encode the relation of $a_1,b_1,a_2,b_2$.
Moreover, we have
$
H^\rb(\varphi,v)=i(v)\cdot\alpha(\varphi)=i(v)\cdot2^{\frac{d-3}{2}}\cdot\left(\frac{d-5}{2}\right)!.
$
The involution $\iota$ and the image of $F$ are depicted in Figure $\ref{fig:3-local1}(2)$.
If $a_1,b_1,a_2,b_2$ satisfy the relation encoded in Figure $\ref{fig:3-local-brief}(2)$, $(3)$ or $(4)$,
the real simple resolution of a real ramified covering, which contributes non-trivially to $H^\rb(\varphi,v)$,
defines a real tropical cover depicted in Figure $\ref{fig:3-local1-res1}$, 
Figure $\ref{fig:3-local1-res2}$, or Figure $\ref{fig:3-local1-res3}$, respectively.
\end{enumerate}
\end{lemma}

\begin{figure}[ht]
    \centering
    \begin{tikzpicture}
    \draw (0,0.5)--(2,0)--(4,0.5);
    \draw (0.7,-1)--(2,0)--(0.5,-1);
    \draw (0,-0.2)--(2,0)--(4,-0.3);
    \draw (0,-1)--(2,0)--(0.3,-1);
    \draw (1,-1)--(2,0)--(1.4,-1);
    \draw (2.6,-1)--(2,0)--(3,-1);
    \draw[line width=0.9pt,dotted] (1.8,-0.5)--(2.2,-0.5);
    \draw[<->,bend right] (0.9,-1.2) to (1.4,-1.2);
    \draw[<->,bend right] (2.6,-1.2) to (3.1,-1.2);
    \draw (0.5,0.25) node{\tiny$a_1$} (3.3,0.45) node{\tiny$a_2$} (0.2,0) node{\tiny$b_1$} (3.5,-0.4) node{\tiny$b_2$} (0.5,-1.1) node{\tiny$1$} (1,-1.1) node{\tiny$1$} (0,-1.1) node{\tiny$3$} (0.3,-1.1) node{\tiny$1$} (0.7,-1.1) node{\tiny$1$} (1.4,-1.1) node{\tiny$1$} (2.6,-1.1) node{\tiny$1$} (3,-1.1) node{\tiny$1$};
    \foreach \Point in {(6,0),(8,0)}
    \draw \Point ellipse (0.2 and 0.5);
    \draw (7,-1.1)  ellipse(0.4 and 0.2);
    \foreach \Point in {(6,-0.5),(8,-0.5)}
    \draw[fill=red] \Point circle (0.05);
    \draw[bend right] (6,0.5) to (8,0.5);
    \draw (6,-0.5)..controls(6.5,-0.5) and (6.6,-0.5)..(6.6,-1.1);
    \draw (7.4,-1.1)..controls(7.4,-0.5) and (7.5,-0.5)..(8,-0.5);
    \draw (5.2,0) node{\tiny$(a_1,b_1)$} (8.7,0) node{\tiny$(a_2,b_2)$} (6,-0.7) node{\color{red}\tiny$2\times$} (8,-0.7) node{\color{red}\tiny$2\times$} (8.1,-1.2) node{\tiny$(3,1,\ldots,1)$} (4.6,-1.6) node{$(1)$ $a_1,b_1,a_2,b_2$ are even, and $j(v)\geq1$.};
    \draw (0,-2)--(2,-2.5)--(4,-2);
    \draw (0,-2.7)--(2,-2.5)--(4,-2.8);
    \draw (2,-2.5)--(0.7,-3.5);
    \draw (2,-2.5)--(0.2,-3.5);
    \draw (1,-3.5)--(2,-2.5)--(1.4,-3.5);
    \draw (2.6,-3.5)--(2,-2.5)--(3,-3.5);
    \draw[line width=0.9pt,dotted] (1.8,-3)--(2.2,-3);
    \draw[<->,bend right] (0.9,-3.7) to (1.4,-3.7);
    \draw[<->,bend right] (2.6,-3.7) to (3.1,-3.7);
    \draw (0.5,-2.25) node{\tiny$a_1$} (3.3,-2.05) node{\tiny$a_2$} (0.4,-2.5) node{\tiny$b_1$} (3.5,-2.9) node{\tiny$b_2$}  (0.2,-3.6) node{\tiny$3$} (0.7,-3.6) node{\tiny$1$} (1,-3.6) node{\tiny$1$} (1.4,-3.6) node{\tiny$1$} (2.6,-3.6) node{\tiny$1$} (3,-3.6) node{\tiny$1$};
    \foreach \Point in {(6,-2.5),(8,-2.5)}
    \draw \Point ellipse (0.2 and 0.5);
    \draw (7,-3.6)  ellipse(0.4 and 0.2);
    \draw[bend right] (6,-2) to (8,-2);
    \draw (6,-3)..controls(6.5,-3) and (6.6,-3)..(6.6,-3.6);
    \draw (7.4,-3.6)..controls(7.4,-3) and (7.5,-3)..(8,-3);
    \foreach \Point in {(6,-3),(6,-2)}
    \draw[fill=red] \Point circle (0.05);
    \draw (5.2,-2.5) node{\tiny$(a_1,b_1)$} (8.7,-2.5) node{\tiny$(a_2,b_2)$} (8.1,-3.7) node{\tiny$(3,1,\ldots,1)$} (4.6,-4.1) node{$(2)$ $a_1,b_1$ are even, $a_2,b_2$ are odd, and $j(v)\geq1$.};
    \draw (0,-5)--(2,-5.5)--(4,-5);
    \draw (0,-5.7)--(2,-5.5)--(4,-5.8);
    \draw (2,-5.5)--(0.7,-6.5);
    \draw (2,-5.5)--(0.2,-6.5);
    \draw (1,-6.5)--(2,-5.5)--(1.4,-6.5);
    \draw (2.6,-6.5)--(2,-5.5)--(3,-6.5);
    \draw[line width=0.9pt,dotted] (1.8,-6)--(2.2,-6);
    \draw[<->,bend right] (0.9,-6.7) to (1.4,-6.7);
    \draw[<->,bend right] (2.6,-6.7) to (3.1,-6.7);
    \draw (0.5,-5.25) node{\tiny$a_1$} (3.3,-5.05) node{\tiny$a_2$} (0.4,-5.5) node{\tiny$b_1$} (3.5,-5.9) node{\tiny$b_2$}  (0.2,-6.6) node{\tiny$3$} (0.7,-6.6) node{\tiny$1$} (1,-6.6) node{\tiny$1$} (1.4,-6.6) node{\tiny$1$} (2.6,-6.6) node{\tiny$1$} (3,-6.6) node{\tiny$1$};
    \foreach \Point in {(6,-5.5),(8,-5.5)}
    \draw \Point ellipse (0.2 and 0.5);
    \draw (7,-6.6)  ellipse(0.4 and 0.2);
    \draw[bend right] (6,-5) to (8,-5);
    \draw (6,-6)..controls(6.5,-6) and (6.6,-6)..(6.6,-6.6);
    \draw (7.4,-6.6)..controls(7.4,-6) and (7.5,-6)..(8,-6);
    \foreach \Point in {(8,-6),(8,-5)}
    \draw[fill=red] \Point circle (0.05);
    \draw (5.2,-5.5) node{\tiny$(a_1,b_1)$} (8.7,-5.5) node{\tiny$(a_2,b_2)$} (8.1,-6.7) node{\tiny$(3,1,\ldots,1)$} (4.6,-7.1) node{$(3)$ $a_1,b_1$ are odd, $a_2,b_2$ are even, and $j(v)\geq1$.};
    \draw (0,-8)--(2,-8.5)--(4,-8);
    \draw (0,-8.7)--(2,-8.5)--(4,-8.8);
    \draw (2,-8.5)--(0.7,-9.5);
    \draw (2,-8.5)--(0.2,-9.5);
    \draw (1,-9.5)--(2,-8.5)--(1.4,-9.5);
    \draw (2.6,-9.5)--(2,-8.5)--(3,-9.5);
    \draw[line width=0.9pt,dotted] (1.8,-9)--(2.2,-9);
    \draw[<->,bend right] (0.9,-9.7) to (1.4,-9.7);
    \draw[<->,bend right] (2.6,-9.7) to (3.1,-9.7);
    \draw (0.5,-8.25) node{\tiny$a_1$} (3.3,-8.05) node{\tiny$a_2$} (0.4,-8.5) node{\tiny$b_1$} (3.5,-8.9) node{\tiny$b_2$}  (0.2,-9.6) node{\tiny$3$} (0.7,-9.6) node{\tiny$1$} (1,-9.6) node{\tiny$1$} (1.4,-9.6) node{\tiny$1$} (2.6,-9.6) node{\tiny$1$} (3,-9.6) node{\tiny$1$};
    \foreach \Point in {(6,-8.5),(8,-8.5)}
    \draw \Point ellipse (0.2 and 0.5);
    \draw (7,-9.6)  ellipse(0.4 and 0.2);
    \draw[bend right] (6,-8) to (8,-8);
    \draw (6,-9)..controls(6.5,-9) and (6.6,-9)..(6.6,-9.6);
    \draw (7.4,-9.6)..controls(7.4,-9) and (7.5,-9)..(8,-9);
    \draw (5.2,-8.5) node{\tiny$(a_1,b_1)$} (8.7,-8.5) node{\tiny$(a_2,b_2)$} (8.1,-9.7) node{\tiny$(3,1,\ldots,1)$} (4.6,-10.1) node{$(4)$ $a_1,b_1,a_2,b_2$ are odd, and $j(v)\geq1$.};
    \end{tikzpicture}
    \caption{$d$ is even and the points in $\im(F)$ are coloured in red.}
    \label{fig:4-local}
\end{figure}

\begin{lemma}
\label{lem:4-local-brief-version}
Fix an even integer $d\geq4$, and two partitions $(a_1,b_1)$ and $(a_2,b_2)$ of $d$.
Let $\varphi:C\to T$ be a genus $0$ real tropical cover of an open $3$-valent line $T$ with ramification
profiles $(a_1,b_1)$, $(a_2,b_2)$ and $(3,1,\ldots,1)$.
Let $j(v)$ denote the number of pairs, whose first entry is greater than the second,
within the sequence $(a_1,a_2),(a_1,b_2),(b_1,a_2),(b_1,b_2)$.
Then, the following statements hold.
\begin{enumerate}[$(1)$]
\item If $a_1\neq b_1$ and $a_2\neq b_2$, the number $|\rl_0((a_1,b_1),(a_2,b_2))|=j(v)$,
and every ramified covering in $\rl_0((a_1,b_1),(a_2,b_2))$ contributes to $H^\rb(\varphi,v)$ non-trivially.
\item If $a_1=b_1$ and $a_2\neq b_2$, or $a_1\neq b_1$ and $a_2=b_2$,  the number $|\rl_0((a_1,b_1),(a_2,b_2))|=1$.
\item All non-zero local Hurwitz numbers $H^\rb(\varphi,v)$ are listed as follows
(The map $F$ is encoded by marking the points in $\im(F)$ in red.
If fixed points of $s$ real even edges in $C$ are mapped to a same fixed point of an even end of $T$,
we use the symbol `` $s\times$" to indicate that, where $s=2,3$).
\begin{itemize}
\item If $a_1,b_1,a_2,b_2$ satisfy the relation encoded in Figure $\ref{fig:4-local-brief}(1)$,
$
H^\rb(\varphi,v)=j(v)\cdot\alpha(\varphi)=j(v)\cdot3!\cdot2^{\frac{d-6}{2}}\cdot\left(\frac{d-6}{2}\right)!.
$
The involution $\iota$ and the image of $F$ are depicted in Figure $\ref{fig:4-local}(1)$. 
The real simple resolution of a real ramified covering in $\rl_0((a_1,b_1),(a_2,b_2))$
defines a real tropical cover depicted in Figure $\ref{fig:3-local2-res}$.
\item If $a_1,b_1,a_2,b_2$ satisfy the relation encoded in Figure $\ref{fig:4-local-brief}(2)$, $(3)$ or $(4)$,
the number
$
H^\rb(\varphi,v)=j(v)\cdot\alpha(\varphi)=j(v)\cdot2^{\frac{d-4}{2}}\cdot\left(\frac{d-4}{2}\right)!.
$
The involution $\iota$ and the image of $F$ are depicted in Figure $\ref{fig:4-local}(2), (3)$ or $4$, respectively. 
The real simple resolution of the real ramified covering in $\rl_0((a_1,b_1),(a_2,b_2))$
defines a real tropical cover depicted in Figure $\ref{fig:3-local2-res1}$, 
Figure $\ref{fig:3-local2-res2}$, or Figure $\ref{fig:3-local2-res3}$,
respectively.
\end{itemize}
\end{enumerate}
\end{lemma}

\section{The tropical correspondence theorem}
\label{sec:4}

\subsection{Effectively coloured resolving tropical cover}

Fix two integers $d\geq1$, $g\geq0$. Let $\lambda$, $\mu$ be two partitions of $d$,
and let $s,t$ be two non-negative integers such that $2s+t=l(\lambda)+l(\mu)+2g-2$.
Let $T\pb^1=\rb\cup\{\pm\infty\}$ and we consider it as the \textit{tropical projective line}.
\begin{definition}
\label{def:ram-split}
Let $\undl x'=\{x'_1,\ldots,x'_{2s}\}\subset\rb$ be a set whose points satisfy $x'_1<\cdots<x'_{2s}$.
Denote by $x_i=(x'_{2i-1},x'_{2i})$ the pair consisting of two points $x'_{2i-1},x'_{2i}$, where $i=1,\ldots,s$.
Let $\undl x=\{x_1,\ldots,x_s\}$ be the set consisting of pairs, and let
$\undl y=\{y_1,\ldots,y_{t}\}\subset\rb\setminus\undl x'$ be a set of $t$ points with $y_1<\ldots<y_{t}$.
The set $\undl z=\undl x\sqcup\undl y$ is called an {\it $(s,t)$-distribution}, if
$[x'_{2i-1},x'_{2i}]\cap\undl y=\emptyset$ for any $i=1,\ldots,s$.
\end{definition}
Let $\undl z=\undl x\sqcup\undl y$ be an $(s,t)$-distribution. 
We use $x_i<y_j$ to denote that $x'_{2i-1}<x'_{2i}<y_j$.
With this notation, we may suppose that
the elements of $\undl z=\{z_1,\ldots,z_{s+t}\}$ satisfy $z_1<\cdots<z_{s+t}$.

\begin{definition}
\label{def:compatible2}
Let $\undl z=\undl x\sqcup\undl y$ be an $(s,t)$-distribution
whose elements satisfy $z_1<\cdots<z_{s+t}$,
and $\Lambda_{s,t}=(\Lambda_1,\ldots,\Lambda_{s+t})$ be an $(s,t)$-tuple.
The $(s,t)$-tuple $\Lambda_{s,t}$ is \textit{compatible} with the $(s,t)$-distribution $\undl z=\undl x\sqcup\undl y$,
if the following relation holds.
\begin{itemize}
    \item $\Lambda_i=(3,1,\ldots,1)$ iff $z_i\in\undl x$, and $\Lambda_j=(2,1,\ldots,1)$ iff $z_j\in\undl y$.
\end{itemize}
\end{definition}

\begin{definition}
\label{def:d-tro-cover}
Let $\undl z=\undl x\sqcup\undl y$ be an $(s,t)$-distribution.
A \textit{tropical cover} $\varphi:C\to T\pb^1$ of type $(g,\lambda,\mu,\undl z)$
consists of a tropical curve $C$ and a continuous map $\varphi$ such that
\begin{enumerate}
    \item $C$ is a genus $g$ tropical curve with only $1$-valent and $3$-valent vertices.
    \item $\undl x'\sqcup\undl y$ is the set of images of inner vertices of $C$ under $\varphi$.
    Vertices in $\undl x'\sqcup\undl y$ are called {\it inner vertices} of $T\pb^1$.
    \item The set of leaves of $C$ is mapped onto $\{\pm\infty\}$ by $\varphi$.
    \item $\varphi$ is piecewise integer affine linear,
    and the slope $\w(e)\in\nb_{>0}$ of $\varphi$ on an edge $e$ is called the {\it weight} of $e$.
    \item $\varphi$ satisfies the balancing condition at any inner vertex of $C$.
    \item $\lambda$ (resp. $\mu$) is the tuple of weights of ends whose images under $\varphi$
    are adjacent to $-\infty$ (resp. $+\infty$).
\end{enumerate}
\end{definition}

Let $\varphi:C\to T\pb^1$ be a tropical cover.
For any edge $e'$ of $T\pb^1$,
it follows from the balancing condition that the sum
$$
\deg(\varphi):=\sum_{
    \substack{e\text{ edge of } C\\
    e'\subset\varphi(e)}}\w(e)
$$
is independent of $e'$, and it is called the \textit{degree} of $\varphi$.

A \textit{symmetric cycle} (resp. \textit{symmetric fork}) of a
tropical cover $\varphi:C\to T\pb^1$ is a pair of inner edges (resp. ends) of the same weight
and adjacent to the same two vertices (resp. one vertex).
If the ends of a symmetric fork of 
$\varphi:C\to T\pb^1$ can be oriented pointing from leaves to inner vertices (resp. from inner vertices to leaves),
the symmetric fork is called an \textit{inward symmetric fork} (resp. \textit{outward symmetric fork}).
We denote by $\sym(\varphi)$ the set of symmetric cycles and symmetric forks of $\varphi:C\to T\pb^1$.

\begin{figure}[ht]
    \centering
    \begin{tikzpicture}
    \draw[line width=0.4mm] (0,-1)--(-1,-1)--(-3,0);
    \draw[line width=0.4mm] (-2,-1)--(-1,-1);
    \draw[line width=0.4mm] (-3,-1.5)--(-2,-1)--(-3,-0.5);
    \draw (-2.8,-0.75) node {\tiny $k$} (-2.8,-1.25) node {\tiny $k$} (-1.6,-1.2) node {\tiny $2k$} (-1.6,-0.5) node {\tiny $e$} (-0.6,-1.2) node {\tiny $e+2k$};
    \draw (-1.6,-2) node {(\textrm{i})};
    \draw[line width=0.4mm] (3.5,-1)--(2.5,-1)--(0.5,0);
    \draw[line width=0.4mm] (1.5,-1)--(2.5,-1);
    \draw[line width=0.4mm] (0.5,-1.5)--(1.5,-1)--(0.5,-0.5);
    \draw (0.7,-0.75) node {\tiny $k$} (0.7,-1.25) node {\tiny $k$} (1.9,-1.2) node {\tiny $2k$} (1.9,-0.5) node {\tiny $o$} (2.9,-1.2) node {\tiny $o+2k$};
    \draw (1.9,-2) node {(\textrm{ii}) };
    \draw[line width=0.4mm] (4.3,-1)--(5.3,-1)--(7.3,0);
    \draw[line width=0.4mm] (5.3,-1)--(6.3,-1);
    \draw[line width=0.4mm] (7.3,-1.5)--(6.3,-1)--(7.3,-0.5);
    \draw (7,-0.8) node {\tiny $k$} (7,-1.2) node {\tiny $k$} (5.9,-1.2) node {\tiny $2k$} (5.9,-0.5) node {\tiny $e$} (4.7,-1.2) node {\tiny $e+2k$};
    \draw (5.7,-2) node {(\textrm{iii})};
    \draw[line width=0.4mm] (7.8,-1)--(8.8,-1)--(10.8,0);
    \draw[line width=0.4mm] (8.8,-1)--(9.8,-1);
    \draw[line width=0.4mm] (10.8,-1.5)--(9.8,-1)--(10.8,-0.5);
    \draw (10.5,-0.8) node {\tiny $k$} (10.5,-1.2) node {\tiny $k$} (9.4,-1.2) node {\tiny $2k$} (9.4,-0.5) node {\tiny $o$} (8.2,-1.2) node {\tiny $o+2k$};
    \draw (9.2,-2) node {(\textrm{iv})};
    \draw[line width=0.4mm] (-3,-3.5)--(-2,-3.5)--(0,-2.7);
    \draw[line width=0.4mm] (-2,-3.5)--(-1,-4);
    \draw[line width=0.4mm] (-3,-4.3)--(-1,-4)--(0,-4);
    \draw (-2.8,-4.1) node{\tiny $e$} (-2.8,-3.3) node{\tiny $e'$}
    (-1.5,-3.1) node{\tiny $e_1'$} (-1.5,-3.6) node{\tiny $e_2'$} (-0.5,-3.85) node{\tiny $e_2'+e$}
    (-1.6,-4.8) node {(\textrm{v})};
    \draw[line width=0.4mm] (0.5,-3.5)--(1.5,-3.5);
    \draw[line width=0.4mm] (3.5,-2.7)--(1.5,-3.5)--(2.5,-4)--(3.5,-4);
    \draw[line width=0.4mm] (0.5,-4.3)--(2.5,-4);
    \draw (0.7,-4.1) node{\tiny $e_2$} (0.7,-3.3) node{\tiny $e_1$}
    (2,-3.1) node{\tiny $o_1$} (2,-3.6) node{\tiny $o_2$} (3,-3.85) node{\tiny $o_2+e_2$};
    \draw (1.9,-4.8) node {(\textrm{vi})};
    \draw[line width=0.4mm] (4.3,-3.5)--(5.3,-3.5);
    \draw[line width=0.4mm] (7.3,-2.7)--(5.3,-3.5)--(6.3,-4)--(4.3,-4.3);
    \draw[line width=0.4mm] (6.3,-4)--(7.3,-4);
    \draw (4.5,-4.1) node{\tiny $o$} (4.5,-3.3) node{\tiny $e$}
    (5.8,-3.1) node{\tiny $o_1$} (5.8,-3.6) node{\tiny $o_2$} (6.8,-3.85) node{\tiny $o_2+o$};
    \draw (5.7,-4.8) node {(\textrm{vii}) };
    \draw[line width=0.4mm] (7.8,-3.5)--(8.8,-3.5)--(10.8,-2.7);
    \draw[line width=0.4mm] (8.8,-3.5)--(9.8,-4);
    \draw[line width=0.4mm] (7.8,-4.3)--(9.8,-4)--(10.8,-4);
    \draw (8,-4.1) node{\tiny $o$} (8,-3.3) node{\tiny $e$}
    (9.3,-3.1) node{\tiny $e_1$} (9.3,-3.6) node{\tiny $e_2$} (10.3,-3.85) node{\tiny $e_2+o$};
    \draw (9.2,-4.8) node {(\textrm{viii})};
    \draw[line width=0.4mm] (-3,-6.3)--(-2,-6.3);
    \draw[line width=0.4mm] (-2,-6.3)--(-1,-6.8);
    \draw[line width=0.4mm] (-2,-6.3)--(0,-5.5);
    \draw[line width=0.4mm] (-3,-7.1)--(-1,-6.8);
    \draw[line width=0.4mm] (0,-6.8)--(-1,-6.8);
    \draw (-2.8,-6.9) node{\tiny $e$} (-2.8,-6.1) node{\tiny $o$}
    (-1.5,-5.9) node{\tiny $e_1$} (-1.5,-6.4) node{\tiny $o_1$} (-0.5,-6.65) node{\tiny $o_1+e$};
    \draw (-1.6,-7.6) node {(\textrm{ix})};
    \draw[line width=0.4mm] (0.5,-6.3)--(1.5,-6.3);
    \draw[line width=0.4mm] (1.5,-6.3)--(2.5,-6.8);
    \draw[line width=0.4mm] (1.5,-6.3)--(3.5,-5.5);
    \draw[line width=0.4mm] (0.5,-7.1)--(2.5,-6.8);
    \draw[line width=0.4mm] (3.5,-6.8)--(2.5,-6.8);
    \draw (0.7,-6.9) node{\tiny $e$} (0.7,-6.1) node{\tiny $o$}
    (2,-5.9) node{\tiny $o_1$} (2,-6.4) node{\tiny $e_1$} (3,-6.65) node{\tiny $e_1+e$};
    \draw (1.9,-7.6) node {(\textrm{x})};
    \draw[line width=0.4mm] (4.3,-6.3)--(5.3,-6.3);
    \draw[line width=0.4mm] (5.3,-6.3)--(6.3,-6.8);
    \draw[line width=0.4mm] (5.3,-6.3)--(7.3,-5.5);
    \draw[line width=0.4mm] (4.3,-7.1)--(6.3,-6.8);
    \draw[line width=0.4mm] (7.3,-6.8)--(6.3,-6.8);
    \draw (4.5,-6.9) node{\tiny $o_2$} (4.5,-6.1) node{\tiny $o_1$}
    (5.8,-5.9) node{\tiny $e_1'$} (5.8,-6.4) node{\tiny $o_1'$} (6.8,-6.65) node{\tiny $o_1'+o_2$};
    \draw (5.7,-7.6) node {(\textrm{xi})};
    \draw[line width=0.4mm] (7.8,-6.3)--(8.8,-6.3);
    \draw[line width=0.4mm] (8.8,-6.3)--(9.8,-6.8);
    \draw[line width=0.4mm] (8.8,-6.3)--(10.8,-5.5);
    \draw[line width=0.4mm] (7.8,-7.1)--(9.8,-6.8);
    \draw[line width=0.4mm] (10.8,-6.8)--(9.8,-6.8);
    \draw (8,-6.9) node{\tiny $o_2$} (8,-6.1) node{\tiny $o_1$}
    (9.3,-5.9) node{\tiny $o_1'$} (9.3,-6.4) node{\tiny $e_1'$} (10.3,-6.65) node{\tiny $e_1'+o_2$};
    \draw (9.2,-7.6) node {(\textrm{xii})};
    \draw[line width=0.4mm] (-0.75,-8.5)--(-1.75,-8.5);
    \draw[line width=0.4mm] (0.25,-8.5)--(1.25,-8.5);
    \draw[line width=0.4mm] (-0.25,-8.5) +(0:0.5 and 0.25) arc (0:180:0.5 and 0.25);
    \draw[line width=0.4mm] (-0.25,-8.5) +(0:0.5 and 0.25) arc (0:-180:0.5 and 0.25);
    \draw (-1.25,-8.7) node {\tiny $o$} (-0.25,-8.1) node {\tiny $o_1$} (-0.25,-9) node {\tiny $e_1$} (0.65,-8.7) node {\tiny $o$};
    \draw (-0.25,-9.6) node {(\textrm{xiii})};
    \draw[line width=0.4mm] (7.05,-8.5)--(6.05,-8.5);
    \draw[line width=0.4mm] (8.05,-8.5)--(9.05,-8.5);
    \draw[line width=0.4mm] (7.55,-8.5) ellipse (0.5 and 0.25);
    \draw (6.45,-8.7) node {\tiny $e$} (7.55,-8.1) node {\tiny $e_1$} (7.55,-9) node {\tiny $e_2$} (8.45,-8.7) node {\tiny $e$};
    \draw (7.55,-9.6) node {(\textrm{xiv})};
    \end{tikzpicture}
    \caption{Pairs of vertices of resolving tropical covers.}
    \label{fig:pair-EC}
\end{figure}
\begin{definition}
\label{def:resolving-cover}
A tropical cover $\varphi:C\to T\pb^1$ of type $(g,\lambda,\mu,\undl z)$
is a \textit{resolving tropical cover} if
$\varphi:C\to T\pb^1$ satisfies the following conditions.
\begin{enumerate}
    \item For any $i\in\{1,\ldots,s\}$, the pair of vertices $x_i=(x_{2i-1}', x_{2i}')$
    are images of the adjacent two $3$-valent vertices of the figures in Figure \ref{fig:pair-EC}, up to reflection along a horizontal line.
    \item The pair of two flags with the same weight $k$ in the first row of Figure \ref{fig:pair-EC}
    is a symmetric fork of $C$ or is contained in a symmetric cycle of $C$.
\end{enumerate}
\end{definition}

Let $\varphi$ be a resolving tropical cover of type $(g,\lambda,\mu,\undl z)$.
The edge connecting the two vertices in any figure of Figure \ref{fig:pair-EC}
is called the \textit{contractible edge} of that pair of vertices.
A symmetric cycle or fork is called {\it non-contractible} if it does not contain any contractible edge.
We use the following notations.
\begin{itemize}
    \item $\symc(\varphi)$ is the set of symmetric cycles of $\varphi:C\to T\pb^1$.
    \item $E_c(\varphi)$ is the set of even contractible edges in $C$.
    \item $\sym_3(\varphi)$ (resp. $\sym_2(\varphi)$) is the set of non-contractible symmetric cycles and non-contractible symmetric forks which are adjacent
to a contractible edge (resp. are not adjacent to any contractible edge) .
    \item $\symc_c(\varphi)$ is the set of contractible symmetric cycles.
    \item $\nsym_c(\varphi)$ is the set of nonsymmetric
    cycles consisting of two even contractible edges.
\end{itemize}
Note that $\sym_3(\varphi)$ is the set of symmetric forks and symmetric cycles containing
the pair of two flags with the same weight $k$ in the first row of Figure \ref{fig:pair-EC}.

\begin{definition}
\label{def:coloured-TC}
A \textit{coloured} tropical cover $(\varphi:C\to T\pb^1,\rho)$ consists of a tropical cover and
a subset $I_\rho(\varphi)\subset \sym(\varphi)$ together with a
choice of colour red or blue for every component
of the subgraph of edges with even weights in $C\setminus (I_\rho(\varphi))^\circ$.
\end{definition}

Let $(\varphi_1:C_1\to T\pb^1;\rho_1)$ and
$(\varphi_2:C_2\to T\pb^1;\rho_2)$ be two coloured tropical covers.
An isomorphism $\varPhi:C_1\to C_2$ of tropical curves
preserving the colourings such that $\varphi_2\circ\varPhi=\varphi_1$
is called an {\it isomorphism} of two coloured tropical covers
$(\varphi_1,\rho_1)$ and $(\varphi_2,\rho_2)$.
A {\it positive} or {\it negative} point of
a coloured tropical cover $(\varphi,\rho)$ of type
$(g,\lambda,\mu,\undl z)$ is
the image of an inner vertex of $C$
depicted in the left two columns of Figure \ref{fig:coloured-vertices} or the right two columns of
Figure \ref{fig:coloured-vertices},
respectively, up to reflection along a vertical line.
\begin{figure}[ht]
    \centering
    \begin{tikzpicture}
    \draw[line width=0.3mm] (-3,0)--(-2,0)--(-1,0.5);
    \draw[line width=0.3mm,blue] (-2,0)--(-1,-0.5);
    \draw[line width=0.3mm,blue] (-0.5,0,0)--(0.5,0)--(1.5,0.5);
    \draw[line width=0.3mm,blue] (0.5,0)--(1.5,-0.5);
    \draw[line width=0.3mm,red] (-3,-1.3)--(-2,-1.3);
    \draw[line width=0.3mm] (-1,-0.8)--(-2,-1.3)--(-1,-1.8);
    \draw[line width=0.3mm,blue] (-0.5,-1.3)--(0.5,-1.3);
    \draw[line width=0.3mm,dotted] (1.5,-0.8)--(0.5,-1.3)--(1.5,-1.8);
    \draw[line width=0.3mm] (4,0)--(5,0)--(6,0.5);
    \draw[line width=0.3mm,red] (5,0)--(6,-0.5);
    \draw[line width=0.3mm,red] (6.5,0)--(7.5,0)--(8.5,0.5);
    \draw[line width=0.3mm,red] (7.5,0)--(8.5,-0.5);
    \draw[line width=0.3mm,blue] (4,-1.3)--(5,-1.3);
    \draw[line width=0.3mm] (6,-0.8)--(5,-1.3)--(6,-1.8);
    \draw[line width=0.3mm,red] (6.5,-1.3)--(7.5,-1.3);
    \draw[line width=0.3mm,dotted] (8.5,-0.8)--(7.5,-1.3)--(8.5,-1.8);
    \end{tikzpicture}
    \caption{Positive and negative vertices: even edges are coloured, odd edges are in black, and edges in $I_\rho$ are dotted.}
    \label{fig:coloured-vertices}
\end{figure}

\begin{definition}
\label{def:E-real-TC1}
A colouring $\rho$ of a resolving tropical cover $\varphi:C\to T\pb^1$
is \textit{effective}, if the two vertices $x'_{2i-1},x'_{2i}$ in a pair $x_i$
have the same sign,
$\sym_3(\varphi)\subset I_\rho(\varphi)$ and $\symc_c(\varphi)\cap I_\rho(\varphi)=\emptyset$.
\end{definition}
Let $(\varphi:C\to T\pb^1,\rho)$ be an effectively coloured resolving tropical cover.
From Definition \ref{def:E-real-TC1}, any pair of vertices $x'_{2i-1}, x'_{2i}$,
$i=1,\ldots,s$, are images of
two adjacent $3$-valent vertices of $C$ depicted in Figure \ref{fig:enhanced-vertices1},
up to reflection along a vertical line.
\begin{figure}[ht]
    \centering
    \begin{tikzpicture}
    \draw[line width=0.4mm,blue] (-3,-1)--(-2,-1)--(0,0);
    \draw[line width=0.4mm,blue] (-2,-1)--(-1,-1);
    \draw[line width=0.4mm,dotted] (0,-1.5)--(-1,-1)--(0,-0.5);
    \draw[line width=0.4mm] (0.5,-1)--(1.5,-1)--(3.5,0);
    \draw[line width=0.4mm,blue] (1.5,-1)--(2.5,-1);
    \draw[line width=0.4mm,dotted] (3.5,-1.5)--(2.5,-1)--(3.5,-0.5);
    \draw[line width=0.4mm] (-3,-2.5)--(-2,-2.5);
    \draw[line width=0.4mm] (-1,-2.5)--(0,-2.5);
    \draw[line width=0.4mm] (-1.5,-2.5) +(0:0.5 and 0.25) arc (0:180:0.5 and 0.25);
    \draw[line width=0.4mm,blue] (-1.5,-2.5) +(0:0.5 and 0.25) arc (0:-180:0.5 and 0.25);
    \draw[line width=0.4mm,blue] (0.5,-2.5)--(1.5,-2.5);
    \draw[line width=0.4mm,blue] (2.5,-2.5)--(3.5,-2.5);
    \draw[line width=0.4mm,blue] (2,-2.5) ellipse (0.5 and 0.25);
    \draw[line width=0.4mm,blue] (-3,-4)--(-2,-4)--(0,-3.2);
    \draw[line width=0.4mm,blue] (-2,-4)--(-1,-4.5);
    \draw[line width=0.4mm,blue] (-3,-4.8)--(-1,-4.5)--(0,-4.5);
    \draw[line width=0.4mm,red] (0.5,-4)--(1.5,-4);
    \draw[line width=0.4mm] (3.5,-3.2)--(1.5,-4)--(2.5,-4.5)--(3.5,-4.5);
    \draw[line width=0.4mm,blue] (0.5,-4.8)--(2.5,-4.5);
    \draw[line width=0.4mm,red] (-3,-6)--(-2,-6);
    \draw[line width=0.4mm] (0,-5.2)--(-2,-6)--(-1,-6.5)--(-3,-6.8);
    \draw[line width=0.4mm,red] (-1,-6.5)--(0,-6.5);
    \draw[line width=0.4mm,blue] (0.5,-6)--(1.5,-6);
    \draw[line width=0.4mm,blue] (1.5,-6)--(2.5,-6.5);
    \draw[line width=0.4mm,blue] (1.5,-6)--(3.5,-5.2);
    \draw[line width=0.4mm] (0.5,-6.8)--(2.5,-6.5);
    \draw[line width=0.4mm] (3.5,-6.5)--(2.5,-6.5);
    \draw[line width=0.4mm] (-3,-8)--(-2,-8);
    \draw[line width=0.4mm] (-2,-8)--(-1,-8.5);
    \draw[line width=0.4mm,blue] (-2,-8)--(0,-7.2);
    \draw[line width=0.4mm,blue] (-3,-8.8)--(-1,-8.5);
    \draw[line width=0.4mm] (0,-8.5)--(-1,-8.5);
    \draw[line width=0.4mm] (0.5,-8)--(1.5,-8);
    \draw[line width=0.4mm,blue] (1.5,-8)--(2.5,-8.5);
    \draw[line width=0.4mm] (1.5,-8)--(3.5,-7.2);
    \draw[line width=0.4mm] (0.5,-8.8)--(2.5,-8.5);
    \draw[line width=0.4mm] (3.5,-8.5)--(2.5,-8.5);
    \draw (0,-9.2) node{\tiny $(1)$ Positive pairs};
    \draw[dotted] (3.9,0)--(3.9,-9);
    \draw[line width=0.4mm,red] (4.3,-1)--(5.3,-1)--(7.3,0);
    \draw[line width=0.4mm,red] (5.3,-1)--(6.3,-1);
    \draw[line width=0.4mm,dotted] (7.3,-1.5)--(6.3,-1)--(7.3,-0.5);
    \draw[line width=0.4mm] (7.8,-1)--(8.8,-1)--(10.8,0);
    \draw[line width=0.4mm,red] (8.8,-1)--(9.8,-1);
    \draw[line width=0.4mm,dotted] (10.8,-1.5)--(9.8,-1)--(10.8,-0.5);
    \draw[line width=0.4mm] (4.3,-2.5)--(5.3,-2.5);
    \draw[line width=0.4mm] (6.3,-2.5)--(7.3,-2.5);
    \draw[line width=0.4mm] (5.8,-2.5) +(0:0.5 and 0.25) arc (0:180:0.5 and 0.25);
    \draw[line width=0.4mm,red] (5.8,-2.5) +(0:0.5 and 0.25) arc (0:-180:0.5 and 0.25);
    \draw[line width=0.4mm,red] (7.8,-2.5)--(8.8,-2.5);
    \draw[line width=0.4mm,red] (9.8,-2.5)--(10.8,-2.5);
    \draw[line width=0.4mm,red] (9.3,-2.5) ellipse (0.5 and 0.25);
    \draw[line width=0.4mm,red] (4.3,-4)--(5.3,-4)--(7.3,-3.2);
    \draw[line width=0.4mm,red] (5.3,-4)--(6.3,-4.5);
    \draw[line width=0.4mm,red] (4.3,-4.8)--(6.3,-4.5)--(7.3,-4.5);
    \draw[line width=0.4mm,blue] (7.8,-4)--(8.8,-4);
    \draw[line width=0.4mm] (10.8,-3.2)--(8.8,-4)--(9.8,-4.5)--(10.8,-4.5);
    \draw[line width=0.4mm,red] (7.8,-4.8)--(9.8,-4.5);
    \draw[line width=0.4mm,blue] (4.3,-6)--(5.3,-6);
    \draw[line width=0.4mm] (7.3,-5.2)--(5.3,-6)--(6.3,-6.5)--(4.3,-6.8);
    \draw[line width=0.4mm,blue] (6.3,-6.5)--(7.3,-6.5);
    \draw[line width=0.4mm,red] (7.8,-6)--(8.8,-6);
    \draw[line width=0.4mm,red] (8.8,-6)--(9.8,-6.5);
    \draw[line width=0.4mm,red] (8.8,-6)--(10.8,-5.2);
    \draw[line width=0.4mm] (7.8,-6.8)--(9.8,-6.5);
    \draw[line width=0.4mm] (10.8,-6.5)--(9.8,-6.5);
    \draw[line width=0.4mm] (4.3,-8)--(5.3,-8);
    \draw[line width=0.4mm] (5.3,-8)--(6.3,-8.5);
    \draw[line width=0.4mm,red] (5.3,-8)--(7.3,-7.2);
    \draw[line width=0.4mm,red] (4.3,-8.8)--(6.3,-8.5);
    \draw[line width=0.4mm] (7.3,-8.5)--(6.3,-8.5);
    \draw[line width=0.4mm] (7.8,-8)--(8.8,-8);
    \draw[line width=0.4mm,red] (8.8,-8)--(9.8,-8.5);
    \draw[line width=0.4mm] (8.8,-8)--(10.8,-7.2);
    \draw[line width=0.4mm] (7.8,-8.8)--(9.8,-8.5);
    \draw[line width=0.4mm] (10.8,-8.5)--(9.8,-8.5);
    \draw (7.3,-9.2) node{\tiny $(2)$ Negative pairs};
    \end{tikzpicture}
    \caption{Signed pairs: even edges are drawn in colours, odd edges in black. Dotted edges are
    the symmetric cycles or forks contained in $I_\rho$.}
    \label{fig:enhanced-vertices1}
\end{figure}
Denote by $\undl z^-$ (resp. $\undl z^+$) the set of negative (resp. positive) points or pairs of $(\varphi,\rho)$.
Then we have a natural signed splitting of $\undl z=\undl x\sqcup\undl y=\undl z^-\sqcup\undl z^+$
into positive and negative points or pairs.

We define the multiplicity of an effectively coloured
resolving tropical cover $(\varphi,\rho)$ as
\begin{equation}
\label{eq:mult-1}
\mult^\rb(\varphi,\rho)=\frac{2^{|E(I_\rho)\setminus E_c(\varphi)|+|\symc(\varphi)\cap\sym_3(\varphi)|+|\nsym_c(\varphi)|}}{2^{|\sym_2(\varphi)|}}\prod_{e\in \symc(\varphi)\cap I_\rho}\omega(e),
\end{equation}
where $\omega(e)$ is the weight of an edge in the dotted symmetric cycle $e$,
and $E(I_\rho)$ is the set of inner edges with even weights in $C\setminus I_\rho^\circ$.

%
%

\begin{definition}
\label{def:compatible-signed}
Let $\Lambda_{s,t}=(\Lambda_1,\ldots,\Lambda_{s+t})$ be an $(s,t)$-tuple,
and $\Lambda_{s,t}=(\Lambda_{s,t}^-,\Lambda_{s,t}^+)$ be a signed splitting of $\Lambda_{s,t}$.
Suppose that $\undl z=\undl x\sqcup\undl y$ is an $(s,t)$-distribution,
where $\undl x=\{x_1,\ldots,x_{s}\}$ is a set consisting of $s$ pairs, and $\undl y=\{y_1,\ldots,y_t\}$ is a set consisting of $t$ points.
Let $\undl z^-=\{z^-_1,\ldots,z^-_{|\undl z^-|}\}$ and $\undl z^+=\{z^+_1,\ldots,z^+_{|\undl z^+|}\}$
be a signed splitting $\undl z^-\sqcup\undl z^+$ of $\undl z$.
The signed splitting $\undl z=\undl z^-\sqcup\undl z^+$
is \textit{compatible} with the signed splitting $(\Lambda_{s,t}^-,\Lambda_{s,t}^+)$ of $\Lambda_{s,t}$,
if $|\undl z^-|=|\Lambda^-_{s,t}|$ and the following correspondence relation holds.
\begin{itemize}
    \item $\Lambda_i^*=(3,1,\ldots,1)$ iff $z_i^*\in\undl x$ ({\it i.e.}, $z_i^*$ is a pair),
     and $\Lambda_j^*=(2,1,\ldots,1)$ iff $z_j^*\in\undl y$ ({\it i.e.}, $z_j^*$ is a point), where $*=\pm$.
\end{itemize}
\end{definition}

\subsection{The modified correspondence theorem}
We use Markwig--Rau's theorem \cite[Theorem 4.1]{mr-2015} and the real simple resolution
to compute real double Hurwitz numbers with triple ramification.

Let $\undl z=\{z_1,\ldots,z_{s+t}\}\subset T\pb^1\setminus\{\pm\infty\}$ be a set whose elements satisfy
$z_1<\cdots<z_{s+t}$.
Each $z_i$ is equipped with a sign $s(z_i)\in\{+1,-1\}$ and a ramification index $r(z_i)\in\{2,3\}$
such that $|\{z_i|r(z_i)=3\}|=s$.
We apply the construction described on page 513 of \cite{mr-2015} to obtain a real tropical curve $D$, and sketch the construction as follows.
For each $z_i\in\undl z$, we attach an end $l_i$ to $T\pb^1$ at the vertex $z_i$.
Consequently, we obtain a genus $0$ tropical curve $D$.
We continue to use $z_i$ to denote the inner vertex of $D$.
The identity map on $D$ is a real involution on $D$.
Let $f_i^{-\infty}$ and $f_i^{+\infty}$ be the flags
pointing from $z_i$ to $-\infty$ and $+\infty$, respectively.
Let $b_i\in\rb\pb^1\setminus\{0,\infty\}$ be a point whose sign is the same as $s(z_i)$.
Denote by $S_{z_i}$ the punctured surface $S_{0,3}$
obtained by removing three small open real discs centered at $0,b_i,\infty$ from $\cb\pb^1$.
Equip $S_{z_i}$ with the standard real structure of $\cb\pb^1$.
We assign to every vertex $z_i$ the real oriented surface $S_{z_i}$.
The circle $B_{f_i^{-\infty}}$ centered at $0$ (resp. $B_{f_i^{+\infty}}$ centered at $\infty$) in $S_{z_i}$
corresponds to the flag of $D$ pointing from $z_i$ to $-\infty$ (resp. $+\infty$).
Let $e\subset D$ be an inner edge given by two flags $f$ and $f'$.
We identify the positive fixed point on $B_f$ with the positive fixed point on $B_{f'}$
and the negative fixed  point on $B_f$ with the negative fixed point on $B_{f'}$.
This gives the tropical curve $D$ a real structure $(D,\id,(S_{z_i})_{z_i\in\undl z})$.
Let $(\sal,\iota_\sal)$ be the real topological surface obtained from the real tropical curve
$(D,\id,(S_{z_i})_{z_i\in\undl z})$ according to \cite[Remark 3.3]{mr-2015}.
It is a punctured surface obtained from $\cb\pb^1$ with the standard real structure
and punctures at $0,\infty$, and at $|\{z_i|s(z_i)=+1\}|$ positive points
and $|\{z_i|s(z_i)=-1\}|$ negative points.
We still use $b_i$ to denote the point where $\sal$ is punctured, corresponding to $S_{z_i}$.
The natural order on the real parts of $\sal$ induces a total order on the set
$\undl b=\{b_1,\ldots,b_{s+t}\}\subset\rb\pb^1\setminus\{0,\infty\}$.
We suppose that $b_{i_1}<\cdots<b_{i_{s+t}}$, where
$i_1,\ldots,i_{s+t}$ is a permutation of $1,\ldots,s+t$.
The ramification index of a point $b_i\in\undl b$, which corresponds to a vertex $z_i$, is defined to be $r(z_i)$.
The set $\undl b$ is an $(s,t)$-branch set.
Let $\Lambda_{s,t}=(\Lambda_1,\ldots,\Lambda_{s+t})$ be an $(s,t)$-tuple
compatible with the $(s,t)$-branch set $\undl b$,
and $(\Lambda^-_{s,t}, \Lambda^+_{s,t})$ be the signed splitting compatible with $\undl b$.
Let $\xi=(\lambda,\Lambda_{s,t},\mu)$ be a collection of $s+t+2$ partitions of $d$.
The partition $\Lambda_i$ is fixed on the leaf of $l_{z_i}$, and $\lambda,\mu$
are fixed on $-\infty,+\infty$ respectively.
From Theorem \ref{thm:markwig-rau},
$H^\rb_g(\lambda,\mu;\Lambda^-_{s,t}, \Lambda^+_{s,t})=H^{\trop}_g((D,\id,(S_{z_i})_{z_i\in\undl z}),\xi)$.

\begin{construction}\label{def:simple-resolution}
Let $\varphi:C\to D$ be a real tropical cover contributing to the tropical real Hurwitz number
$H^{\trop}_g((D,\id,(S_{z_i})_{z_i\in\undl z}),\xi)$.
A {\it simple resolution} $\tilde\varphi:\widetilde C\to\widetilde D$ of $\varphi:C\to D$ is a real tropical cover
obtained from $\varphi:C\to D$
after a series of surgeries as follows.

Let $z\in\undl z$ be a vertex with ramification index $r(z)=3$.
We first consider the case $s(z)=+1$ ({\it i.e.}, $l_{z}$ points down from $z$).
Cut the two edges in $D$ that are adjacent to $z$ and different from $l_{z}$.
Let $T$ be the open three-valent line containing $z$.
We have a real tropical curve $(T,\id,S_{z})$,
and the complement $(D\setminus T,\id,(S_{z_i})_{z_i\in\undl z\setminus\{z\}})$.
Suppose that $\varphi^{-1}(z)=\{v_1,\ldots,v_k\}$.
Cut the tropical curve $(C,\iota_C,(S_v)_{v\in\vt(C)})$ into two parts.
The first part consists of the collection $(C_{v_i},\iota_{C_{v_i}}, S_{v_i})_{i=1}^k$,
where $C_{v_1},\ldots,C_{v_k}$ are the preimages $\varphi^{-1}(T)$,
$\iota_{C_{v_i}}$ is the involution restricted to $C_{v_i}$.
The second part is the complement, denoted as $(C\setminus\cup_iC_{v_i},\iota_{C\setminus\cup_iC_{v_i}},(S_v)_{v\in\vt(C)\setminus\cup_iv_i})$.

The restriction $\varphi_{v_i}:=\varphi|_{C_{v_i}}:C_{v_i}\to T$ defines a real tropical cover of degree $d_{\varphi}(v_i)$.
Among $\varphi_{v_1},\ldots,\varphi_{v_k}$, there are $k-1$ tropical covers
that ramify with ramification profiles $(d_{\varphi}(v_i)),(d_{\varphi}(v_i)),(1,\ldots,1)$.
Suppose $\varphi_{v_1},\ldots,\varphi_{v_{k-1}}$ are these $k-1$ real tropical covers.
It follows from \cite[Lemma 5.3]{mr-2015} that the real tropical covers $\varphi_{v_1},\ldots,\varphi_{v_{k-1}}$
are depicted in Figure \ref{fig:local-simple}.
Let $\tilde T$ be a genus 0 open tropical curve with only two inner vertices $z',z''$.
Moreover, $z',z''$ are three-valent vertices adjacent to two ends $l',l''$, respectively.
Then $(\tilde T,\id,(S_{z'},S_{z''}))$ is a real tropical curve. 
If $\varphi_{v_i}$ is depicted in Figure \ref{fig:local-simple}(1), (2) or (3),
let $\tilde\varphi_{v_i}:\tilde C_{v_i}\to \tilde T$ be the real tropical cover depicted
in Figure \ref{fig:local-simple1}, Figure \ref{fig:local-simple2} or Figure \ref{fig:local-simple3},
respectively.
Suppose that $\varphi_{v_i}':C_{v_i}'\to T$ is a real tropical cover isomorphic to $\varphi_{v_i}$ for $i=1,\ldots,k-1$.
Let $\tilde C_{v_i}$ be the real tropical curve obtained by gluing $C_{v_i}$ and $C_{v_i}'$
at a weighted $d_{\varphi}(v_i)$ leaf. 
In fact, $\tilde\varphi_{v_i}:\tilde C_{v_i}\to \tilde T$ is the real tropical cover resulting
from gluing $\varphi_{v_i}'$ and $\varphi_{v_i}$ at a weighted $d_{\varphi}(v_i)$ leaf.
This $\tilde\varphi_{v_i}$ is a genus 0 real tropical cover of $\tilde T$ with ramification profiles
$(d_{\varphi}(v_i)),(d_{\varphi}(v_i)),(1,\ldots,1), (1,\ldots,1)$.
The real tropical cover $\varphi_{v_k}$ is one of the type considered in Lemma \ref{lem:1-local}--\ref{lem:4-local-brief-version}.
Several real ramified coverings contribute to the local Hurwitz number $H^\rb(\varphi_{v_k},v_k)$ non-trivially.
Let $\tilde\varphi_{v_k}:\tilde C_{v_k}\to T\pb^1$ be the real tropical cover corresponding to a simple resolution
of such a real ramified covering.

Define $F_{\tilde C_{v_i}}$ as the set of flags in $\tilde C_{v_i}$ not mapped to $l'$ and $l''$ by $\tilde\varphi_{v_i}$.
Similarly, let $F_{C_{v_i}}$ be the set of flags in $C_{v_i}$ not mapped to $l_z$ by $\varphi_{v_i}$.
There exists a bijection between $F_{\tilde C_{v_i}}$ and $F_{C_{v_i}}$ that preserves weights and the fixed point set $F_e$.
Now, replace the collection
$C_{v_1},\ldots,C_{v_k}$ in $C$ with $\tilde C_{v_1},\ldots,\tilde C_{v_k}$, respectively.
Let $\tilde C_z$ be the resulting tropical curve,
and it is a real tropical curve.
Note that the two flags in $\tilde T$ share the same fixed point sets with the two flags in $T$. 
Replace $T$ in $D$ with $\tilde T$ to obtain a real tropical curve $\tilde D_z$.
The collection of real tropical covers $\tilde\varphi_{v_1},\ldots,\tilde\varphi_{v_k}$
and the restriction of $\varphi$ on
$(C\setminus\cup_iC_{v_i},\iota_{C\setminus\cup_iC_{v_i}},(S_v)_{v\in\vt(C)\setminus\cup_iv_i})$
induce a real tropical cover $\tilde\varphi_{z}:\tilde C_z\to\tilde D_z$.

If $r(z)=-1$ ({\it i.e.} $l_{z}$ points up from $z$), local real covers $\varphi_{v_i}$ and its resolutions $\tilde\varphi_{v_i}$ are
given by reflections of the local real covers and resolutions in the above case along a horizontal line.
By the same cut-glue process as above, we resolve $\varphi$ locally at $z$.
$\tilde\varphi:\widetilde C\to\widetilde D$ is the real tropical cover after resolving all $z\in\undl z$ with $r(z)=3$.
\end{construction}

Let $\varphi:C\to D$ be a real tropical cover
that has non-trivial contribution to the tropical real Hurwitz number
$H^{\trop}_g((D,\id,(S_{z_i})_{z_i\in\undl z}),\xi)$.
Note that $\undl z$ is a set of $s+t$ points, where each point $z_i$ has a sign index $s(z_i)$ and a ramification index $r(z_i)$.
Let $\tilde\varphi:\widetilde C\to\widetilde D$ be the simple resolution of $\varphi:C\to D$.
We use $\tilde{\undl z}$ to denote the inner vertices of $\widetilde D$ with the following meaning:
If $r(z_i)=3$, we use $\tilde z_i$ to denote the pair of vertices in $\tilde D_{z_i}$.
Otherwise, let $\tilde z_i$ be the point $z_i$ itself.
This gives an $(s,t)$-distribution $\tilde{\undl z}=\undl x\sqcup\undl y$,
and $\tilde{\undl z}$ is also equipped with a sign index $s(\tilde z_i)=s(z_i)$ and a ramification
index $r(\tilde z_i)=r(z_i)$.

\begin{theorem}
\label{thm:mr}
Let $g\geq0$, $d\geq1$, $s\geq0$ and $t\geq0$ be four integers,
and suppose that $\lambda$, $\mu$ are two partitions of $d$
such that $2s+t=l(\lambda)+l(\mu)+2g-2$.
Let $\tilde{\undl z}=\undl x\sqcup\undl y$ be an $(s,t)$-distribution.
Suppose that the set $\tilde{\undl z}=\undl x\sqcup\undl y$ possesses a splitting
$\tilde{\undl z}^-\sqcup\tilde{\undl z}^+$ compatible with a signed splitting $(\Lambda^-_{s,t}, \Lambda^+_{s,t})$
of an $(s,t)$-tuple $\Lambda_{s,t}$.
Then, we have
\begin{equation}\label{eq:mr-resolving}
H^\rb_g(\lambda,\mu;\Lambda^-_{s,t}, \Lambda^+_{s,t})=\sum_{[(\varphi,\rho)]}\mult^\rb(\varphi,\rho),
\end{equation}
where we sum over all isomorphism classes $[(\varphi,\rho)]$
of effectively coloured resolving tropical covers of type $(g,\lambda,\mu,\tilde{\undl z})$
whose positive and negative points or pairs of points reproduce
the splitting $\tilde{\undl z}^+\sqcup\tilde{\undl z}^-$.
\end{theorem}

\begin{proof}
The proof of this theorem modifies the strategy of the proof of \cite[Proposition 5.8]{mr-2015}.
The splitting $\tilde{\undl z}=\undl x\sqcup\undl y=\tilde{\undl z}^-\sqcup\tilde{\undl z}^+$
gives a sign index $s$ and a ramification index $r$ on $\tilde{\undl z}$ in an obvious way:
$s(\tilde z_i)=+1$ if $\tilde z_i\in\tilde{\undl z}^+$. Otherwise, $s(\tilde z_i)=-1$.
$r(\tilde z_i)=3$ if $\tilde z_i\in\undl x$ ({\it i.e.}, $\tilde z_i$ is a pair of points). Otherwise, $r(\tilde z_i)=2$.
Let $\undl z$ be a set consisting of $s+t$ points corresponding to $\tilde{\undl z}$, {\it i.e.}
if $\tilde z_i=(\tilde z_i',\tilde z_i'')$ is a pair of points, let $z_i=\tilde z_i'$;
if $\tilde z_i$ is a point, let $z_i=\tilde z_i$.
The sign index and ramification index on $\tilde{\undl z}$ induce a sign index and a ramification index on $\undl z$.
Let $D$ be the real tropical curve constructed above.
Let $\varphi:C\to D$ be a real tropical cover contributing to $H^{\trop}_g((D,\id,(S_{z_i})_{z_i\in\undl z}),\xi)$ with multiplicity $\mult^\rb(\varphi,F)$.
Denote by $\rl_\varphi$ the set of isomorphism classes of simple resolutions of $\varphi:C\to D$.
Let $\tilde\varphi:\widetilde C\to\widetilde D\in\rl_\varphi$ be a simple resolution of $\varphi$.
If $\tilde{z}_i$ is a pair, we use $l_{\tilde{z_i}}$ to denote the two ends in $\widetilde D$ adjacent to $\tilde z_i$.
Shrink the ends $l_{\tilde z_1},\ldots,l_{\tilde z_{s+t}}$ in $\widetilde D$ to the vertices (or pairs of vertices) $\tilde z_1,\ldots,\tilde z_{s+t}$,
and shrink the ends on $\widetilde C$ which are mapped to these ends.
We obtain a prereal tropical curve $\widetilde C'$ and a harmonic map $\widetilde\varphi':\widetilde C'\to T\pb^1$.
Note that the involution of $\widetilde C$ induces an involution on $\widetilde C'$.
We neglect the genus zero $2$-valent vertices in $\widetilde C'$ which are produced by shrinking ends,
and merge their two adjacent edges to one edge with proper length.
From the Riemann-Hurwitz condition, \cite[Lemma 5.3]{mr-2015}, \cite[Lemma 5.5]{mr-2015}, and
Lemma \ref{lem:1-local}--\ref{lem:4-local-brief-version}, the tropical curve $\widetilde C'$ is a
graph with $1$-valent vertices and $3$-valent vertices.

We choose the conjugated edges in the non-fixed locus of $\widetilde C'$
as the set $I_\rho(\widetilde\varphi')$ of the colouring $\rho$.
For a fixed even edge $e$ of $\widetilde C'$, we colour it in red if the fixed points of $e$ are mapped to the fixed point at the bottom of
the boundary circle, otherwise, we colour it in blue (see Figure \ref{fig:1-local-res}, Figure \ref{fig:2-local-res}, and Figure \ref{fig:3-local1-res}-\ref{fig:3-local2-res3}).
If an even fixed edge of $\widetilde C'$ emerges from merging two adjacent edges of $\widetilde C$,
it follows from \cite[Lemma 5.3]{mr-2015} that the two adjacent edges of $\widetilde C$ have the same fixed points.
Hence, this colouring is well-defined.
Note that $I_\rho(\widetilde\varphi')$ consists of $\sym_3(\widetilde\varphi')$
and part of $\sym_2(\widetilde\varphi')$, so $I_\rho(\widetilde\varphi')\cap\symc_c(\widetilde\varphi')=\emptyset$.
Moreover, according to the above colouring rule
a vertex $v$ in $\widetilde C'$ is coloured as one of the positive vertices in Figure \ref{fig:coloured-vertices} if
$\widetilde\varphi'(v)\in \tilde{\undl z}^+$,
otherwise, it is coloured as one of the negative vertices in Figure \ref{fig:coloured-vertices}.
Therefore, We obtain a set $\rl'_\varphi$ of isomorphism classes of effectively coloured resolving tropical covers
of type $(g,\lambda,\mu,\undl x\sqcup\undl y)$.
Moreover, positive and negative points or pairs of a cover in $\rl'_\varphi$ reproduce
a signed splitting $\tilde{\undl z}=\tilde{\undl z}^+\sqcup\tilde{\undl z}^-$
which is compatible with the signed splitting $(\Lambda^-_{s,t},\Lambda^+_{s,t})$.
Vice versa, we obtain a real tropical cover $(\tilde\varphi:\widetilde C\to\widetilde D)\in\rl_\varphi$
from an effectively coloured resolving tropical cover $(\tilde\varphi':\widetilde C'\to T\pb^1)\in\rl'_\varphi$
by growing ends at the inner vertices of $T\pb^1$ and their preimages in $\widetilde C'$.
From \cite[Lemma 5.3]{mr-2015}, \cite[Lemma 5.5]{mr-2015} and Lemma \ref{lem:1-local}--\ref{lem:4-local-brief-version},
we have a unique way to grow ends and extend the involution,
up to the choice of markings of the new ends.

Now we calculate the multiplicity.
Let $\varphi_v:C_v\to D_{\varphi(v)}$ be a local real cover at a real fixed vertex $v\in C$,
and $l_{\varphi(v)}$ be the end in $D_{\varphi(v)}$ adjacent to $\varphi(v)$. 
Denote by $\alpha(\varphi,v)$ the number of real automorphisms of $\varphi$
which only exchange the ends mapping to $l_{\varphi(v)}$.
For a non-fixed vertex $v$ of $C$, let $\alpha(\varphi,v)=d_\varphi(v)!$.
Let $V\subset\vt_{\fix}(C)$ be the set of inner vertices $v$ such that the 
local real cover $\varphi_v:C_v\to D_{\varphi(v)}$ at $v$ has ramification profiles $(a_1,b_1)$,
$(a_2,b_2)$ and $(3,1,\ldots,1)$, where $(a_1,b_1)$
and $(a_2,b_2)$ are two partitions of $d_\varphi(v)$.
The automorphism groups of $\varphi$ and $\widetilde\varphi'$ satisfy the following relation:
\begin{equation}
\label{eq:number-auto}
|\aut(\varphi)|=|\aut(\widetilde\varphi')|\prod_{[v]\in\vt(C)/\iota_C}\alpha(\varphi,v)\cdot2^{|\sym_V(C)|-|\symc_c(\widetilde\varphi')|},
\end{equation}
where $\sym_V(C)$ is the set of symmetric cycles and symmetric forks in $C$ that
are adjacent to a vertex in $V$, and that are not mapped to $l_{\varphi(v)}$ for any $v\in V$.
Note that the factor $2^{|\sym_V(C)|-|\symc_c(\widetilde\varphi')|}$ appears because
we resolve the $4$-valent vertices and genus one $2$-valent vertices in $C$.
We summarize the calculation of local Hurwitz numbers as follows.
\begin{enumerate}[$(1)$]
    \item Let $v\in C$ be a vertex with $r(\varphi(v))=2$, and suppose that $v$ produces a $3$-valent vertex in $\widetilde C'$ after shrinking the ends.
    We have $H^\rb(\varphi,v)=\alpha(\varphi,v)$ from \cite[Lemma 5.5]{mr-2015}.
    \item Let $v\in C$ be a vertex that produces a genus zero $2$-valent point $p_v$ in $\widetilde C'$ after shrinking the ends.
    From \cite[Lemma 5.3]{mr-2015}, we have $H^\rb(\varphi,v)=\alpha(\varphi,v)$ 
(resp. $H^\rb(\varphi,v)=\frac{1}{2}\alpha(\varphi,v)$), if the edge in $\widetilde C'$ that contains the point $p_v$ is a fixed edge of odd weight (resp. even weight).
    \item Let $v\in C$ be a non-real vertex that produces a genus zero $2$-valent point $p_v$ in $\widetilde C'$ after shrinking the ends.
    From \cite[Remark 5.4]{mr-2015}, we have $H^\cb(\varphi,v)=\frac{1}{\omega}\alpha(\varphi,v)$,
    where $\omega$ is the weight of the non-fixed edge in $\widetilde C'$ that contains the point $p_v$.
    \item Let $v\in C$ be a vertex with $r(\varphi(v))=3$,
    and suppose that $v$ produces a cycle $C_v$ in $\widetilde C'$ after resolution and shrinking the ends.
    From Lemma \ref{lem:1-local}, we have $H^\rb(\varphi,v)=\frac{d_\varphi(v)-1}{2}\alpha(\varphi,v)$ (resp. $H^\rb(\varphi,v)=\frac{d_\varphi(v)-2}{2}\alpha(\varphi,v)$)
    if the weight $d_\varphi(v)$ of the edge in $\widetilde C'$ adjacent to $C_v$ is odd (resp. even).
    In the case that $d_\varphi(v)$ is odd, there are $\frac{d_\varphi(v)-1}{2}$ ways to resolve $C$ locally at $v$
({\it i.e.} $v$ can produce $\frac{d_\varphi(v)-1}{2}$ contractible non-symmetric cycles (see Lemma \ref{lem:1-local} and Figure \ref{fig:1-local-res}(1))).
    Hence, $v$ contributes a factor $\frac{d_\varphi(v)-1}{2}$ to $|\rl_\varphi'|$.
    If $d_\varphi(v)$ is even, there are $\left\lceil\frac{d_\varphi(v)-2}{4}\right\rceil$ ways to resolve $C$ locally at $v$,
that is, $v$ can produce $\left\lceil\frac{d_\varphi(v)-2}{4}\right\rceil$ contractible cycles (see Lemma \ref{lem:1-local} and Figure \ref{fig:1-local-res}(2)).
    Hence, $v$ contributes a factor $\left\lceil\frac{d_\varphi(v)-2}{4}\right\rceil$ to $|\rl_\varphi'|$.
    \item Let $v\in C$ be a vertex with $r(\varphi(v))=3$, and suppose that $\varphi_v$ ramifies with profiles $(d_\varphi(v))$, $(a,b,c)$ and $(3,1,\ldots,1)$, where $(a,b,c)$ is a partition of $d_\varphi(v)$ with $b=c$.
After resolution and shrinking the ends, $v$ is resolved by one of the graphs in the first row of Figure \ref{fig:enhanced-vertices1}, 
up to reflection along a vertical line.
    From Lemma \ref{lem:2-local}, we have $H^\rb(\varphi,v)=2\alpha(\varphi,v)$.
    Since there is only one way to resolve $C$ locally at $v$, the vertex $v$ contributes a factor $1$ to $|\rl_\varphi'|$.
Each contractible edge in the graphs of the first row of Figure \ref{fig:enhanced-vertices1},
up to reflection along a vertical line, 
contributes a factor $2\alpha(\varphi,v)$ to the multiplicity.
    \item Let $v\in C$ be a vertex such that $\varphi_v$ ramifies with profiles $(a_1,b_1)$, $(a_2,b_2)$ and $(3,1,\ldots,1)$,
    where $(a_1,b_1)$ and $(a_2,b_2)$ are two partitions of $d_\varphi(v)$.
    The local Hurwitz number $H^\rb(\varphi,v)$ is computed in
    Lemma \ref{lem:3-local-brief-version} and Lemma \ref{lem:4-local-brief-version}.
    Let $H^\rb(\varphi,v)=f(v)\cdot\alpha(\varphi,v)$, where $f(v)$ is the factor
    given in those lemmas.
    When $d_\varphi(v)$ is even and $a_1=b_1$ or $a_2=b_2$, the number $f(v)=2$.
    Moreover, there is only one local real cover $\varphi_v:C_v\to D_{\varphi(v)}$ contributing to $H^\rb(\varphi,v)$ non-trivially
    and one simple resolution of $\varphi_v:C_v\to D_{\varphi(v)}$ (see Lemma \ref{lem:4-local-brief-version}).
    We denote by $f_1,f_2$ the two ends in $C_v$ weighted by $a_i,b_i$ with $a_i=b_i$, where $i=1$ or $2$.
The following cases describe how to glue the simple resolution $\tilde\varphi_v:\widetilde C_v\to\widetilde D_{\varphi(v)}$ 
with the remaining part $\varphi|_{C\setminus C_v}:C\setminus C_v\to D\setminus D_{\varphi(v)}$.
    \begin{enumerate}[$(a)$]
    \item If $f_1,f_2$ are not contained in any symmetric cycle or symmetric fork of $C$,
    we have two ways to glue the simple resolution
    $\tilde\varphi_v:\widetilde C_v\to\widetilde D_{\varphi(v)}$ with the remaining part
    $\varphi|_{C\setminus C_v}:C\setminus C_v\to D\setminus D_{\varphi(v)}$.
    \item If $f_1,f_2$ are contained in a symmetric fork of $C$, we only have one way to glue the simple resolution.
    \item Suppose that $f_1,f_2$ are contained in a symmetric cycle in $C$ with endpoints $v,v'$,
    and $r(\varphi(v'))=2$.
    We have one way to glue the simple resolution.
    \item Assume that $f_1,f_2$ are contained in a symmetric cycle in $C$ with endpoints $v,v'$,
    and $r(\varphi(v'))=3$.
    We have two ways to glue the simple resolution
    $\tilde\varphi_v:\widetilde C_v\to\widetilde D_{\varphi(v)}$
    with $\tilde\varphi_{v'}:\widetilde C_{v'}\to\widetilde D_{\varphi(v')}$.
    Let $\tilde\varphi_{v,v'}:\widetilde C_{v,v'}\to\widetilde D_{\varphi(v),\varphi(v')}$ be a glued cover.
   Then we have one way to glue $\tilde\varphi_{v,v'}:\widetilde C_{v,v'}\to\widetilde D_{\varphi(v),\varphi(v')}$ 
   with the remaining part.
    \end{enumerate}
    In the other cases, there are $f(v)$ local simple resolutions for $\varphi_v:C_v\to D_{\varphi(v)}$ (see Lemma \ref{lem:3-local-brief-version} 
    and Lemma \ref{lem:4-local-brief-version}).
    And every local simple resolution has a unique way to glue with the remaining part.
\end{enumerate}
Let $V_4'$ (resp. $V_4''$) be the set of type $(4)$ vertices mentioned above that have an odd local degree $d_\varphi(v)$ (resp. an even $d_\varphi(v)$).
Denote by $V_i$, $i=5,6$, the set of type $(i)$ vertices mentioned above.
Furthermore, let $V_6^{(1)}\subset V_6$ represent the subset of type $(b)$ vertices $v\in V_6$. 
Similarly, we use $V_6^{(2)}\subset V_6$ (resp. $V_6^{(3)}\subset V_6$) to denote the subset of type $(c)$ (resp. type $(d)$)
vertices $v\in V_6$.
Note that the factors in $(2)$ and $(3)$ cancel out with parts of
the factors $2^{|E(\varphi)|}, \prod\omega_\varphi(e),\frac{1}{|\aut(\varphi)|}$ in $\mult^\rb(\varphi,F)$ (see equation (\ref{eq:real-mult1})).
After shrinking all ends, we get
$$
\begin{aligned}
\mult^\rb(\varphi,F)=&\frac{2^{|E(I_\rho(\widetilde\varphi'))\setminus
E_c(\widetilde\varphi')|}}{|\aut(\widetilde\varphi')|\cdot2^{|\sym_V(C)|-|\symc_c(\widetilde\varphi')|}}
\cdot\prod_{e\in\symc(\tilde\varphi')\cap I_\rho(\widetilde\varphi')}\omega(e)\\
&\cdot\prod_{v\in V_4'}\frac{d_\varphi(v)-1}{2}\cdot\prod_{v\in V_4''}\frac{d_\varphi(v)-2}{2}
\cdot2^{|V_5|}\cdot\prod_{v\in V_6}f(v),
\end{aligned}
$$
where $\widetilde\varphi'$ is any cover in $\rl'_\varphi$.
The right hand part of the above equation does not depend on the choice of $\widetilde\varphi'$.
In the effectively coloured resolving tropical cover $\widetilde\varphi':\widetilde C\to T\pb^1$,
$$
\sym(\widetilde\varphi')=\sym_3(\widetilde\varphi')\cup\sym_2(\widetilde\varphi')\cup\symc_c(\widetilde\varphi'),
$$
so we have $|\aut(\widetilde\varphi')|\cdot2^{|\sym_V(C)|-|\symc_c(\widetilde\varphi')|}
=2^{|\sym_3(\widetilde\varphi')|+|\sym_2(\widetilde\varphi')|+|\sym_V(C)|}$.
The number $|V_5|=|\sym_3(\widetilde\varphi')|+|\sym_3(\widetilde\varphi')\cap\symc(\widetilde\varphi')|$.
From Lemma \ref{lem:3-local-brief-version}, Lemma \ref{lem:4-local-brief-version} and item $(6)$,
the product
$$
\prod_{v\in V_6}f(v)=\prod_{v\in V_6\setminus(V_6^{(1)}\cup V_6^{(2)}\cup V_6^{(3)})}f(v)\cdot2^{|V_6^{(1)}|+|V_6^{(2)}|+|V_6^{(3)}|}.
$$
The number $|\sym_V(C)|=|V_6^{(1)}|+|V_6^{(2)}|+\frac{1}{2}|V_6^{(3)}|$.
Hence, we get
$$
\begin{aligned}
\mult^\rb(\varphi,F)=&\frac{2^{|E(I_\rho(\widetilde\varphi'))\setminus
E_c(\widetilde\varphi')|+|\sym_3(\widetilde\varphi')\cap\symc(\widetilde\varphi')|}}{2^{|\sym_2(\widetilde\varphi')|}}
\cdot\prod_{e\in\symc(\tilde\varphi')\cap I_\rho(\widetilde\varphi')}\omega(e)\\
&\cdot\prod_{v\in V_4'}\frac{d_\varphi(v)-1}{2}\cdot\prod_{v\in V_4''}\frac{d_\varphi(v)-2}{2}
\cdot\prod_{v\in V_6\setminus(V_6^{(1)}\cup V_6^{(2)}\cup V_6^{(3)})}f(v)\cdot2^{\frac{1}{2}|V_6^{(3)}|},
\end{aligned}
$$

If $v,v'$ is a pair of vertices in $V_6^{(3)}$ connected by a symmetric cycle in $\varphi$,
we have two ways to resolve this symmetric cycle (see the item (6)).
If $v\in V_4''$, we have $\left\lceil\frac{d_\varphi(v)-2}{4}\right\rceil$ ways to resolve $C$ locally at $v$.
Denote by $\rl_v$ the set of local resolutions of a vertex $v$ in $V_4''$.
Let $c_1(v)$ be the number of pairs of integers $(2a_1(v),2b_1(v))$ such that $0<a_1(v)<b_1(v)$ and
$d_\varphi(v)=2a_1(v)+2b_1(v)$.
Let $c_2(v)$ be the number of pairs of integers $(2a_1(v),2a_1(v))$ such that $d_\varphi(v)=4a_1(v)$.
It is easy to see that $|\rl_v|=\left\lceil\frac{d_\varphi(v)-2}{4}\right\rceil=c_1(v)+c_2(v)$ and $\frac{d_\varphi(v)-2}{2}=2c_1(v)+c_2(v)$.
Note that 
$$
2c_1(v)+c_2(v)=\sum_{\widetilde\varphi_v\in\rl_v}2^{m(\widetilde\varphi_v)},
$$
where $m(\widetilde\varphi_v)=0$ if $\widetilde\varphi_v$ produces a symmetric cycle, otherwise let $m(\widetilde\varphi_v)=1$.
Then we have
$$
\begin{aligned}
\prod_{v\in V_4''}\frac{d_\varphi(v)-2}{2}&=\prod_{v\in V_4''}(\sum_{\widetilde\varphi_v\in\rl_v}2^{m(\widetilde\varphi_v)})\\
&=\sum_{x\in\prod_{v\in V_4''}\rl_v}2^{k(x)},
\end{aligned}
$$
where $k(x)$ is the number of elements in the sequence of local resolutions $x$ that yield non-symmetric cycles.
Note that $|\nsym_c(\widetilde\varphi')|$ is the number of local resolutions that contain non-symmetric cycles in $\widetilde\varphi'$.
Therefore, we obtain
$$
\mult^\rb(\varphi,F)=\sum_{[(\widetilde\varphi',\rho)]\in\rl_\varphi'}\frac{2^{|E(I_\rho(\widetilde\varphi'))\setminus
E_c(\widetilde\varphi')|+|\sym_3(\widetilde\varphi')\cap\symc(\widetilde\varphi')|+|\nsym_c(\widetilde\varphi')|}}{2^{|\sym_2(\widetilde\varphi')|}}
\cdot\prod_{e\in\symc(\tilde\varphi')\cap I_\rho(\widetilde\varphi')}\omega(e).
$$
Here, $|\rl_\varphi'|=\prod_{v\in V_4'}\frac{d_\varphi(v)-1}{2}\cdot\prod_{v\in V_4''}\left\lceil\frac{d_\varphi(v)-2}{4}\right\rceil
\cdot\prod_{v\in V_6\setminus(V_6^{(1)}\cup V_6^{(2)}\cup V_6^{(3)})}f(v)\cdot2^{\frac{1}{2}|V_6^{(3)}|}$.
Note that any effectively coloured resolving tropical cover $(\widetilde\varphi',\rho)$
can be obtained by resolving a real tropical cover contributing to
$H^{\trop}_g((D,\id,(S_{z_i})_{z_i\in\undl z}),\xi)$ and shrinking the ends.
If $(\widetilde\varphi',\rho)\notin\rl_\varphi'$, after growing ends and reversing the simple
resolution procedure, one can not get the real tropical cover $\varphi$.
Hence, we have $H^{\trop}_g((D,\id,(S_{z_i})_{z_i\in\undl z}),\xi)=\sum_{[(\widetilde\varphi',\rho)]}\mult^\rb(\widetilde\varphi',\rho)$,
where $\mult^\rb(\widetilde\varphi',\rho)$ is given in equation (\ref{eq:mult-1}).
\end{proof}

\begin{remark}
When $s=0$, $H^\rb_g(\lambda,\mu;\Lambda^-_{s,t}, \Lambda^+_{s,t})$ is the
ordinary real double Hurwitz number. Our computation in equation $(\ref{eq:mr-resolving})$
coincides with \cite[Corollary 5.9]{mr-2015} and \cite[Theorem 3.17]{gpmr-2015}.
When $t=0$, $H^\rb_g(\lambda,\mu;\Lambda^-_{s,t}, \Lambda^+_{s,t})$ is the
real double Hurwitz number with only $3$-cycles. Our computation in equation $(\ref{eq:mr-resolving})$
coincides with the pure combinatorial computation carried out
in a previous version \cite[Theorem 3.8]{dlly-2023} of this manuscript.
\end{remark}

\appendix
\section{Proofs of several lemmas}
\label{sec:a}
In this appendix, we give the proofs of Lemma \ref{lem:2-local}, Lemma \ref{lem:3-local-brief-version} and Lemma \ref{lem:4-local-brief-version}.

\begin{proof}[Proof of Lemma $\ref{lem:2-local}$]
We calculate
the number of tuples in $\fl^\rb_0((d),(a,b,c))$.
We have $(d-1)!$ possibilities to choose $\sigma_1$.
If $d$ is odd (resp. even), we have $d$ (resp. $\frac{d}{2}$) possibilities to choose $\gamma$.
Once partitions $(d)$ and $(a,b,c)$ are fixed, there is only one monodromy graph in Figure \ref{fig:2-local-brief} whose ends are weighted by the partitions $(d)$ and $(a,b,c)$.
From Lemma \ref{lem:computation-2}(3), there are two permutations (resp. one permutation) $\tau$
of type $(3,1,\ldots,1)$ such that $\gamma\circ(\tau\circ\sigma_1)\circ\gamma=(\tau\circ\sigma_1)^{-1}$ and
$\tau\circ\sigma_1$ is of type $(a,b,c)$, if $d$ is even (resp. odd).
Moreover, two entries in $(a,b,c)$ are equal (for example, $b=c$), and the two leaves in $C$ that are labelled by these two entries are conjugated leaves.
Note that the $d!$ tuples in $\fl^\rb_0((d),(a,b,c))$ are in a same orbit of the $S_d$ action, and the stabilizer of any tuple $(\sigma_1,\gamma,\tau,\sigma_2)\in\fl^\rb_0((d),(a,b,c))$ is trivial.
Therefore, there is only one real cover in $\rl_0((d),(a,b,c))$.
From Lemma \ref{lem:fixed-points}(2), the set $\pi^{-1}(p)$
consists of one real point and $\frac{d-3}{2}$ pairs (resp. two real points and $\frac{d-4}{2}$ pairs) of conjugated points, if $d$ is odd (resp. even).
Hence, when $d$ is odd (resp. even), the number of labels of the leaves in $C$ is
$2\alpha(\varphi)=2^{\frac{d-1}{2}}\cdot\left(\frac{d-3}{2}\right)!$ (resp. $2^{\frac{d-2}{2}}\cdot\left(\frac{d-4}{2}\right)!$),
where the factor $2$ comes from automorphisms of conjugated leaves weighted by
the two equal entries in $(a,b,c)$.

Let $A=(\sigma_1,\gamma,\tau,\sigma_2)$ be a tuple in $\fl^\rb_0((d),(a,b,c))$ with $b=c$, and $(\pi,\iota)=\Psi(A)$.
Suppose that $\tilde A=(\sigma_1,\gamma,\bar\tau_1,\bar\tau_2,\sigma_2)\in\widetilde\fl^\rb_0((d),(a,b,c))$ is the tuple corresponding to $A$ via Lemma \ref{lem:computation-2}(1).
From Lemma \ref{lem:computation-2}(2), the composition $\bar\tau_1\circ\sigma_1$ consists of two disjoint cycles with lengths $a,2b$.
A similar argument as the proof of Lemma \ref{lem:1-local} shows that $\tilde\pi^{-1}(B)$ consists of two circles $B_1$ and $B_2$.
Moreover, degrees of $\tilde\pi|_{B_1}$ and $\tilde\pi|_{B_2}$ are $a$, $2b$.
Riemann-Hurwitz formula implies that $\widetilde\cl$ is divided by $B_1$ and $B_2$ into three genus zero components.
It follows from \cite[Lemma 5.3]{mr-2015} and \cite[Lemma 5.5]{mr-2015} that the real tropical cover corresponding to a real ramified covering in $\widetilde\rl_0((d),(a,b,c))$ is depicted in Figure \ref{fig:2-local-res}.
\end{proof}

\begin{proof}[Proof of Lemma $\ref{lem:3-local-brief-version}$]
We only need to compute the number of tuples in $\fl^\rb_0((a_1,b_1),(a_2,b_2))$.
From Lemma \ref{lem:computation-3}(2), only if $(a_1,b_1)$ and $(a_2,b_2)$ satisfy
the relations depicted in Figure \ref{fig:3-local-brief}, the set $\fl^\rb_0((a_1,b_1), (a_2,b_2))\neq\emptyset$.
Suppose that $a_1,a_2$ are even.

There are $\frac{d!}{a_1\cdot b_1}$ permutations of cycle type $(a_1,b_1)$ in the symmetric group $S_d$,
so we have $\frac{d!}{a_1\cdot b_1}$ possibilities to choose $\sigma_1$.
Let $\sigma_1'$ and $\sigma_1''$ denote the circles that represent the cycles of lengths $a_1$ and $b_1$, respectively.
Let $\gamma\in S_d$ be an involution satisfying $\gamma\circ\sigma_1\circ\gamma=\sigma_1^{-1}$.
From \cite[Lemma $3.12$]{gpmr-2015} (see also Figure \ref{fig:sym-circle1}),
we have $b_1$ possibilities to choose the symmetric axis for the odd circle $\sigma_1''$.
For the even circle $\sigma_1'$, we have $\frac{a_1}{2}$ ways to choose the symmetric axis of the circle $\sigma_1'$.
From Lemma \ref{lem:computation-3}(3), for any relation in Figure \ref{fig:3-local-brief} there are $2$
permutations $\tau$ such that $(\sigma_1,\gamma,\tau,(\tau\circ\sigma_1)^{-1})\in\fl^\rb_0((a_1,b_1), (a_2,b_2))$.
Hence, every relation in Figure \ref{fig:3-local-brief} determines $d!$ tuples in $\fl^\rb_0((a_1,b_1), (a_2,b_2))$.
Note that these $d!$ tuples are in a same orbit of the action $S_d$ on $\fl^\rb_0((a_1,b_1), (a_2,b_2))$.
Since the stabilizer of any tuple $(\sigma_1,\gamma,\tau,\sigma_2)\in\fl^\rb_0((a_1,b_1),(a_2,b_2))$ is trivial,
these $d!$ tuples determine a real ramified covering in the set $\rl_0((a_1,b_1),(a_2,b_2))$ with trivial automorphism group.
Denote by $[(\pi_i,\iota_i)]$ the real ramified covering in $\rl_0((a_1,b_1), (a_2,b_2))$ determined by
the monodromy graph in Figure \ref{fig:3-local-brief}$(i)$, where $i=1,2,3,4$.
From Lemma \ref{lem:fixed-points}(2), the set $\pi_1^{-1}(p)$
consists of one real point and $\frac{d-3}{2}$ pairs of conjugated points (see Figure $\ref{fig:3-local1}(1)$).
The set $\pi_j^{-1}(p)$
consists of three real points and $\frac{d-5}{2}$ pairs of conjugated points (see Figure $\ref{fig:3-local1}(2)$),
where $j=2,3,4$. The image of $F$ is determined by Lemma \ref{lem:fixed-points}$(1)$.
There is only one real ramified covering in $\rl_0((a_1,b_1), (a_2,b_2))$ that
defines the real tropical cover depicted in Figure $\ref{fig:3-local1}(1)$.
The number $H^\rb(\varphi,v)$ for real tropical cover in Figure $\ref{fig:3-local1}(1)$
is $\alpha(\varphi)=2^{\frac{d-3}{2}}\cdot\left(\frac{d-3}{2}\right)!$.
There are $i(v)$ real ramified coverings in $\rl_0((a_1,b_1), (a_2,b_2))$ that 
define the real tropical cover depicted in Figure $\ref{fig:3-local1}(2)$.
Hence, the corresponding local Hurwitz number
$H^\rb(\varphi,v)=i(v)\cdot\alpha(\varphi)=i(v)\cdot2^{\frac{d-3}{2}}\cdot\left(\frac{d-5}{2}\right)!$.

Let $A=(\sigma_1,\gamma,\tau,\sigma_2)$ be a tuple in $\fl^\rb_0((a_1,b_1),(a_2,b_2))$, and $(\pi_1,\iota_1)=\Psi([A])$.
Assume that $a_1,b_1,a_2,b_2$ satisfy the relation encoded in Figure \ref{fig:3-local-brief}$(1)$.
Suppose that $\tilde A=(\sigma_1,\gamma,\bar\tau_1,\bar\tau_2,\sigma_2)\in\widetilde\fl^\rb_0((a_1,b_1),(a_2,b_2))$
is the tuple corresponding to $A$ via Lemma \ref{lem:computation-3}$(1)$.
Let $(\tilde\pi_1:\widetilde\cl_1\to\widetilde \dl,\tilde\iota)$ be the real simple resolution
corresponding to $\tilde A$.
From Lemma \ref{lem:computation-3}$(2)$,
the composition $\bar\tau_1\circ\sigma_1$ consists of three disjoint cycles with lengths $b_1,b_2,a_1-b_2$.
A similar argument as the proof of Lemma \ref{lem:1-local} shows that
$\tilde\pi_1^{-1}(B)$ consists of three circles $B_1$, $B_2$ and $B_3$,
where $B$ is the circle depicted in Figure \ref{fig:3-local1-res}.
Moreover, degrees of $\tilde\pi_1|_{B_1}$, $\tilde\pi_1|_{B_2}$ and $\tilde\pi_1|_{B_3}$ are $b_1$, $b_2$ and $a_1-b_2$.
Riemann-Hurwitz formula implies that $\widetilde\cl_1$ is divided by $B_1,B_2$ and $B_3$ into four genus zero components.
It follows from \cite[Lemma 5.3]{mr-2015} and \cite[Lemma 5.5]{mr-2015} that the real tropical cover corresponding to $\tilde\pi_1\in\widetilde\rl_0((a_1,b_1),(a_2,b_2))$ is depicted in Figure \ref{fig:3-local1-res}.
The construction of the real tropical covers corresponding to
$\tilde\pi_2$, $\tilde\pi_3$ and $\tilde\pi_4$ are similar to that of $\tilde\pi_1$,
so we omit the proof and only depict them in Figure \ref{fig:3-local1-res1}, Figure \ref{fig:3-local1-res2} and Figure \ref{fig:3-local1-res3},
respectively.
\end{proof}

\begin{proof}[Proof of Lemma $\ref{lem:4-local-brief-version}$]
We first consider the case that $a_1,b_1$ are even integers.
As in the proof of Lemma \ref{lem:3-local-brief-version},
we compute the number of tuples in $\fl^\rb_0((a_1,b_1),(a_2,b_2))$, then we find the number of real
coverings contributing to the local Hurwitz number non-trivially.
We have $\frac{d!}{a_1\cdot b_1}$ (resp. $\frac{d!}{2a_1\cdot b_1}$) possibilities
to choose $\sigma_1$, if $a_1\neq b_1$ (resp. if $a_1=b_1$).
Let $\sigma_1'$ and $\sigma_1''$ denote the circles that  represent the cycles of lengths $a_1$ and $b_1$, respectively.
From \cite[Lemma $3.12$]{gpmr-2015} (see also Figure \ref{fig:sym-circle1}),
for the even circles $\sigma_1'$ and $\sigma_1''$, we have $\frac{a_1b_1}{4}$ ways to choose the symmetric axes of $\sigma_1'$ and $\sigma_1''$.
\begin{enumerate}[$(a)$]
\item $a_1\neq b_1$, $a_2\neq b_2$. Without loss of generality, we assume $a_1>a_2$.
If $a_2,b_2$ are even (resp. odd), $j(v)$ equals the number of distinct ways to  assign $a_1,a_2,b_1,b_2$ to the ends of the monodromy graph in Figure \ref{fig:4-local-brief}(1) (resp. in Figure \ref{fig:4-local-brief}(2)).
Once $a_1,b_1,a_2,b_2$ are assigned to the ends of the monodromy graph in Figure \ref{fig:4-local-brief}$(1)$ or $(2)$, there is only one way to label $\sigma_1'$ and $\sigma_1''$ to the left ends of the monodromy graph.
From Lemma \ref{lem:computation-3}(3), we have four ways to choose the permutation $\tau$.
Hence, one way to
assign $a_1,b_1,a_2,b_2$ to the ends of the monodromy graph in Figure \ref{fig:4-local-brief}(1) or (2)
corresponds to $d!$ tuples in $\fl^\rb_0((a_1,b_1),(a_2,b_2))$.
Note that these $d!$ tuples are coordinatewise conjugate with each other.
Since the stabilizer of any tuple $(\sigma_1,\gamma,\tau,\sigma_2)\in\fl^\rb_0((a_1,b_1),(a_2,b_2))$ is trivial,
these $d!$ tuples correspond to a real cover $\pi$ in the set $\rl_0((a_1,b_1),(a_2,b_2))$ with trivial automorphism group.
If $a_2,b_2$ are even integers,
from Lemma \ref{lem:fixed-points}(2) and Figure \ref{fig:4-local-brief}(1), the set $\pi^{-1}(p)$
consists of four real points and $\frac{d-6}{2}$ pairs of conjugated points (see Figure $\ref{fig:4-local}(1)$). The image of $F$ is determined by Lemma \ref{lem:fixed-points}$(1)$.
Therefore, $H^\rb(\varphi,v)=j(v)\cdot\alpha(\varphi,v)=j(v)\cdot3!\cdot2^{\frac{d-6}{2}}(\frac{d-6}{2})!$,
and $|\rl_0((a_1,b_1),(a_2,b_2))|=j(v)$.
If $a_2,b_2$ are odd integers,
from Lemma \ref{lem:fixed-points}(2) and Figure \ref{fig:4-local-brief}(2), the set $\pi^{-1}(p)$
consists of two real points and $\frac{d-4}{2}$ pairs of conjugated points (see Figure $\ref{fig:4-local}(2)$). 
The image of $F$ is determined by Lemma \ref{lem:fixed-points}$(1)$.
Therefore, $H^\rb(\varphi,v)=j(v)\cdot\alpha(\varphi,v)=j(v)\cdot2^{\frac{d-4}{2}}(\frac{d-4}{2})!$,
and $|\rl_0((a_1,b_1),(a_2,b_2))|=j(v)$.
\item $a_1=b_1$ and $a_2\neq b_2$. Without loss of generality, we assume $a_1>a_2$.
We have one way to
assign $a_1,b_1,a_2,b_2$ to the ends of the monodromy graph in Figure \ref{fig:4-local-brief}(1) or (2).
Since $a_1=b_1$, we have two ways to assign $\sigma_1'$ and $\sigma_1''$ to the left ends
of the monodromy graph in Figure \ref{fig:4-local-brief}(1) or (2).
When $\sigma_1'$ and $\sigma_1''$ are assigned to the left ends of the monodromy graph in Figure \ref{fig:4-local-brief}(1) or (2), we have $4$ ways to choose permutation $\tau$ from Lemma \ref{lem:computation-3}(3).
Hence, we also obtain $d!$ tuples in $\fl^\rb_0((a_1,b_1),(a_2,b_2))$, 
which correspond to a real cover $\pi$ in the set $\rl_0((a_1,b_1),(a_2,b_2))$ with trivial automorphism group.
Note that we have two ways to label the ends in $C$ that correspond to the ramification $(a_1,b_1)$.
We have $H^\rb(\varphi,v)=2\alpha(\varphi,v)=3!\cdot2^{\frac{d-4}{2}}(\frac{d-6}{2})!$, if $a_2,b_2$ are even integers. If $a_2,b_2$ are odd integers,  $H^\rb(\varphi,v)=2\alpha(\varphi,v)=2^{\frac{d-2}{2}}(\frac{d-4}{2})!$.
Moreover, $|\rl_0((a_1,b_1),(a_2,b_2))|=1$.
\item $a_1\neq b_1$ and $a_2=b_2$. The same argument as case $(b)$ shows that there is one real cover $\pi$ in the set $\rl_0((a_1,b_1),(a_2,b_2))$, and its automorphism group is trivial.
We have two ways to label the ends in $C$ that correspond to the ramification $(a_2,b_2)$.
We obtain that $H^\rb(\varphi,v)=2\alpha(\varphi,v)=3!\cdot2^{\frac{d-4}{2}}(\frac{d-6}{2})!$ if $a_2,b_2$ are even. When $a_2,b_2$ are odd integers,
$H^\rb(\varphi,v)=2\alpha(\varphi,v)=2^{\frac{d-2}{2}}(\frac{d-4}{2})!$.
\end{enumerate}
Let $A=(\sigma_1,\gamma,\tau,\sigma_2)$ be a tuple in $\fl^\rb_0((a_1,b_1),(a_2,b_2))$, and $(\pi,\iota)=\Psi([A])$.
Suppose that $\tilde A=(\sigma_1,\gamma,\bar\tau_1,\bar\tau_2,\sigma_2)\in\widetilde\fl^\rb_0((a_1,b_1),(a_2,b_2))$
is the tuple corresponding to $A$ via Lemma \ref{lem:computation-3}(1), and $a_1>a_2$.
Let $(\tilde\pi:\widetilde\cl\to\widetilde \dl,\tilde\iota)$ be the real simple resolution
corresponding to $\tilde A$.
From Lemma \ref{lem:computation-3}(2), the composition $\bar\tau_1\circ\sigma_1$ consists of three disjoint cycles.
Suppose the lengths are $b_1,a_2,a_1-a_2$.
A similar argument as the proof of Lemma \ref{lem:1-local} shows that $\tilde\pi^{-1}(B)$ consists of three distinct circles $B_1$, $B_2$ and $B_3$, where $B$ is the circle in Figure \ref{fig:3-local2-res}.
Moreover, degrees of $\tilde\pi|_{B_1}$, $\tilde\pi|_{B_2}$ and $\tilde\pi|_{B_3}$ are $b_1$, $a_2$ and $a_1-a_2$.
Riemann-Hurwitz formula implies that $\widetilde\cl$ is divided by $B_1,B_2$ and $B_3$ into four genus zero components.
It follows from \cite[Lemma 5.3]{mr-2015} and \cite[Lemma 5.5]{mr-2015} that the real tropical cover corresponding to $\tilde\pi\in\widetilde\rl_0((a_1,b_1),(a_2,b_2))$ is depicted in Figure \ref{fig:3-local2-res} (resp. Figure \ref{fig:3-local2-res1}) if $a_2,b_2$ are even (resp. odd).

Now we consider the case that $a_1,b_1$ are odd integers.
We have $\frac{d!}{a_1\cdot b_1}$ (resp. $\frac{d!}{2a_1\cdot b_1}$) possibilities
to choose $\sigma_1$, if $a_1\neq b_1$ (resp. if $a_1=b_1$).
Denote by $\sigma_1'$ and $\sigma_1''$ the circles that  represent the cycles of lengths $a_1$ and $b_1$, respectively.
From \cite[Lemma $3.12$]{gpmr-2015} (see also Figure \ref{fig:sym-circle1}),
we have $a_1b_1$ ways to choose the involution $\gamma$.
\begin{enumerate}[$(i)$]
    \item $a_1\neq b_1$ and $a_2\neq b_2$. Without loss of generality, we assume $a_1>a_2$. 
If $a_2,b_2$ are even integers (resp. odd integers), 
$j(v)$ equals the number of ways to assign $a_1,a_2,b_1,b_2$ to the ends of the monodromy graph in Figure \ref{fig:4-local-brief}(3) (resp. in Figure \ref{fig:4-local-brief}(4)).
Once $a_1,b_1,a_2,b_2$ are assigned to ends of the monodromy graph in Figure \ref{fig:4-local-brief}$(3)$ or $(4)$, there is only one way to label $\sigma_1'$ and $\sigma_1''$ to the left ends of the monodromy graph. 
From Lemma \ref{lem:computation-3}(3), we have one way to choose $\tau$.
    Hence, one way to assign $a_1,b_1,a_2,b_2$ to the ends of the monodromy graph in Figure \ref{fig:4-local-brief}$(3)$ or $(4)$ corresponds to $d!$ tuples in $(\sigma_1,\gamma,\tau,\sigma_2)\in\fl^\rb((a_1,b_1),(a_2,b_2))$.
    These $d!$ tuples correspond to a real ramified covering $\pi$ in the set $\rl_0((a_1,b_1),(a_2,b_2))$ with trivial automorphism group.
    Moreover, from Lemma \ref{lem:fixed-points}(2), the set $\pi^{-1}(p)$
    consists of two real points and $\frac{d-4}{2}$ pairs of conjugated points (see Figure $\ref{fig:4-local}(3)$ or (4)). The image of $F$ is determined by Lemma \ref{lem:fixed-points}(1).
    Therefore, $H^\rb(\varphi,v)=j(v)\cdot\alpha(\varphi,v)=j(v)\cdot2^{\frac{d-4}{2}}(\frac{d-4}{2})!$,
and $|\rl_0((a_1,b_1),(a_2,b_2))|=j(v)$.
    \item $a_1=b_1$ and $a_2\neq b_2$. Without loss of generality, we assume $a_1>a_2$. We have two ways to assign $\sigma_1'$ and $\sigma_1''$ to the left ends
of the monodromy graph in Figure \ref{fig:4-local-brief}(3) or (4).
When $\sigma_1'$ and $\sigma_1''$ are assigned to the left ends of the monodromy graph in Figure \ref{fig:4-local-brief}(3) or (4), we have $1$ way to choose permutation $\tau$ from Lemma \ref{lem:computation-3}(3).
    Hence, we get $d!$ tuples in $(\sigma_1,\gamma,\tau,\sigma_2)\in\fl^\rb((a_1,b_1),(a_2,b_2))$,
    which correspond to a real cover $\pi$ in the set $\rl_0((a_1,b_1),(a_2,b_2))$ with trivial automorphism group.
    Moreover, from Lemma \ref{lem:fixed-points}(2), the set $\pi^{-1}(p)$
    consists of two real points and $\frac{d-4}{2}$ pairs of conjugated points (see Figure $\ref{fig:4-local}(3)$ or (4)).
    We have two ways to label the ends in $C$ that  correspond to the ramification $(a_1,b_1)$, so $H^\rb(\varphi,v)=2\alpha(\varphi,v)=2^{\frac{d-2}{2}}(\frac{d-4}{2})!$.
Moreover, $|\rl_0((a_1,b_1),(a_2,b_2))|=1$.
    \item $a_1\neq b_1$ and $a_2=b_2$. The same argument as case $(ii)$ shows that
    there is one real cover $\pi$ in the set $\rl_0((a_1,b_1),(a_2,b_2))$, and its automorphism group is trivial. The set $\pi^{-1}(p)$
    consists of two real points and $\frac{d-4}{2}$ pairs of conjugated points (see Figure $\ref{fig:4-local}(3)$ or (4)).
    Since we have two ways to label the ends in $C$ that  correspond to the ramification $(a_2,b_2)$,
    $H^\rb(\varphi,v)=2\alpha(\varphi,v)=2^{\frac{d-2}{2}}(\frac{d-4}{2})!$.
\end{enumerate}
The construction of real tropical covers corresponding to the real simple resolution is similar to the case that $a_1,a_2$ are even,
so we omit the proof and only depict them in Figure \ref{fig:3-local2-res2} and Figure \ref{fig:3-local2-res3}.
\end{proof}

\section{Figures}
\label{sec:b}
In this appendix we give some figures which are used in the above sections.

\begin{figure}[H]
    \centering
    \begin{tikzpicture}
    \draw (0,0)--(1.5,0)--(4,1);
    \draw (1.5,0)--(2.5,-0.3)--(4,-0.3);
    \draw (0,-0.8)--(2.5,-0.3);
    \draw (2.5,0.4)--(2.4,1);
    \draw (2.6,1)--(2.5,0.4)--(2.9,1);
    \draw (3.3,1)--(2.5,0.4)--(3.6,1);
    \draw (2.4,1.1) node{\tiny$1$} (2.6,1.1) node{\tiny$1$}
    (2.9,1.1) node{\tiny$1$} (3.3,1.1) node{\tiny$1$} (3.6,1.1) node{\tiny$1$};
    \draw[<->,bend left] (2.6,1.2) to (2.9,1.2);
    \draw[<->,bend left] (3.3,1.2) to (3.6,1.2);
    \draw[line width=0.9pt,dotted] (2.95,1)--(3.25,1);
    \draw (0,-1)--(1.5,-0.5);
    \draw (0.6,-1)--(1.5,-0.5)--(0.9,-1);
    \draw (1.5,-1)--(1.5,-0.5)--(1.8,-1);
    \draw[line width=0.9pt,dotted] (1,-1)--(1.4,-1);
    \draw[<->,bend right] (0.6,-1.2) to (0.9,-1.2);
    \draw[<->,bend right] (1.5,-1.2) to (1.8,-1.2);
    \draw (0.5,-0.4) node{\tiny$b_1$} (0,-1.1) node{\tiny$1$}
    (0.6,-1.1) node{\tiny$1$} (0.9,-1.1) node{\tiny$1$}
    (1.5,-1.1) node{\tiny$1$} (1.8,-1.1) node{\tiny$1$};
    \draw (0,1)--(1.5,0);
    \draw (0.6,1)--(1.5,0)--(0.9,1);
    \draw (1.5,1)--(1.5,0)--(1.8,1);
    \draw[line width=0.9pt,dotted] (1.2,0.5)--(1.5,0.5);
    \draw[<->,bend left] (0.6,1.2) to (0.9,1.2);
    \draw[<->,bend left] (1.5,1.2) to (1.8,1.2);
    \draw (2,-1)--(2.5,-0.3);
    \draw (2.6,-1)--(2.5,-0.3)--(2.9,-1);
    \draw (3.5,-1)--(2.5,-0.3)--(3.8,-1);
    \draw[line width=0.9pt,dotted] (2.95,-1)--(3.45,-1);
    \draw[<->,bend right] (2.6,-1.2) to (2.9,-1.2);
    \draw[<->,bend right] (3.5,-1.2) to (3.8,-1.2);
    \draw (0.5,0.15) node{\tiny$a_1$} (2,0.3) node{\tiny$b_2$} (2.4,-0.1) node{\tiny$a_1-b_2$} (3.8,0.8) node{\tiny$b_2$} (3.5,-0.2) node{\tiny$a_2$} (0,1.1) node{\tiny$2$}
    (0.6,1.1) node{\tiny$1$} (0.9,1.1) node{\tiny$1$}
    (1.5,1.1) node{\tiny$1$} (1.8,1.1) node{\tiny$1$} (2,-1.1) node{\tiny$2$}
    (2.6,-1.1) node{\tiny$1$} (2.9,-1.1) node{\tiny$1$}
    (3.5,-1.1) node{\tiny$1$} (3.8,-1.1) node{\tiny$1$};
    \draw[bend right] (6,0.5) to(6,-0.5);
    \draw[bend left,dotted] (6,0.5) to(6,-0.5);
    \draw[bend right] (8,0.5) to(8,-0.5);
    \draw[bend left,dotted] (8,0.5) to(8,-0.5);
    \draw (7,-1.1)  ellipse(0.4 and 0.2);
    \draw[bend right] (6,0.5) to (8,0.5);
    \draw (6,-0.5)..controls(6.5,-0.5) and (6.6,-0.5)..(6.6,-1.1);
    \draw (7.4,-1.1)..controls(7.4,-0.5) and (7.5,-0.5)..(8,-0.5);
    \draw (5.6,0) node{\tiny$(a_1,b_1)$} (8.3,0) node{\tiny$B$}
    (10.4,0) node{\tiny$(a_2,b_2)$} (7.1,-1.4) node{\tiny$(2,1,\ldots,1)$} (9.1,-1.4) node{\tiny$(2,1,\ldots,1)$};
    \draw (10,0) ellipse (0.2 and 0.5);
    \draw (9,-1.1)  ellipse(0.4 and 0.2);
    \draw[bend right] (8,0.5) to (10,0.5);
    \draw (8,-0.5)..controls(8.5,-0.5) and (8.6,-0.5)..(8.6,-1.1);
    \draw (9.4,-1.1)..controls(9.4,-0.5) and (9.5,-0.5)..(10,-0.5);
    \foreach \Point in {(6,0.5),(10,0.5)}
    \draw[fill=red] \Point circle (0.05);
    \draw[fill=red] (7.4,-1.1) circle (0.05);
    \draw[fill=red] (8.6,-1.1) circle (0.05);
    \end{tikzpicture}
    \caption{Real tropical cover corresponding to $\tilde\pi_1$: the map $F$ used to define real tropical cover is encoded by marking the points in $\im(F)$ in red. $d$ is odd, and $a_1, a_2$ are even.}
    \label{fig:3-local1-res}
\end{figure}

\begin{figure}[H]
\begin{tikzpicture}
    \draw (0,0)--(1.5,0)--(4,1);
    \draw (1.5,0)--(2.5,-0.3)--(4,-0.3);
    \draw (0,-0.8)--(2.5,-0.3);
    \draw (2.2,1)--(2.5,0.4)--(2.4,1);
    \draw (2.6,1)--(2.5,0.4)--(2.9,1);
    \draw (3.3,1)--(2.5,0.4)--(3.6,1);
    \draw (2.2,1.1) node{\tiny$1$} (2.4,1.1) node{\tiny$1$} (2.6,1.1) node{\tiny$1$}
    (2.9,1.1) node{\tiny$1$} (3.3,1.1) node{\tiny$1$} (3.6,1.1) node{\tiny$1$};
    \draw[<->,bend left] (2.6,1.2) to (2.9,1.2);
    \draw[<->,bend left] (3.3,1.2) to (3.6,1.2);
    \draw[line width=0.9pt,dotted] (2.95,1)--(3.25,1);
    \draw (0,-1)--(1.5,-0.5);
    \draw (0.6,-1)--(1.5,-0.5)--(0.9,-1);
    \draw (1.5,-1)--(1.5,-0.5)--(1.8,-1);
    \draw[line width=0.9pt,dotted] (1,-1)--(1.4,-1);
    \draw[<->,bend right] (0.6,-1.2) to (0.9,-1.2);
    \draw[<->,bend right] (1.5,-1.2) to (1.8,-1.2);
    \draw (0.5,-0.45) node{\tiny$b_1$} (0,-1.1) node{\tiny$1$}
    (0.6,-1.1) node{\tiny$1$} (0.9,-1.1) node{\tiny$1$}
    (1.5,-1.1) node{\tiny$1$} (1.8,-1.1) node{\tiny$1$};
    \draw (0,1)--(1.5,0);
    \draw (0.2,1)--(1.5,0)--(0.4,1);
    \draw (0.6,1)--(1.5,0)--(0.9,1);
    \draw (1.5,1)--(1.5,0)--(1.8,1);
    \draw[line width=0.9pt,dotted] (1.2,0.5)--(1.5,0.5);
    \draw[<->,bend left] (0.6,1.2) to (0.9,1.2);
    \draw[<->,bend left] (1.5,1.2) to (1.8,1.2);
    \draw (2,-1)--(2.5,-0.3)--(2.3,-1);
    \draw (2.6,-1)--(2.5,-0.3)--(2.9,-1);
    \draw (3.5,-1)--(2.5,-0.3)--(3.8,-1);
    \draw[line width=0.9pt,dotted] (2.95,-1)--(3.45,-1);
    \draw[<->,bend right] (2.6,-1.2) to (2.9,-1.2);
    \draw[<->,bend right] (3.5,-1.2) to (3.8,-1.2);
    \draw (0.5,0.15) node{\tiny$a_1$} (2,0.3) node{\tiny$a_2$} (2.4,-0.1) node{\tiny$a_1-a_2$} (3.8,0.8) node{\tiny$a_2$} (3.5,-0.2) node{\tiny$b_2$} (0,1.1) node{\tiny$2$}
    (0.2,1.1) node{\tiny$1$} (0.4,1.1) node{\tiny$1$}
    (0.6,1.1) node{\tiny$1$} (0.9,1.1) node{\tiny$1$}
    (1.5,1.1) node{\tiny$1$} (1.8,1.1) node{\tiny$1$} (2,-1.1) node{\tiny$2$}
    (2.3,-1.1) node{\tiny$1$}
    (2.6,-1.1) node{\tiny$1$} (2.9,-1.1) node{\tiny$1$}
    (3.5,-1.1) node{\tiny$1$} (3.8,-1.1) node{\tiny$1$};
    \draw[bend right] (6,0.5) to(6,-0.5);
    \draw[bend left,dotted] (6,0.5) to(6,-0.5);
    \draw[bend right] (8,0.5) to(8,-0.5);
    \draw[bend left,dotted] (8,0.5) to(8,-0.5);
    \draw (7,-1.1)  ellipse(0.4 and 0.2);
    \draw[bend right] (6,0.5) to (8,0.5);
    \draw (6,-0.5)..controls(6.5,-0.5) and (6.6,-0.5)..(6.6,-1.1);
    \draw (7.4,-1.1)..controls(7.4,-0.5) and (7.5,-0.5)..(8,-0.5);
    \draw (5.6,0) node{\tiny$(a_1,b_1)$} (8.3,0) node{\tiny$B$}
    (10.4,0) node{\tiny$(a_2,b_2)$} (7.1,-1.4) node{\tiny$(2,1,\ldots,1)$} (9.1,-1.4) node{\tiny$(2,1,\ldots,1)$}
    (8,-0.65) node{\tiny\color{red}$2\times$};
    \draw (10,0) ellipse (0.2 and 0.5);
    \draw (9,-1.1)  ellipse(0.4 and 0.2);
    \draw[bend right] (8,0.5) to (10,0.5);
    \draw (8,-0.5)..controls(8.5,-0.5) and (8.6,-0.5)..(8.6,-1.1);
    \draw (9.4,-1.1)..controls(9.4,-0.5) and (9.5,-0.5)..(10,-0.5);
    \foreach \Point in {(8,-0.5),(6,-0.5),(10,-0.5),(7.4,-1.1),(8.6,-1.1)}
    \draw[fill=red] \Point circle (0.05);
\end{tikzpicture}
    \caption{Real tropical cover corresponding to $\tilde\pi_2$: the map $F$ used to define real tropical cover is encoded by marking the points in $\im(F)$ in red. $d$ is odd, and $a_1,a_2$ are even.}
    \label{fig:3-local1-res1}
\end{figure}

\begin{figure}[H]
\begin{tikzpicture}
    \draw (0,0)--(1.5,0)--(4,1);
    \draw (1.5,0)--(2.5,-0.3)--(4,-0.3);
    \draw (0,-0.8)--(2.5,-0.3);
    \draw (2.2,1)--(2.5,0.4)--(2.4,1);
    \draw (2.6,1)--(2.5,0.4)--(2.9,1);
    \draw (3.3,1)--(2.5,0.4)--(3.6,1);
    \draw (2.2,1.1) node{\tiny$1$} (2.4,1.1) node{\tiny$1$} (2.6,1.1) node{\tiny$1$}
    (2.9,1.1) node{\tiny$1$} (3.3,1.1) node{\tiny$1$} (3.6,1.1) node{\tiny$1$};
    \draw[<->,bend left] (2.6,1.2) to (2.9,1.2);
    \draw[<->,bend left] (3.3,1.2) to (3.6,1.2);
    \draw[line width=0.9pt,dotted] (2.95,1)--(3.25,1);
    \draw (0,-1)--(1.5,-0.5)--(0.3,-1);
    \draw (0.6,-1)--(1.5,-0.5)--(0.9,-1);
    \draw (1.5,-1)--(1.5,-0.5)--(1.8,-1);
    \draw[line width=0.9pt,dotted] (1,-1)--(1.4,-1);
    \draw[<->,bend right] (0.6,-1.2) to (0.9,-1.2);
    \draw[<->,bend right] (1.5,-1.2) to (1.8,-1.2);
    \draw (0.5,-0.45) node{\tiny$a_1$} (0,-1.1) node{\tiny$1$}
    (0.3,-1.1) node{\tiny$1$}
    (0.6,-1.1) node{\tiny$1$} (0.9,-1.1) node{\tiny$1$}
    (1.5,-1.1) node{\tiny$1$} (1.8,-1.1) node{\tiny$1$};
    \draw (0,1)--(1.5,0);
    \draw (0.2,1)--(1.5,0);
    \draw (0.6,1)--(1.5,0)--(0.9,1);
    \draw (1.5,1)--(1.5,0)--(1.8,1);
    \draw[line width=0.9pt,dotted] (1.2,0.5)--(1.5,0.5);
    \draw[<->,bend left] (0.6,1.2) to (0.9,1.2);
    \draw[<->,bend left] (1.5,1.2) to (1.8,1.2);
    \draw (2,-1)--(2.5,-0.3)--(2.3,-1);
    \draw (2.6,-1)--(2.5,-0.3)--(2.9,-1);
    \draw (3.5,-1)--(2.5,-0.3)--(3.8,-1);
    \draw[line width=0.9pt,dotted] (2.95,-1)--(3.45,-1);
    \draw[<->,bend right] (2.6,-1.2) to (2.9,-1.2);
    \draw[<->,bend right] (3.5,-1.2) to (3.8,-1.2);
    \draw (0.5,0.15) node{\tiny$b_1$} (2,0.3) node{\tiny$a_2$} (2.4,-0.1) node{\tiny$b_1-a_2$} (3.8,0.8) node{\tiny$a_2$} (3.5,-0.2) node{\tiny$b_2$} (0,1.1) node{\tiny$2$}
    (0.2,1.1) node{\tiny$1$}
    (0.6,1.1) node{\tiny$1$} (0.9,1.1) node{\tiny$1$}
    (1.5,1.1) node{\tiny$1$} (1.8,1.1) node{\tiny$1$} (2,-1.1) node{\tiny$2$}
    (2.3,-1.1) node{\tiny$1$}
    (2.6,-1.1) node{\tiny$1$} (2.9,-1.1) node{\tiny$1$}
    (3.5,-1.1) node{\tiny$1$} (3.8,-1.1) node{\tiny$1$};
    \draw[bend right] (6,0.5) to(6,-0.5);
    \draw[bend left,dotted] (6,0.5) to(6,-0.5);
    \draw[bend right] (8,0.5) to(8,-0.5);
    \draw[bend left,dotted] (8,0.5) to(8,-0.5);
    \draw (7,-1.1)  ellipse(0.4 and 0.2);
    \draw[bend right] (6,0.5) to (8,0.5);
    \draw (6,-0.5)..controls(6.5,-0.5) and (6.6,-0.5)..(6.6,-1.1);
    \draw (7.4,-1.1)..controls(7.4,-0.5) and (7.5,-0.5)..(8,-0.5);
    \draw (5.6,0) node{\tiny$(a_1,b_1)$} (8.3,0) node{\tiny$B$}
    (10.4,0) node{\tiny$(a_2,b_2)$} (7.1,-1.4) node{\tiny$(2,1,\ldots,1)$} (9.1,-1.4) node{\tiny$(2,1,\ldots,1)$}
    (8,-0.65) node{\tiny\color{red}$2\times$};
    \draw (10,0) ellipse (0.2 and 0.5);
    \draw (9,-1.1)  ellipse(0.4 and 0.2);
    \draw[bend right] (8,0.5) to (10,0.5);
    \draw (8,-0.5)..controls(8.5,-0.5) and (8.6,-0.5)..(8.6,-1.1);
    \draw (9.4,-1.1)..controls(9.4,-0.5) and (9.5,-0.5)..(10,-0.5);
    \foreach \Point in {(8,-0.5),(6,-0.5),(10,-0.5),(7.4,-1.1),(8.6,-1.1)}
    \draw[fill=red] \Point circle (0.05);
\end{tikzpicture}
    \caption{Real tropical cover corresponding to $\tilde\pi_3$: the map $F$ used to define real tropical cover is encoded by marking the points in $\im(F)$ in red. $d$ is odd and $a_1,a_2$ are even.}
    \label{fig:3-local1-res2}
\end{figure}

\begin{figure}[H]
\begin{tikzpicture}
    \draw (0,0)--(1.5,0)--(4,1);
    \draw (1.5,0)--(2.5,-0.3)--(4,-0.3);
    \draw (0,-0.8)--(2.5,-0.3);
    \draw (2.2,1)--(2.5,0.4);
    \draw (2.6,1)--(2.5,0.4)--(2.9,1);
    \draw (3.3,1)--(2.5,0.4)--(3.6,1);
    \draw (2.2,1.1) node{\tiny$1$} 
    (2.6,1.1) node{\tiny$1$}
    (2.9,1.1) node{\tiny$1$} (3.3,1.1) node{\tiny$1$} (3.6,1.1) node{\tiny$1$};
    \draw[<->,bend left] (2.6,1.2) to (2.9,1.2);
    \draw[<->,bend left] (3.3,1.2) to (3.6,1.2);
    \draw[line width=0.9pt,dotted] (2.95,1)--(3.25,1);
    \draw (0,-1)--(1.5,-0.5)--(0.3,-1);
    \draw (0.6,-1)--(1.5,-0.5)--(0.9,-1);
    \draw (1.5,-1)--(1.5,-0.5)--(1.8,-1);
    \draw[line width=0.9pt,dotted] (1,-1)--(1.4,-1);
    \draw[<->,bend right] (0.6,-1.2) to (0.9,-1.2);
    \draw[<->,bend right] (1.5,-1.2) to (1.8,-1.2);
    \draw (0.5,-0.45) node{\tiny$a_1$} (0,-1.1) node{\tiny$1$}
    (0.3,-1.1) node{\tiny$1$}
    (0.6,-1.1) node{\tiny$1$} (0.9,-1.1) node{\tiny$1$}
    (1.5,-1.1) node{\tiny$1$} (1.8,-1.1) node{\tiny$1$};
    \draw (0,1)--(1.5,0);
    \draw (0.2,1)--(1.5,0);
    \draw (0.6,1)--(1.5,0)--(0.9,1);
    \draw (1.5,1)--(1.5,0)--(1.8,1);
    \draw[line width=0.9pt,dotted] (1.2,0.5)--(1.5,0.5);
    \draw[<->,bend left] (0.6,1.2) to (0.9,1.2);
    \draw[<->,bend left] (1.5,1.2) to (1.8,1.2);
    \draw (2.5,-0.3)--(2.15,-1);
    \draw (2,-1)--(2.5,-0.3)--(2.3,-1);
    \draw (2.6,-1)--(2.5,-0.3)--(2.9,-1);
    \draw (3.5,-1)--(2.5,-0.3)--(3.8,-1);
    \draw[line width=0.9pt,dotted] (2.95,-1)--(3.45,-1);
    \draw[<->,bend right] (2.6,-1.2) to (2.9,-1.2);
    \draw[<->,bend right] (3.5,-1.2) to (3.8,-1.2);
    \draw (0.5,0.15) node{\tiny$b_1$} (2,0.3) node{\tiny$b_2$} (2.4,-0.1) node{\tiny$b_1-b_2$} (3.8,0.8) node{\tiny$b_2$} (3.5,-0.2) node{\tiny$a_2$} (0,1.1) node{\tiny$2$}
    (0.2,1.1) node{\tiny$1$}
    (0.6,1.1) node{\tiny$1$} (0.9,1.1) node{\tiny$1$}
    (1.5,1.1) node{\tiny$1$} (1.8,1.1) node{\tiny$1$} (2,-1.1) node{\tiny$2$}
    (2.3,-1.1) node{\tiny$1$} (2.15,-1.1) node{\tiny$1$}
    (2.6,-1.1) node{\tiny$1$} (2.9,-1.1) node{\tiny$1$}
    (3.5,-1.1) node{\tiny$1$} (3.8,-1.1) node{\tiny$1$};
    \draw[bend right] (6,0.5) to(6,-0.5);
    \draw[bend left,dotted] (6,0.5) to(6,-0.5);
    \draw[bend right] (8,0.5) to(8,-0.5);
    \draw[bend left,dotted] (8,0.5) to(8,-0.5);
    \draw (7,-1.1)  ellipse(0.4 and 0.2);
    \draw[bend right] (6,0.5) to (8,0.5);
    \draw (6,-0.5)..controls(6.5,-0.5) and (6.6,-0.5)..(6.6,-1.1);
    \draw (7.4,-1.1)..controls(7.4,-0.5) and (7.5,-0.5)..(8,-0.5);
    \draw (5.6,0) node{\tiny$(a_1,b_1)$} (8.3,0) node{\tiny$B$}
    (10.4,0) node{\tiny$(a_2,b_2)$} (7.1,-1.4) node{\tiny$(2,1,\ldots,1)$} (9.1,-1.4) node{\tiny$(2,1,\ldots,1)$}
    (8,-0.65) node{\tiny\color{red}$2\times$};
    \draw (10,0) ellipse (0.2 and 0.5);
    \draw (9,-1.1)  ellipse(0.4 and 0.2);
    \draw[bend right] (8,0.5) to (10,0.5);
    \draw (8,-0.5)..controls(8.5,-0.5) and (8.6,-0.5)..(8.6,-1.1);
    \draw (9.4,-1.1)..controls(9.4,-0.5) and (9.5,-0.5)..(10,-0.5);
    \foreach \Point in {(8,-0.5),(6,-0.5),(10,-0.5),(7.4,-1.1),(8.6,-1.1)}
    \draw[fill=red] \Point circle (0.05);
\end{tikzpicture}
    \caption{Real tropical cover corresponding to $\tilde\pi_4$: the map $F$ used to define real tropical cover is encoded by marking the points in $\im(F)$ in red. $d$ is odd and $a_1, a_2$ are even.}
    \label{fig:3-local1-res3}
\end{figure}

\begin{figure}[H]
    \centering
    \begin{tikzpicture}
    \draw (0,0)--(1.5,0)--(4,1);
    \draw (1.5,0)--(2.5,-0.3)--(4,-0.3);
    \draw (0,-0.8)--(2.5,-0.3);
    \draw (2.2,1)--(2.5,0.4)--(2.4,1);
    \draw (2.6,1)--(2.5,0.4)--(2.9,1);
    \draw (3.3,1)--(2.5,0.4)--(3.6,1);
    \draw (2.2,1.1) node{\tiny$1$} (2.4,1.1) node{\tiny$1$} (2.6,1.1) node{\tiny$1$}
    (2.9,1.1) node{\tiny$1$} (3.3,1.1) node{\tiny$1$} (3.6,1.1) node{\tiny$1$};
    \draw[<->,bend left] (2.6,1.2) to (2.9,1.2);
    \draw[<->,bend left] (3.3,1.2) to (3.6,1.2);
    \draw[line width=0.9pt,dotted] (2.95,1)--(3.25,1);
    \draw (0,-1)--(1.5,-0.5);
    \draw (0.2,-1)--(1.5,-0.5);
    \draw (0.6,-1)--(1.5,-0.5)--(0.9,-1);
    \draw (1.5,-1)--(1.5,-0.5)--(1.8,-1);
    \draw[line width=0.9pt,dotted] (1,-1)--(1.4,-1);
    \draw[<->,bend right] (0.6,-1.2) to (0.9,-1.2);
    \draw[<->,bend right] (1.5,-1.2) to (1.8,-1.2);
    \draw (0.5,-0.5) node{\tiny$b_1$} (0,-1.1) node{\tiny$1$}
    (0.2,-1.1) node{\tiny$1$}
    (0.6,-1.1) node{\tiny$1$} (0.9,-1.1) node{\tiny$1$}
    (1.5,-1.1) node{\tiny$1$} (1.8,-1.1) node{\tiny$1$};
    \draw (0,1)--(1.5,0);
    \draw (0.4,1)--(1.5,0)--(0.2,1);
    \draw (0.6,1)--(1.5,0)--(0.9,1);
    \draw (1.5,1)--(1.5,0)--(1.8,1);
    \draw[line width=0.9pt,dotted] (1.2,0.5)--(1.5,0.5);
    \draw[<->,bend left] (0.6,1.2) to (0.9,1.2);
    \draw[<->,bend left] (1.5,1.2) to (1.8,1.2);
    \draw (2,-1)--(2.5,-0.3);
    \draw (2.2,-1)--(2.5,-0.3)--(2.4,-1);
    \draw (2.6,-1)--(2.5,-0.3)--(2.9,-1);
    \draw (3.5,-1)--(2.5,-0.3)--(3.8,-1);
    \draw[line width=0.9pt,dotted] (2.95,-1)--(3.45,-1);
    \draw[<->,bend right] (2.6,-1.2) to (2.9,-1.2);
    \draw[<->,bend right] (3.5,-1.2) to (3.8,-1.2);
    \draw (0.5,0.15) node{\tiny$a_1$} (2,0.3) node{\tiny$a_2$} (2.4,-0.1) node{\tiny$a_1-a_2$} (3.8,0.8) node{\tiny$a_2$} (3.5,-0.2) node{\tiny$b_2$} (0,1.1) node{\tiny$2$}
    (0.2,-1.1) node{\tiny$1$} (0.4,1.1) node{\tiny$1$} (0.2,1.1) node{\tiny$1$}
    (0.6,1.1) node{\tiny$1$} (0.9,1.1) node{\tiny$1$}
    (1.5,1.1) node{\tiny$1$} (1.8,1.1) node{\tiny$1$} (2,-1.1) node{\tiny$2$}
    (2.2,-1.1) node{\tiny$1$} (2.4,-1.1) node{\tiny$1$}
    (2.6,-1.1) node{\tiny$1$} (2.9,-1.1) node{\tiny$1$}
    (3.5,-1.1) node{\tiny$1$} (3.8,-1.1) node{\tiny$1$};
    \draw[bend right] (6,0.5) to(6,-0.5);
    \draw[bend left,dotted] (6,0.5) to(6,-0.5);
    \draw[bend right] (8,0.5) to(8,-0.5);
    \draw[bend left,dotted] (8,0.5) to(8,-0.5);
    \draw (7,-1.1)  ellipse(0.4 and 0.2);
    \draw[bend right] (6,0.5) to (8,0.5);
    \draw (6,-0.5)..controls(6.5,-0.5) and (6.6,-0.5)..(6.6,-1.1);
    \draw (7.4,-1.1)..controls(7.4,-0.5) and (7.5,-0.5)..(8,-0.5);
    \draw (5.6,0) node{\tiny$(a_1,b_1)$} (6,-0.7) node{\tiny\color{red}$2\times$} (8,-0.7) node{\tiny\color{red}$3\times$} (10,-0.7) node{\tiny\color{red}$2\times$} (8.3,0) node{\tiny$B$}
    (10.4,0) node{\tiny$(a_2,b_2)$} (7.1,-1.4) node{\tiny$(2,1,\ldots,1)$} (9.1,-1.4) node{\tiny$(2,1,\ldots,1)$};
    \draw (10,0) ellipse (0.2 and 0.5);
    \draw (9,-1.1)  ellipse(0.4 and 0.2);
    \draw[bend right] (8,0.5) to (10,0.5);
    \draw (8,-0.5)..controls(8.5,-0.5) and (8.6,-0.5)..(8.6,-1.1);
    \draw (9.4,-1.1)..controls(9.4,-0.5) and (9.5,-0.5)..(10,-0.5);
    \foreach \Point in {(8,-0.5),(6,-0.5),(10,-0.5)}
    \draw[fill=red] \Point circle (0.05);
    \draw[fill=red] (7.4,-1.1) circle (0.05);
    \draw[fill=red] (8.6,-1.1) circle (0.05);
    \end{tikzpicture}
    \caption{$a_1,a_2,b_1,b_2$ are even: the map $F$ used to define real tropical cover is encoded by marking the points in $\im(F)$ in red.}
    \label{fig:3-local2-res}
\end{figure}

\begin{figure}[H]
\begin{tikzpicture}
    \draw (0,0)--(1.5,0)--(4,1);
    \draw (1.5,0)--(2.5,-0.3)--(4,-0.3);
    \draw (0,-0.8)--(2.5,-0.3);
    \draw (2.2,1)--(2.5,0.4);
    \draw (2.6,1)--(2.5,0.4)--(2.9,1);
    \draw (3.3,1)--(2.5,0.4)--(3.6,1);
    \draw (2.2,1.1) node{\tiny$1$} 
    (2.6,1.1) node{\tiny$1$}
    (2.9,1.1) node{\tiny$1$} (3.3,1.1) node{\tiny$1$} (3.6,1.1) node{\tiny$1$};
    \draw[<->,bend left] (2.6,1.2) to (2.9,1.2);
    \draw[<->,bend left] (3.3,1.2) to (3.6,1.2);
    \draw[line width=0.9pt,dotted] (2.95,1)--(3.25,1);
    \draw (0,-1)--(1.5,-0.5);
    \draw (0.2,-1)--(1.5,-0.5);
    \draw (0.6,-1)--(1.5,-0.5)--(0.9,-1);
    \draw (1.5,-1)--(1.5,-0.5)--(1.8,-1);
    \draw[line width=0.9pt,dotted] (1,-1)--(1.4,-1);
    \draw[<->,bend right] (0.6,-1.2) to (0.9,-1.2);
    \draw[<->,bend right] (1.5,-1.2) to (1.8,-1.2);
    \draw (0.5,-0.5) node{\tiny$b_1$} (0,-1.1) node{\tiny$1$}
    (0.2,-1.1) node{\tiny$1$}
    (0.6,-1.1) node{\tiny$1$} (0.9,-1.1) node{\tiny$1$}
    (1.5,-1.1) node{\tiny$1$} (1.8,-1.1) node{\tiny$1$};
    \draw (0,1)--(1.5,0);
    \draw (0.6,1)--(1.5,0)--(0.9,1);
    \draw (1.5,1)--(1.5,0)--(1.8,1);
    \draw[line width=0.9pt,dotted] (1.2,0.5)--(1.5,0.5);
    \draw[<->,bend left] (0.6,1.2) to (0.9,1.2);
    \draw[<->,bend left] (1.5,1.2) to (1.8,1.2);
    \draw (2,-1)--(2.5,-0.3);
    \draw (2.2,-1)--(2.5,-0.3);
    \draw (2.6,-1)--(2.5,-0.3)--(2.9,-1);
    \draw (3.5,-1)--(2.5,-0.3)--(3.8,-1);
    \draw[line width=0.9pt,dotted] (2.95,-1)--(3.45,-1);
    \draw[<->,bend right] (2.6,-1.2) to (2.9,-1.2);
    \draw[<->,bend right] (3.5,-1.2) to (3.8,-1.2);
    \draw (0.5,0.15) node{\tiny$a_1$} (2,0.3) node{\tiny$a_2$} (2.4,-0.1) node{\tiny$a_1-a_2$} (3.8,0.8) node{\tiny$a_2$} (3.5,-0.2) node{\tiny$b_2$} (0,1.1) node{\tiny$2$}
    (0.2,-1.1) node{\tiny$1$}
    (0.6,1.1) node{\tiny$1$} (0.9,1.1) node{\tiny$1$}
    (1.5,1.1) node{\tiny$1$} (1.8,1.1) node{\tiny$1$} (2,-1.1) node{\tiny$2$}
    (2.2,-1.1) node{\tiny$1$}
    (2.6,-1.1) node{\tiny$1$} (2.9,-1.1) node{\tiny$1$}
    (3.5,-1.1) node{\tiny$1$} (3.8,-1.1) node{\tiny$1$};
    \draw[bend right] (6,0.5) to(6,-0.5);
    \draw[bend left,dotted] (6,0.5) to(6,-0.5);
    \draw[bend right] (8,0.5) to(8,-0.5);
    \draw[bend left,dotted] (8,0.5) to(8,-0.5);
    \draw (7,-1.1)  ellipse(0.4 and 0.2);
    \draw[bend right] (6,0.5) to (8,0.5);
    \draw (6,-0.5)..controls(6.5,-0.5) and (6.6,-0.5)..(6.6,-1.1);
    \draw (7.4,-1.1)..controls(7.4,-0.5) and (7.5,-0.5)..(8,-0.5);
    \draw (5.6,0) node{\tiny$(a_1,b_1)$}
    (6.6,0.7) node{\tiny\color{red}$F(e)$ with $\omega(e)=a_1$}
    (8.3,0) node{\tiny$B$}
    (10.4,0) node{\tiny$(a_2,b_2)$} (7.1,-1.4) node{\tiny$(2,1,\ldots,1)$} (9.1,-1.4) node{\tiny$(2,1,\ldots,1)$};
    \draw (10,0) ellipse (0.2 and 0.5);
    \draw (9,-1.1)  ellipse(0.4 and 0.2);
    \draw[bend right] (8,0.5) to (10,0.5);
    \draw (8,-0.5)..controls(8.5,-0.5) and (8.6,-0.5)..(8.6,-1.1);
    \draw (9.4,-1.1)..controls(9.4,-0.5) and (9.5,-0.5)..(10,-0.5);
    \foreach \Point in {(8,-0.5),(6,-0.5),(6,0.5),(7.4,-1.1),(8.6,-1.1)}
    \draw[fill=red] \Point circle (0.05);
\end{tikzpicture}
    \caption{$a_1,b_1$ are even, and $a_2,b_2$ are odd: the map $F$ used to define real tropical cover is encoded by marking the points in $\im(F)$ in red.}
    \label{fig:3-local2-res1}
\end{figure}

\begin{figure}[H]
\begin{tikzpicture}
    \draw (0,0)--(1.5,0)--(4,1);
    \draw (1.5,0)--(2.5,-0.3)--(4,-0.3);
    \draw (0,-0.8)--(2.5,-0.3);
    \draw (2.2,1)--(2.5,0.4)--(2.4,1);
    \draw (2.6,1)--(2.5,0.4)--(2.9,1);
    \draw (3.3,1)--(2.5,0.4)--(3.6,1);
    \draw (2.2,1.1) node{\tiny$1$} (2.4,1.1) node{\tiny$1$} (2.6,1.1) node{\tiny$1$}
    (2.9,1.1) node{\tiny$1$} (3.3,1.1) node{\tiny$1$} (3.6,1.1) node{\tiny$1$};
    \draw[<->,bend left] (2.6,1.2) to (2.9,1.2);
    \draw[<->,bend left] (3.3,1.2) to (3.6,1.2);
    \draw[line width=0.9pt,dotted] (2.95,1)--(3.25,1);
    \draw (0,-1)--(1.5,-0.5);
    \draw (0.6,-1)--(1.5,-0.5)--(0.9,-1);
    \draw (1.5,-1)--(1.5,-0.5)--(1.8,-1);
    \draw[line width=0.9pt,dotted] (1,-1)--(1.4,-1);
    \draw[<->,bend right] (0.6,-1.2) to (0.9,-1.2);
    \draw[<->,bend right] (1.5,-1.2) to (1.8,-1.2);
    \draw (0.5,-0.45) node{\tiny$b_1$} (0,-1.1) node{\tiny$1$}
    (0.6,-1.1) node{\tiny$1$} (0.9,-1.1) node{\tiny$1$}
    (1.5,-1.1) node{\tiny$1$} (1.8,-1.1) node{\tiny$1$};
    \draw (0,1)--(1.5,0);
    \draw (0.2,1)--(1.5,0);
    \draw (0.6,1)--(1.5,0)--(0.9,1);
    \draw (1.5,1)--(1.5,0)--(1.8,1);
    \draw[line width=0.9pt,dotted] (1.2,0.5)--(1.5,0.5);
    \draw[<->,bend left] (0.6,1.2) to (0.9,1.2);
    \draw[<->,bend left] (1.5,1.2) to (1.8,1.2);
    \draw (2,-1)--(2.5,-0.3);
    \draw (2.6,-1)--(2.5,-0.3)--(2.9,-1);
    \draw (3.5,-1)--(2.5,-0.3)--(3.8,-1);
    \draw[line width=0.9pt,dotted] (2.95,-1)--(3.45,-1);
    \draw[<->,bend right] (2.6,-1.2) to (2.9,-1.2);
    \draw[<->,bend right] (3.5,-1.2) to (3.8,-1.2);
    \draw (0.5,0.15) node{\tiny$a_1$} (2,0.3) node{\tiny$a_2$} (2.4,-0.1) node{\tiny$a_1-a_2$} (3.8,0.8) node{\tiny$a_2$} (3.5,-0.2) node{\tiny$b_2$} (0,1.1) node{\tiny$2$}
    (0.2,1.1) node{\tiny$1$}
    (0.6,1.1) node{\tiny$1$} (0.9,1.1) node{\tiny$1$}
    (1.5,1.1) node{\tiny$1$} (1.8,1.1) node{\tiny$1$} (2,-1.1) node{\tiny$2$}
    (2.6,-1.1) node{\tiny$1$} (2.9,-1.1) node{\tiny$1$}
    (3.5,-1.1) node{\tiny$1$} (3.8,-1.1) node{\tiny$1$};
    \draw[bend right] (6,0.5) to(6,-0.5);
    \draw[bend left,dotted] (6,0.5) to(6,-0.5);
    \draw[bend right] (8,0.5) to(8,-0.5);
    \draw[bend left,dotted] (8,0.5) to(8,-0.5);
    \draw (7,-1.1)  ellipse(0.4 and 0.2);
    \draw[bend right] (6,0.5) to (8,0.5);
    \draw (6,-0.5)..controls(6.5,-0.5) and (6.6,-0.5)..(6.6,-1.1);
    \draw (7.4,-1.1)..controls(7.4,-0.5) and (7.5,-0.5)..(8,-0.5);
    \draw (5.6,0) node{\tiny$(a_1,b_1)$} (8.3,0) node{\tiny$B$}
    (10.4,0) node{\tiny$(a_2,b_2)$} (7.1,-1.4) node{\tiny$(2,1,\ldots,1)$} (9.1,-1.4) node{\tiny$(2,1,\ldots,1)$}
    (10.6,0.7) node{\tiny\color{red}$F(e)$ with $\omega(e)=b_2$};
    \draw (10,0) ellipse (0.2 and 0.5);
    \draw (9,-1.1)  ellipse(0.4 and 0.2);
    \draw[bend right] (8,0.5) to (10,0.5);
    \draw (8,-0.5)..controls(8.5,-0.5) and (8.6,-0.5)..(8.6,-1.1);
    \draw (9.4,-1.1)..controls(9.4,-0.5) and (9.5,-0.5)..(10,-0.5);
    \foreach \Point in {(8,-0.5),(10,0.5),(10,-0.5),(7.4,-1.1),(8.6,-1.1)}
    \draw[fill=red] \Point circle (0.05);
\end{tikzpicture}
    \caption{$a_1,b_1$ are odd, and $a_2,b_2$ are even: the map $F$ used to define real tropical cover is encoded by marking the points in $\im(F)$ in red.}
    \label{fig:3-local2-res2}
\end{figure}

\begin{figure}[H]
\begin{tikzpicture}
    \draw (0,0)--(1.5,0)--(4,1);
    \draw (1.5,0)--(2.5,-0.3)--(4,-0.3);
    \draw (0,-0.8)--(2.5,-0.3);
    \draw (2.2,1)--(2.5,0.4);
    \draw (2.6,1)--(2.5,0.4)--(2.9,1);
    \draw (3.3,1)--(2.5,0.4)--(3.6,1);
    \draw (2.2,1.1) node{\tiny$1$} 
    (2.6,1.1) node{\tiny$1$}
    (2.9,1.1) node{\tiny$1$} (3.3,1.1) node{\tiny$1$} (3.6,1.1) node{\tiny$1$};
    \draw[<->,bend left] (2.6,1.2) to (2.9,1.2);
    \draw[<->,bend left] (3.3,1.2) to (3.6,1.2);
    \draw[line width=0.9pt,dotted] (2.95,1)--(3.25,1);
    \draw (0,-1)--(1.5,-0.5);
    \draw (0.6,-1)--(1.5,-0.5)--(0.9,-1);
    \draw (1.5,-1)--(1.5,-0.5)--(1.8,-1);
    \draw[line width=0.9pt,dotted] (1,-1)--(1.4,-1);
    \draw[<->,bend right] (0.6,-1.2) to (0.9,-1.2);
    \draw[<->,bend right] (1.5,-1.2) to (1.8,-1.2);
    \draw (0.5,-0.5) node{\tiny$b_1$} (0,-1.1) node{\tiny$1$}
    (0.6,-1.1) node{\tiny$1$} (0.9,-1.1) node{\tiny$1$}
    (1.5,-1.1) node{\tiny$1$} (1.8,-1.1) node{\tiny$1$};
    \draw (0,1)--(1.5,0);
    \draw (0.2,1)--(1.5,0);
    \draw (0.6,1)--(1.5,0)--(0.9,1);
    \draw (1.5,1)--(1.5,0)--(1.8,1);
    \draw[line width=0.9pt,dotted] (1.2,0.5)--(1.5,0.5);
    \draw[<->,bend left] (0.6,1.2) to (0.9,1.2);
    \draw[<->,bend left] (1.5,1.2) to (1.8,1.2);
    \draw (2,-1)--(2.5,-0.3)--(2.3,-1);
    \draw (2.6,-1)--(2.5,-0.3)--(2.9,-1);
    \draw (3.5,-1)--(2.5,-0.3)--(3.8,-1);
    \draw[line width=0.9pt,dotted] (2.95,-1)--(3.45,-1);
    \draw[<->,bend right] (2.6,-1.2) to (2.9,-1.2);
    \draw[<->,bend right] (3.5,-1.2) to (3.8,-1.2);
    \draw (0.5,0.15) node{\tiny$a_1$} (2,0.3) node{\tiny$a_2$} (2.4,-0.1) node{\tiny$a_1-a_2$} (3.8,0.8) node{\tiny$a_2$} (3.5,-0.2) node{\tiny$b_2$} (0,1.1) node{\tiny$2$}
    (0.2,1.1) node{\tiny$1$}
    (0.6,1.1) node{\tiny$1$} (0.9,1.1) node{\tiny$1$}
    (1.5,1.1) node{\tiny$1$} (1.8,1.1) node{\tiny$1$} (2,-1.1) node{\tiny$2$}
    (2.3,-1.1) node{\tiny$1$}
    (2.6,-1.1) node{\tiny$1$} (2.9,-1.1) node{\tiny$1$}
    (3.5,-1.1) node{\tiny$1$} (3.8,-1.1) node{\tiny$1$};
    \draw[bend right] (6,0.5) to(6,-0.5);
    \draw[bend left,dotted] (6,0.5) to(6,-0.5);
    \draw[bend right] (8,0.5) to(8,-0.5);
    \draw[bend left,dotted] (8,0.5) to(8,-0.5);
    \draw (7,-1.1)  ellipse(0.4 and 0.2);
    \draw[bend right] (6,0.5) to (8,0.5);
    \draw (6,-0.5)..controls(6.5,-0.5) and (6.6,-0.5)..(6.6,-1.1);
    \draw (7.4,-1.1)..controls(7.4,-0.5) and (7.5,-0.5)..(8,-0.5);
    \draw (5.6,0) node{\tiny$(a_1,b_1)$} (8.3,0) node{\tiny$B$}
    (10.4,0) node{\tiny$(a_2,b_2)$} (7.1,-1.4) node{\tiny$(2,1,\ldots,1)$} (9.1,-1.4) node{\tiny$(2,1,\ldots,1)$};
    \draw (10,0) ellipse (0.2 and 0.5);
    \draw (9,-1.1)  ellipse(0.4 and 0.2);
    \draw[bend right] (8,0.5) to (10,0.5);
    \draw (8,-0.5)..controls(8.5,-0.5) and (8.6,-0.5)..(8.6,-1.1);
    \draw (9.4,-1.1)..controls(9.4,-0.5) and (9.5,-0.5)..(10,-0.5);
    \foreach \Point in {(8,-0.5),(7.4,-1.1),(8.6,-1.1)}
    \draw[fill=red] \Point circle (0.05);
\end{tikzpicture}
    \caption{$a_1,b_1,a_2,b_2$ are odd: the map $F$ used to define real tropical cover is encoded by marking the points in $\im(F)$ in red.}
    \label{fig:3-local2-res3}
\end{figure}

\begin{figure}[H]
    \centering
    \begin{tikzpicture}
    \draw (0,0.5)--(2,0)--(4,0.5);
    \draw (0.7,-1)--(2,0);
    \draw (1,-1)--(2,0)--(1.4,-1);
    \draw (2.6,-1)--(2,0)--(3,-1);
    \draw[line width=0.9pt,dotted] (1.8,-0.5)--(2.2,-0.5);
    \draw[<->,bend right] (0.9,-1.2) to (1.4,-1.2);
    \draw[<->,bend right] (2.6,-1.2) to (3.1,-1.2);
    \draw (0.5,0.25) node{\tiny$d$} (3.3,0.45) node{\tiny$d$}
    (1,-1.1) node{\tiny$1$} 
    (0.7,-1.1) node{\tiny$1$} (1.4,-1.1) node{\tiny$1$} (2.6,-1.1) node{\tiny$1$} (3,-1.1) node{\tiny$1$};
    \foreach \Point in {(6,0),(8,0)}
    \draw \Point ellipse (0.2 and 0.5);
    \draw (7,-1.1)  ellipse(0.4 and 0.2);
    \draw[bend right] (6,0.5) to (8,0.5);
    \draw (6,-0.5)..controls(6.5,-0.5) and (6.6,-0.5)..(6.6,-1.1);
    \draw (7.4,-1.1)..controls(7.4,-0.5) and (7.5,-0.5)..(8,-0.5);
    \draw (5.2,0) node{\tiny$(d)$} (8.7,0) node{\tiny$(d)$}
    (8.1,-1.2) node{\tiny$(1,\ldots,1)$} (4.6,-1.6) node{$(1)$ $d$ is odd.};
    \draw (0,-2)--(2,-2.5)--(4,-2);
    \draw (1,-3.5)--(2,-2.5)--(1.4,-3.5);
    \draw (2.6,-3.5)--(2,-2.5)--(3,-3.5);
    \draw[line width=0.9pt,dotted] (1.8,-3)--(2.2,-3);
    \draw[<->,bend right] (0.9,-3.7) to (1.4,-3.7);
    \draw[<->,bend right] (2.6,-3.7) to (3.1,-3.7);
    \draw (0.5,-2.25) node{\tiny$d$} (3.3,-2.05) node{\tiny$d$}
    (1,-3.6) node{\tiny$1$} (1.4,-3.6) node{\tiny$1$} (2.6,-3.6) node{\tiny$1$} (3,-3.6) node{\tiny$1$};
    \foreach \Point in {(6,-2.5),(8,-2.5)}
    \draw \Point ellipse (0.2 and 0.5);
    \draw (7,-3.6)  ellipse(0.4 and 0.2);
    \draw[bend right] (6,-2) to (8,-2);
    \draw (6,-3)..controls(6.5,-3) and (6.6,-3)..(6.6,-3.6);
    \draw (7.4,-3.6)..controls(7.4,-3) and (7.5,-3)..(8,-3);
    \foreach \Point in {(8,-2),(6,-2)}
    \draw[fill=red] \Point circle (0.05);
    \draw (5.2,-2.5) node{\tiny$(d)$} (8.7,-2.5) node{\tiny$(d)$} (8.1,-3.7) node{\tiny$(1,\ldots,1)$} (4.6,-4.1) node{$(2)$ $d$ is even, and two fixed points for even fixed edges are connected by a real arc.};
    \draw (0,-5)--(2,-5.5)--(4,-5);
    \draw (2,-5.5)--(0.7,-6.5);
    \draw (2,-5.5)--(0.2,-6.5);
    \draw (1,-6.5)--(2,-5.5)--(1.4,-6.5);
    \draw (2.6,-6.5)--(2,-5.5)--(3,-6.5);
    \draw[line width=0.9pt,dotted] (1.8,-6)--(2.2,-6);
    \draw[<->,bend right] (0.9,-6.7) to (1.4,-6.7);
    \draw[<->,bend right] (2.6,-6.7) to (3.1,-6.7);
    \draw (0.5,-5.25) node{\tiny$d$} (3.3,-5.05) node{\tiny$d$}
    (0.2,-6.6) node{\tiny$1$} (0.7,-6.6) node{\tiny$1$} (1,-6.6) node{\tiny$1$} (1.4,-6.6) node{\tiny$1$} (2.6,-6.6) node{\tiny$1$} (3,-6.6) node{\tiny$1$};
    \foreach \Point in {(6,-5.5),(8,-5.5)}
    \draw \Point ellipse (0.2 and 0.5);
    \draw (7,-6.6)  ellipse(0.4 and 0.2);
    \draw[bend right] (6,-5) to (8,-5);
    \draw (6,-6)..controls(6.5,-6) and (6.6,-6)..(6.6,-6.6);
    \draw (7.4,-6.6)..controls(7.4,-6) and (7.5,-6)..(8,-6);
    \foreach \Point in {(8,-6),(6,-6)}
    \draw[fill=red] \Point circle (0.05);
    \draw (5.2,-5.5) node{\tiny$(d)$} (8.7,-5.5) node{\tiny$(d)$} (8.1,-6.7) node{\tiny$(1,1,\ldots,1)$} (4.6,-7.1) node{$(3)$ $d$ is even, and two fixed points for even fixed edges can not be connected by a real arc.};
    \end{tikzpicture}
    \caption{The points in $\im(F)$ are coloured in red.}
    \label{fig:local-simple}
\end{figure}

\begin{figure}[H]
\begin{tikzpicture}
    \draw (0,0.5)--(1.5,0)--(2.5,0)--(4,0.5);
    \draw (0.5,-1)--(1.5,0);
    \draw (0.8,-1)--(1.5,0)--(1.1,-1);
    \draw (1.5,-1)--(1.5,0)--(1.8,-1);
    \draw (2,-1)--(2.5,0);
    \draw (2.3,-1)--(2.5,0)--(2.6,-1);
    \draw (3.5,-1)--(2.5,0)--(3.8,-1);
    \draw (0.5,-1.1) node{\tiny$1$} (0.8,-1.1) node{\tiny$1$} (1.1,-1.1) node{\tiny$1$}
    (1.5,-1.1) node{\tiny$1$} (1.8,-1.1) node{\tiny$1$}(2,-1.1) node{\tiny$1$} (2.3,-1.1) node{\tiny$1$} (2.6,-1.1) node{\tiny$1$}
    (3.5,-1.1) node{\tiny$1$} (3.8,-1.1) node{\tiny$1$} (0.1,0.7) node{\tiny$d$} (2,0.2)node{\tiny$d$} (3.9,0.7) node{\tiny$d$};
    \draw[<->,bend right] (0.8,-1.2) to (1.1,-1.2);
    \draw[<->,bend right] (1.5,-1.2) to (1.8,-1.2);
    \draw[line width=0.9pt,dotted] (1.2,-0.8)--(1.45,-0.8);
    \draw[<->,bend right] (2.3,-1.2) to (2.6,-1.2);
    \draw[<->,bend right] (3.5,-1.2) to (3.8,-1.2);
    \draw[line width=0.9pt,dotted] (2.7,-0.8)--(3.2,-0.8);
    \draw[bend right] (6,0.5) to(6,-0.5);
    \draw[bend left,dotted] (6,0.5) to(6,-0.5);
    \draw[bend right] (8,0.5) to(8,-0.5);
    \draw[bend left,dotted] (8,0.5) to(8,-0.5);
    \draw (7,-1.1)  ellipse(0.4 and 0.2);
    \draw[bend right] (6,0.5) to (8,0.5);
    \draw (6,-0.5)..controls(6.5,-0.5) and (6.6,-0.5)..(6.6,-1.1);
    \draw (7.4,-1.1)..controls(7.4,-0.5) and (7.5,-0.5)..(8,-0.5);
    \draw (5.6,0) node{\tiny$(d)$} (8.3,0) node{\tiny$B$}
    (10.4,0) node{\tiny$(d)$} (7.1,-1.4) node{\tiny$(1,\ldots,1)$} (9.1,-1.4) node{\tiny$(1,\ldots,1)$};
    \draw (10,0) ellipse (0.2 and 0.5);
    \draw (9,-1.1)  ellipse(0.4 and 0.2);
    \draw[bend right] (8,0.5) to (10,0.5);
    \draw (8,-0.5)..controls(8.5,-0.5) and (8.6,-0.5)..(8.6,-1.1);
    \draw (9.4,-1.1)..controls(9.4,-0.5) and (9.5,-0.5)..(10,-0.5);
\end{tikzpicture}
    \caption{$d$ is odd.}
    \label{fig:local-simple1}
\end{figure}

\begin{figure}[H]
\begin{tikzpicture}
    \draw (0,0.5)--(1.5,0)--(2.5,0)--(4,0.5);
    \draw (0.8,-1)--(1.5,0)--(1.1,-1);
    \draw (1.5,-1)--(1.5,0)--(1.8,-1);
    \draw (2.3,-1)--(2.5,0)--(2.6,-1);
    \draw (3.5,-1)--(2.5,0)--(3.8,-1);
    \draw (0.8,-1.1) node{\tiny$1$} (1.1,-1.1) node{\tiny$1$}
    (1.5,-1.1) node{\tiny$1$} (1.8,-1.1) node{\tiny$1$} (2.3,-1.1) node{\tiny$1$} (2.6,-1.1) node{\tiny$1$}
    (3.5,-1.1) node{\tiny$1$} (3.8,-1.1) node{\tiny$1$} (0.1,0.7) node{\tiny$d$} (2,0.2)node{\tiny$d$} (3.9,0.7) node{\tiny$d$};
    \draw[<->,bend right] (0.8,-1.2) to (1.1,-1.2);
    \draw[<->,bend right] (1.5,-1.2) to (1.8,-1.2);
    \draw[line width=0.9pt,dotted] (1.2,-0.8)--(1.45,-0.8);
    \draw[<->,bend right] (2.3,-1.2) to (2.6,-1.2);
    \draw[<->,bend right] (3.5,-1.2) to (3.8,-1.2);
    \draw[line width=0.9pt,dotted] (2.7,-0.8)--(3.2,-0.8);
    \draw[bend right] (6,0.5) to(6,-0.5);
    \draw[bend left,dotted] (6,0.5) to(6,-0.5);
    \draw[bend right] (8,0.5) to(8,-0.5);
    \draw[bend left,dotted] (8,0.5) to(8,-0.5);
    \draw (7,-1.1)  ellipse(0.4 and 0.2);
    \draw[bend right] (6,0.5) to (8,0.5);
    \draw (6,-0.5)..controls(6.5,-0.5) and (6.6,-0.5)..(6.6,-1.1);
    \draw (7.4,-1.1)..controls(7.4,-0.5) and (7.5,-0.5)..(8,-0.5);
    \draw (5.6,0) node{\tiny$(d)$} (8.3,0) node{\tiny$B$}
    (10.4,0) node{\tiny$(d)$} (7.1,-1.4) node{\tiny$(1,\ldots,1)$} (9.1,-1.4) node{\tiny$(1,\ldots,1)$};
    \draw (10,0) ellipse (0.2 and 0.5);
    \draw (9,-1.1)  ellipse(0.4 and 0.2);
    \draw[bend right] (8,0.5) to (10,0.5);
    \draw (8,-0.5)..controls(8.5,-0.5) and (8.6,-0.5)..(8.6,-1.1);
    \draw (9.4,-1.1)..controls(9.4,-0.5) and (9.5,-0.5)..(10,-0.5);
    \foreach \Point in {(8,0.5),(6,0.5),(10,0.5)}
    \draw[fill=red] \Point circle (0.05);
\end{tikzpicture}
    \caption{The points in $\im(F)$ are in red. $d$ is even, and fixed points for even fixed edges are connected by a real arc.}
    \label{fig:local-simple2}
\end{figure}

\begin{figure}[H]
\begin{tikzpicture}
     \draw (0,0.5)--(1.5,0)--(2.5,0)--(4,0.5);
    \draw (0.5,-1)--(1.5,0)--(0.3,-1);
    \draw (0.8,-1)--(1.5,0)--(1.1,-1);
    \draw (1.5,-1)--(1.5,0)--(1.8,-1);
    \draw (2,-1)--(2.5,0)--(2.15,-1);
    \draw (2.3,-1)--(2.5,0)--(2.6,-1);
    \draw (3.5,-1)--(2.5,0)--(3.8,-1);
    \draw (0.5,-1.1) node{\tiny$1$} (0.3,-1.1) node{\tiny$1$} (0.8,-1.1) node{\tiny$1$} (1.1,-1.1) node{\tiny$1$}
    (1.5,-1.1) node{\tiny$1$} (1.8,-1.1) node{\tiny$1$}(2,-1.1) node{\tiny$1$} (2.15,-1.1) node{\tiny$1$} (2.3,-1.1) node{\tiny$1$} (2.6,-1.1) node{\tiny$1$}
    (3.5,-1.1) node{\tiny$1$} (3.8,-1.1) node{\tiny$1$} (0.1,0.7) node{\tiny$d$} (2,0.2)node{\tiny$d$} (3.9,0.7) node{\tiny$d$};
    \draw[<->,bend right] (0.8,-1.2) to (1.1,-1.2);
    \draw[<->,bend right] (1.5,-1.2) to (1.8,-1.2);
    \draw[line width=0.9pt,dotted] (1.2,-0.8)--(1.45,-0.8);
    \draw[<->,bend right] (2.3,-1.2) to (2.6,-1.2);
    \draw[<->,bend right] (3.5,-1.2) to (3.8,-1.2);
    \draw[line width=0.9pt,dotted] (2.7,-0.8)--(3.2,-0.8);
    \draw[bend right] (6,0.5) to(6,-0.5);
    \draw[bend left,dotted] (6,0.5) to(6,-0.5);
    \draw[bend right] (8,0.5) to(8,-0.5);
    \draw[bend left,dotted] (8,0.5) to(8,-0.5);
    \draw (7,-1.1)  ellipse(0.4 and 0.2);
    \draw[bend right] (6,0.5) to (8,0.5);
    \draw (6,-0.5)..controls(6.5,-0.5) and (6.6,-0.5)..(6.6,-1.1);
    \draw (7.4,-1.1)..controls(7.4,-0.5) and (7.5,-0.5)..(8,-0.5);
    \draw (5.6,0) node{\tiny$(d)$} (8.3,0) node{\tiny$B$}
    (10.4,0) node{\tiny$(d)$} (7.1,-1.4) node{\tiny$(1,\ldots,1)$} (9.1,-1.4) node{\tiny$(1,\ldots,1)$};
    \draw (10,0) ellipse (0.2 and 0.5);
    \draw (9,-1.1)  ellipse(0.4 and 0.2);
    \draw[bend right] (8,0.5) to (10,0.5);
    \draw (8,-0.5)..controls(8.5,-0.5) and (8.6,-0.5)..(8.6,-1.1);
    \draw (9.4,-1.1)..controls(9.4,-0.5) and (9.5,-0.5)..(10,-0.5);
    \foreach \Point in {(8,-0.5),(6,-0.5),(10,-0.5)}
    \draw[fill=red] \Point circle (0.05);
\end{tikzpicture}
    \caption{The points in $\im(F)$ are in red. $d$ is even, and fixed points for even fixed edges can not be connected by a real arc.}
    \label{fig:local-simple3}
\end{figure}

\section*{Statements and Declarations}
The authors declare that they have no known competing financial interests or personal
relationships that could have appeared to influence the work reported in this paper.

\section*{Acknowledgements}
The first author is grateful to Ilia Itenberg for
pointing out the problem of asymptotic behavior of real double Hurwitz numbers
when non-simple branch points are added, during his visit to IMJ-PRG in 2020.
The authors thank Di Yang for suggesting the problem of uniform asymptotics of double Hurwitz numbers,
as well as for valuable discussions and suggestions.
The authors would like to thank Chenglang Yang and Zhiyuan Wang for their helpful discussions,
and they also thank Jianfeng Wu for carefully reading the manuscript.
Y. Ding was supported by the National Natural Science Foundation of China (No.12101565), and
the Natural Science Foundation of Henan (No. 212300410287).
H. Liu was supported by the National Natural Science Foundation of China (No.12171439).


\end{document}